\documentclass[10pt]{article}
\usepackage{amsmath}
\usepackage{amssymb}
\usepackage{amsthm}
\usepackage{graphicx}
\usepackage{subfigure}
\usepackage{verbatim}
\usepackage[square]{natbib}

\font\elevenss=cmss11

\font\eightss=cmss8

\font\sixss=cmss8 at 6pt

\newfam\ssfam
\textfont\ssfam=\elevenss \scriptfont\ssfam=\eightss
 \scriptscriptfont\ssfam=\sixss
\def\ss{\fam\ssfam \elevenss}%

\theoremstyle{plain}
\newtheorem{thm}{Theorem}[section]
\newtheorem{lem}[thm]{Lemma}
\newtheorem{pr}[thm]{Proposition}
\newtheorem{cor}[thm]{Corollary}
\newtheorem{defn}[thm]{Definition}
\newtheorem{hyps}[thm]{Hypotheses}

\newtheorem{example}[thm]{Example}
\theoremstyle{remark}

\newtheorem*{unremark}{Remark}

\newtheorem*{unremarks}{Remarks}

\newcommand{\Em}[1]{\textbf{#1}}

\def\grad{\nabla}
\def\ee{\epsilon}

\def\vv{{\bf v}}
\def\uu{{\bf u}}
\def\xx{{\bf x}}
\def\yy{{\bf y}}
\def\zz{{\bf z}}

\def\ZZ{{\bf Z}}
\def\pp{{\bf p}}
\def\mm{{\bf m}}
\def\kk{{\bf k}}
\def\ww{{\bf w}}
\def\rr{{\bf r}}
\def\xmax{{\xx_{\min}}}
\def\bb{{\bf b}}
\def\nn{{\bf n}}
\def\half{{\frac{1}{2}}}
\def\flattorus{T_\R}
\def\tor3{{\mathbb T}^3}
\def\Proj{{\mathbb P}}
\def\cK{{\bf K}}
\def\cKd{{\bf K^*}}
\def\cL{{\bf L}}
\def\cLd{{\bf L^*}}
\def\CC{{\cal C}}
\def\dchain{{\CC_\delta}}
\def\proj{{\overline{\CC}}}
\def\generic{{\CC}}

\def\sing{{\cal V}}

\def\domain{{\cal D}}
\def\one{{\bf 1}}
\def\zero{{\bf 0}}
\def\rhat{\hat{\rr}}

\def\rt{\hat{r}}
\def\st{\hat{s}}
\def\tat{\hat{t}}
\def\nbd{{\cal N}}

\def\Cox{\hfill \Box}
\def\C{{\mathbb C}}
\def\Z{{\mathbb Z}}

\def\R{{\mathbb R}}
\def\CP{{\mathbb C}{\mathbb P}}
\def\RP{{\mathbb R}{\mathbb P}}
\def\Real{{\rm Re\,}}

\def\Arg{{\rm Arg\,}}
\def\Log{{\rm ReLog\,}}
\def\Res{{\rm Res}}
\def\dblres{{\Res}^{(2)}}
\def\thom{{\delta}}

\def\M{{\cal M}}

\def\form{A}
\def\dual{\form^*}
\def\lorentz{S}
\def\lorentzd{\lorentz^*}
\def\amoeba{{\mbox{\ss amoeba}}}
\def\At{\tilde{A}}
\def\Rdd{(\R^d)^*}
\def\qt{{\tilde{q}}}
\def\hht{{\tilde{h}}}
\def\contrib{{\mbox{\ss contrib}}}
\def\op{{\cal S}}
\def\fth{{\hat h}}
\def\normal{{\bf N^*}}
\def\normalx{{\bf N}_\xx^*}
\def\tangent{\mbox{\ss tan}}
\def\critical{{\mbox{\ss crit}}}
\def\logcrit{{\mbox{\ss W}}}

\def\loggrad{{\grad_{\log}}}
\def\ft{{\overline{f}}}
\def\gauss{\kappa}
\def\vhat{{\hat{\bf v}}}

\def\fiberchain{- i \uu + \R^d}
\def\wchain{\CC_\ww}
\def\newton{{\mbox{\ss P}}}
\def\hull{\mbox{\ss hull}\,}
\def\ctorus{(\C^*)^d}
\def\joy{{\mathbf \xi}}

\def\homog{\mbox{\ss hom}}
\def\symp{S}
\def\disp{\displaystyle}
\def\cs{C_0 (\R^{d*})}
\def\rd{C_{\rm RD} (\R^d)}
\def\locint{\mbox{\ss loc-int}}
\def\polybd{\mbox{\ss poly-bd}} 
\def\gf{{\cal L}}

\def\gfspace{{\cal G}^*}
\def\gfs{{\cal G}}
\def\Sup{\mbox{\ss supp}\,}
\def\IFT{{\cal F}^{-1}}

\def\imf{\iota f}
\def\term{u}
\def\ellipse{{\cal E}}

\makeatletter \def\romenumi{ \def\theenumi{\roman{enumi}}
\def\p@enumi{\theenumi} \def\labelenumi{(\@roman\c@enumi)}} 
\makeatother

\begin{document}
\renewcommand{\thepage}{\roman{page}}

\begin{titlepage}
\begin{center}
{\large \bf Asymptotics of multivariate sequences, part III:
quadratic points} \\
\end{center}
\vspace{5ex}
\begin{flushright}
Yuliy Baryshnikov \footnote{Bell Laboratories, Lucent Technologies,
700 Mountain Avenue, Murray Hill, NJ 07974-0636, ymb@research.bell-labs.com}
  ~\\
Robin Pemantle \footnote{Research supported in part by
National Science Foundation grant \# DMS 0603821}$^,$\footnote{University
of Pennsylvania, Department of Mathematics, 209 S. 33rd Street, Philadelphia,
PA 19104 USA, pemantle@math.upenn.edu}
\end{flushright}

\vfill

\noindent{\bf ABSTRACT:} We consider a number of combinatorial
problems in which rational generating functions may be obtained,
whose denominators have factors with certain singularities.
Specifically, there exist points near which one of
the factors is asymptotic to a nondegenerate quadratic. 
We compute the asymptotics of the coefficients of such a 
generating function.  The computation requires some topological 
deformations as well as Fourier-Laplace transforms of generalized 
functions.  We apply the results of the theory to specific 
combinatorial problems, such as Aztec diamond tilings, cube
groves, and multi-set permutations.
\hfill \\[1ex]

\vfill



\end{titlepage}

\tableofcontents

%
{
\begin{tabular}{lll}
&& \\[-0.5in]
& & {\bf Glossary of notation} \\[2ex]
page & symbol & meaning \\[2ex]
\pageref{relog} & $\Log$ & coordinatewise log-modulus \\
\pageref{deg} & $\deg (f , \zz)$ & degree of vanishing of $f$ at $\zz$ \\
\pageref{homog} & $\homog (f , \zz)$ & leading homogeneous part of $f$ at $\zz$ \\
\pageref{sing} & $\sing , \sing_F$ & variety where $F$ vanishes \\
\pageref{amoeba} & $\amoeba (F)$ & amoeba of the Laurent polynomial $F$ \\
\pageref{flattorus} & $\flattorus$ & the torus $\R^d / \Z^d$  \\
\pageref{dual} & $\L^*$ & dual cone to a cone $\cL$ \\
\pageref{tan} & $\tangent_\xx (C)$ & tangent cone to $C$ at $\xx$  \\
\pageref{normalx} & $\normalx (C)$ & dual cone to $-\tan_{\xx^*} (C)$ \\
\pageref{xmax} & $\xmax$ & the point of a given region where $\rr \cdot \xx$ is maximized \\
\pageref{hypcone} & $\cK^\vv (A)$ & cone of hyperbolicity in direction $\vv$ of a homogeneous polynomial $A$ \\
\pageref{hypcone2} & $\cK^{A,B}(\xx)$ & cone of hyperbolicity of $A$ at $\xx$ containing $B$ \\
\pageref{hypcone} & $\cK^{f,B} (\ZZ)$ & cone of hyperbolicity in direction $\vv$ of a polynomial $f$ at $B$ \\
\pageref{ft} & $\ft$ & abbreviation for $\homog (f , \xx + i \yy)$ \\
\pageref{normal} & $\normal (\ZZ) , (\normal)^{f , B} (\ZZ)$ & the normal cone to $f$ at $\ZZ$ that contains $B$ \\
\pageref{sing1} & $\sing_1$ & intersection of $\sing$ with the torus whose image under $\Log$ is $\xmax$ \\
\pageref{crit} & $\critical$ & set of minimal critical points in direction $\rr$ \\
\pageref{logcrit} & $\logcrit, \logcrit (\rr)$ & logarithmic version of $\critical$ \\
\pageref{loggrad} & $\loggrad$ & the logarithmic gradient \\
\pageref{oexp} & $o_{\exp}$ & less by an exponential factor \\
\pageref{lorentz} & $\lorentz$ & the standard Lorentzian quadratic \\
\pageref{form} & $\dual$ & dual to the quadratic form $\form$ \\
\pageref{contrib} & $\contrib$ & contribution to the Cauchy integral from the chain $\CC_\ww$ local to $\ww$ \\
\pageref{gauss} & $\gauss$ & Gaussian curvature \\
\pageref{U} & $U_\ww$ & a neighborhood of $\ww \in \logcrit (\rr)$ \\
\pageref{U} & $U$ & the union of all the $U_\ww$ \\
\pageref{eta} & $\eta_{U^c}$ & an outward vector field on $U^c$ that is a section of the cones $\cK^{f,b} (\cdot)$ \\
\pageref{eta2} & $\eta$ & a section of $\cK^{f,b} (\cdot)$ defined everywhere but pointing outward only on $U^c$ \\
\pageref{Phi} & $\Phi , \Phi^{\ee , \eta}$ & $\ee$-scaled homotopy from a constant inward vector field to $C(\eta)$ \\
\pageref{Ceta} & $\generic (\eta)$ & the cycle resulting from sliding along $\eta$ \\
\pageref{Cw} & $\generic_\ww$ & restriction of $\generic (\eta)$ to a neighborhood of $\ww$ \\
\pageref{proj} & $\proj$ & projective chain \\
\pageref{projd} & $\proj^{(\delta)}$ & projective chain lifted off $\sing$ in the $\delta$-ball \\
\pageref{dchainw} & $\dchain (\ww)$ & $\generic (\eta)$ for a particular $\eta$, restricted to a neighborhood of $\ww$ \\
\pageref{dchain} & $\dchain$ & chain pieced together from local chains $\dchain (\ww)$ \\
\pageref{cs} & $\cs$ & the space of test functions \\
\pageref{gfspace} & $\gfspace$ & the space of generalized functions \\
\pageref{locint} & $\locint$ & the space of locally integrable functions \\
\pageref{rd} & $\rd$ & the space of rapidly decaying functions \\
\pageref{IFT} & $\IFT$ & inverse Fourier transform \\
\pageref{leray} & $\alpha$ & the Leray cycle \\
\pageref{petrovsky} & $\gamma$ & the Petrovsky cycle \\
\pageref{Res} & $\Res [\theta , \overline{\sing_H}]$ & the form $\omega / dH$ \\
\pageref{dblres},~\pageref{dblres 2} & $\dblres$ & the second residue of $\omega / (\qt \hht)$ on the double pole $\sing_{\qt} \cap \sing_{\hht}$ \\
\end{tabular}
}

\renewcommand{\thepage}{\arabic{page}}
\setcounter{page}{1}
\setcounter{equation}{0}
\section{Introduction} \label{sec:intro} 

\subsection{Background and motivation}

Problems in combinatorial enumeration and discrete probability
can often be attacked by means of generating functions.  If one
is lucky enough to obtain a closed form generating function,
then the asymptotic enumeration formula, or probabilistic
limit theorem is often not far behind.  Recently, several 
problems have arisen to which can be associated very nice
generating functions, in fact rational functions of several 
variables, but for which asymptotic estimates have not followed
(although formulae were found in some cases by other means).  
These problems include random tilings (the so-called 
Aztec and Diabolo tilings) and other statistical
mechanical ensembles (cube groves) as well as some 
enumerative and graph theoretic problems discussed
later in the paper.  

A series of recent papers~\cite{PW1,PW2,BP-residues} provides 
a method for asymptotic evaluation of the coefficients of 
multivariate generating functions.  To describe the scope of
this previous work, we set up some notation that will be
in force for the rest of this article.  Throughout,
we will assume that the generating function converges in a domain
defining there a {\it quasirational} function
\begin{equation} \label{eq:F form}
F(\ZZ) = \frac{P(\ZZ)}{Q(\ZZ)^s \prod_{j=1}^k H_j (\ZZ)^{n_j}} 
   = \sum_{\rr} a_\rr \ZZ^\rr
\end{equation}
with polynomial $P, Q$, affine-linear $H_j$'s, integer $n_j$'s and
real $s$. (Here the quantities
in boldface are vectors of dimension $d$ and the notation $\ZZ^\rr$ 
is used to denote $\prod_{j=1}^d Z_j^{r_j}$.) 
In dimension three and below, we use $X,Y,Z$ to denote $Z_1 , Z_2$ 
and $Z_3$ respectively.  We let $\sing := \{ \ZZ : Q(\ZZ) = 0 \} 
\subseteq \C^d$ denote the \Em{pole variety} of $F$, that is, 
the complex algebraic variety where $Q$ vanishes ($Q$ will always be
a polynomial).  Analytic methods for recovering asymptotics of 
$a_\rr$ from $F$ always begin with the multivariate Cauchy integral 
formula
\begin{equation} \label{eq:mv cauchy}
a_\rr = \frac{1}{(2 \pi i)^{d}} \int_{T} \ZZ^{-\rr} F \,  
   \frac{d\ZZ}{\ZZ} \, .
\end{equation}
Here $T$ is a $d$-torus, i.e. the product of circles about the
origin in each coordinate axis (importantly, choice of a torus affects
the corresponding Laurent expansion (\ref{eq:F form})).  The pole set $\sing$ is of 
central importance because the contour $T$ may be deformed 
without affecting the integral as long as one avoids places 
where the integrand is singular.  

When $\sing$ is smooth or has singularities of self-intersection
type (where locally $\sing$ is the union of smooth divisors),
a substantial amount is known.  The case where $\sing$ is
smooth is analyzed in~\cite{PW1}; the existence under further
hypotheses of a local central limit theorem dates back
at least as far as~\cite{BR2}.  The more general 
case where the singular points of $\sing$
are all unions of smooth components with normal
intersections is analyzed using explicit
changes of variables in~\cite{PW2} and by multivariate
residues in~\cite{BP-residues}, the pre-cursor to 
which is~\cite{bertozzi-mckenna}.  Applications in which $\sing$
satisfies these conditions are abundant, and 
a number of examples are worked in~\cite{PW9}.  
For bivariate generating functions, all rational generating 
functions we have seen fall within this class.  
Other local geometries are possible, namely those of irreducible
monomial curve, 
e.g., $X^p - Y^q = 0$, but, as will become clear, they cannot contribute
to the asymptotic expansions, being non-hyperbolic.

\begin{figure}[ht]
\hspace{1.2in} \includegraphics[scale=0.40]{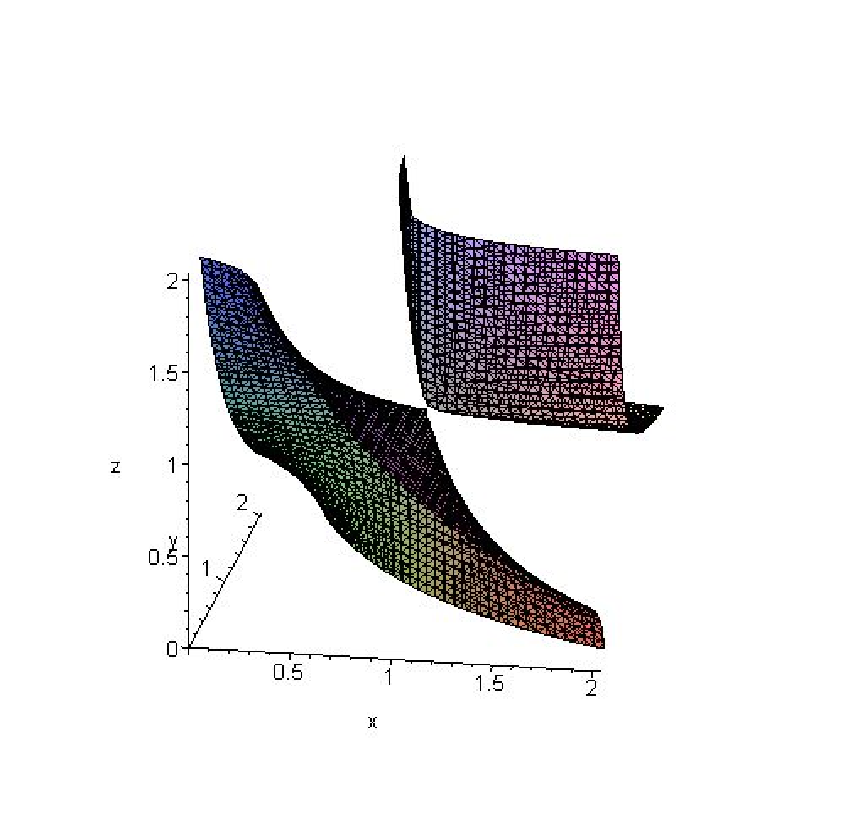}
\caption{an isolated quadratic singularity}
\label{fig:cone}
\end{figure}

In dimension three and above, there are many further possibilities. 
The simplest case not handled by previous techniques is that
of isolated quadratic singularity.  The purpose of this paper
is to address this type of generating function.  All the examples
Section~\ref{sec:applications} are of this type.
In fact, all the rational 3-variable generating functions we know of, 
that are not one of the types previous analyzed, have isolated,
usually quadratic, singularities.  The simplest case is when
the denominator is irreducible and its variety has a single, 
isolated quadratic singularity; a concrete example in dimension $d=3$
is the cube grove creation generating function, whose denominator
$Q = 1 + XYZ - (X+Y+Z+XY+YZ+ZX)/3$ has the zero set illustrated in 
figure~\ref{fig:cone}.  

The main results of this paper, Theorems~\ref{th:no plane} 
and~\ref{th:cone and plane} below, are asymptotic formulae 
for the coefficients of a generating function having a divisor
with this geometry.  In addition to one or more isolated quadratic 
singularities, our most general results allow $Q$ to be taken to 
an arbitrary real power and we allow the possibility of
other smooth divisors passing through the singularities of $Q$.
These generalizations complicate the exposition somewhat but 
are necessary to handle some of the motivating examples.

As a preview of the behavior of the coefficients, consider the 
case $d=3$ and $F = 1/Q$ illustrated in figure~\ref{fig:cone}.
The leading homogeneous term of $Q$ in the variables 
$(x,y,z) = (\log X , \log Y , \log Z)$ is $xy+xz+yz$.  
The outward normal cone (dual cone) at this point is the cone 
$\normal$ on which $Q^* \geq 0$, where $Q^* (r,s,t) = (r+s+t)^2 - 
2(r^2 + s^2 + t^2)$ is the dual quadratic to $Q$ (see
Section~\ref{ss:quadratic} for definitions).
The asymptotics for $a_\rr$ in this example are given
by Corollary~\ref{cor:simple} for $\rr$ in the interior
of $\normal$ and by an easier result (Proposition~\ref{pr:dir})
when $\rr \notin \normal$:
$$a_\rr \sim \begin{cases} C Q^* (\rr)^{-1/2} & \mbox{ if $\rr$ is
   in the interior of } \normal \\ 
   \mbox{exponentially small } & \mbox{ if } \rr \notin \normal 
   \end{cases}$$
The behavior of $a_\rr$ near $\partial \normal$ 
is more complicated and is not dealt with in this paper.  
The generating function $1/Q$ is the \Em{creation rate} generating 
function for cube groves, discussed in Section~\ref{ss:cube-grove-appl}.
The \Em{edge placement} generating function for cube groves (edge 
placement probabilities have more direct interpretations than do
creation rates) has an extra factor of $(1-Z)$
in the denominator.  Theorem~\ref{th:cone and plane} gives
the asymptotics in this case, for $\rr$ interior to $\normal$, as
$$a_\rr \sim C \arctan \theta (\rr)$$
where $\theta$ is a homogeneous degree $0$ function of $\rr$
which can be expressed in terms of dual quadratic form $Q^*$.  
Homogeneity of $\theta$ implies 
that that there is a limit theorem $a_{\lambda \rr}
\to \theta (\rhat)$ as $\lambda \to \infty$, where $\rhat := 
\rr / |\rr|$ is the unit vector in the direction $\rr$.

A total of five motivating applications will be discussed in detail
in Section~\ref{sec:applications}.  All of these may be seen
to have factors with isolated quadratic singularities.  There are
known trivariate rational generating functions with isolated
singularities that are not quadratic.  For example, the \Em{diabolo}
or \Em{fortress} tiling ensemble has an isolated quartic 
singularity~\cite{DGIP}.
Some of our results apply to this case, but a detailed analysis 
will be left for another paper.  The last example goes slightly beyond
what we do in this paper, but we include it because the analysis
follows largely the same methods.
\subsubsection*{Aztec diamond placement probability generating function
\protect{~\cite{JPS}}}
\begin{equation} \label{eq:aztec}
F (X,Y,Z) = \frac{Z/2}{(1-YZ) (1 - \frac{Z}{2} (X + X^{-1} + Y + Y^{-1})
   + Z^2)} \, .
\end{equation}
\subsubsection*{Cube groves edge probability generating function
\protect{~\cite{cube-grove}}}
\begin{equation} \label{eq:cube-grove}
F (X,Y,Z) = \frac{2 Z^2}{3 (1-Z) (1 + XYZ - \frac{1}{3} (X+Y+Z+XY+XZ+YZ))} 
   \, . \\
\end{equation}
\subsubsection*{Quantum random walk space-time probability generating function 
\protect{~\cite{QRW-one-dim,BBBP}}}
\begin{equation} \label{eq:quantumRW}
F (X,Y,Z) 
   = \frac{XZ - (1+XY)Z^2+YZ^3)}{(1-Z^2)(1-(X+X^{-1}+Y+Y^{-1}) Z/2 + Z^2)}
   \, . \\
\end{equation}
\subsubsection*{Friedrichs-Lewy-Szeg\"o graph polynomial
\protect{~\cite{scott-sokal}}}
\begin{equation} \label{eq:graphpoly}
F (X,Y,Z) = \left [ (1-X)(1-Y) + (1-X)(1-Z) + (1-Y)(1-Z) \right ]^{-\beta}
\end{equation}
\subsubsection*{Multi-set permutation generating function
\protect{~\cite{gessel-superballot}}}
\begin{equation} \label{eq:superballot}
G (X,Y,Z) = \frac{1}{1 - X - Y - Z + 4 XYZ}  \, .
\end{equation}

\subsection{Methods and organization} \label{ss:overview}

Our methods of analysis owe a great debt to two bodies of
existing theory.  Our approach to harmonic analysis of
cones is fashioned after the work of~\cite{ABG}.  We not
only quote their results on generalized Fourier transforms,
which date back somewhat farther to computations of~\cite{riesz}
and generalized function theory as described in~\cite{gelfand-generalized},
but we also employ their results on hyperbolic polynomials to produce 
homotopies of various contours.  Secondly, our understanding of the
existence of these homotopies has been greatly informed by
Morse theoretic results of~\cite{GM}.  We do not quote these 
results directly because our setting does not satisfy all their
hypotheses, but the idea to piece together deformations
local to strata is really the central idea behind stratified Morse
theory as explained in~\cite{GM}; see also the discussion of 
stratified critical points in Section~\ref{ss:dir}.

An outline our methods is as follows.  The chain of integration 
in the multivariate Cauchy integral~\eqref{eq:mv cauchy} is a $d$-dimensional
torus $T$ embedded in the \Em{complex torus} $\ctorus$, where
$\C^* := \C - \{0\}$.  Changing variables by $Z_j = \exp (z_j)$, the 
chain of integration becomes a chain $\CC$, the set of points with a fixed
real part.
Under this change of variables, the Cauchy integral~\eqref{eq:mv cauchy} 
becomes 
\begin{equation} \label{eq:new cauchy}
\int_{\CC} \exp (- \rr \cdot \zz) f (\zz) \, d\zz
\end{equation}
where $f := F \circ \exp$.
Letting $\zz := \xx + i \yy$, Morse theoretic considerations tell 
us we can deform the chain of integration so that it is supported 
by the region where $e^{-\rr \cdot \xx}$ is small (for large $|\rr|$) 
except near certain \Em{critical points}.
To elaborate, we can accomplish most of the deformation
by moving $\xx$.  The allowable region for such deformations 
of $\xx$ is a component of the complement to \Em{amoeba} of $F$ 
(see Section~\ref{ss:amoeba} for definitions).  Heuristically, 
we move $\xx$ to the 
support point $\xmax$ on the boundary of this region for
a hyperplane orthogonal to $\rr$ (see Sections~\ref{ss:dir} 
and those preceding for details).  Unfortunately, when $\xmax$ 
is on the boundary,~\eqref{eq:new cauchy} fails to be integrable.  
Ignoring this, however, and continuing with the heuristic, 
we let $q := Q \circ \exp$ and $h_j := H_j \circ \exp$ and we denote
the leading homogeneous parts of $q$ and $h_j$ by $\qt$ and $\hht_j$ 
respectively.  We then express $f$ near $\xmax$ as a series in 
negative powers of $\qt$ and $\hht_j$ (this is carried out in 
Section~\ref{ss:linearization}).  Integrating term by term,
each integral has the form
$$\int_\CC \exp (- \rr \cdot \zz) \frac{\zz^\mm}{\qt (\zz)^s \hht (\zz)^\nn}
   \, d\zz$$
where $\hht^\nn := \prod_{j=1}^k h_j^{n_j}$.
Replacing $\zz$ by $i \zz$, we recognize the Fourier transform of 
a product of a monomial with inverse powers of quadratics and
linear functions.  The Fourier transform of an inverse quadratic
is the dual quadratic and the Fourier transform of a linear
function is the Heaviside function.  The Fourier transform
of a product is a convolution.  These facts tell us what
result to expect.

Much of what has been described thus far is based on known 
methods and results, most of which are collected in the 
preliminary Section~\ref{sec:prelim}.  The bulk of the work, 
however, is in making rigorous these identities which involve 
Fourier integrals that do not converge, taken over regions 
which are not obviously deformations of each other (the part
above where we said, ``ignoring this, $\ldots$'').  
For this purpose, some carefully chosen deformations are 
constructed, based largely on deformations found in
Sections~5 and~6 of~\cite{ABG}.  Specifically, we use
results on hyperbolic polynomials (see Section~\ref{ss:hyp}
for definitions) established in~\cite{ABG} and elsewhere,
to construct certain vector fields on $\C^d$.  These vector 
fields, based on the construction of~\cite[Section~5]{ABG} 
and described in 
our Section~\ref{ss:VF}, then allow us to construct
deformations in Section~\ref{ss:hom}, which 
satisfy several properties.
First, they enact what Morse theory
has guaranteed: they push the chain of integration
to where the integrand of~\eqref{eq:new cauchy} is very small,
except near critical points, as in figure~\ref{fig:cone 1}
below.  Secondly, they do this without intersecting $\sing$, thereby
allowing the integral to remain the same.  This localizes 
the integral to the critical points, and allows us to concentrate 
on one critical point at a time.  The resulting chain of integration
is depicted in figure~\ref{fig:cone 1}.
\begin{figure}[ht] \centering
\subfigure[in log space]
{\includegraphics[scale=0.50]{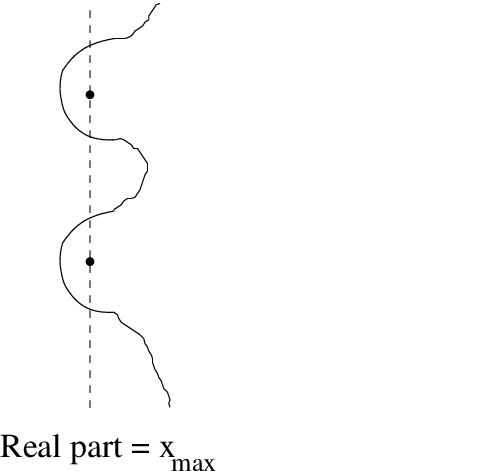} \label{sf:localizing}}
\subfigure[in the original domain]
{\hspace{1.2in} \includegraphics[scale=0.50]{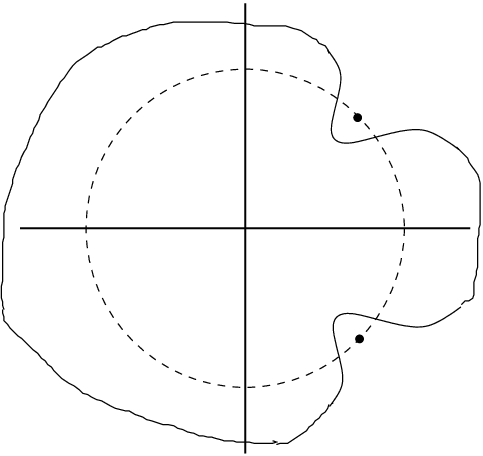}}
\caption{the localizing chain, in logarithmic and original coordinates}
\label{fig:cone 1}
\end{figure}

\begin{figure}[ht]
\hspace{1.2in} \includegraphics[scale=0.50]{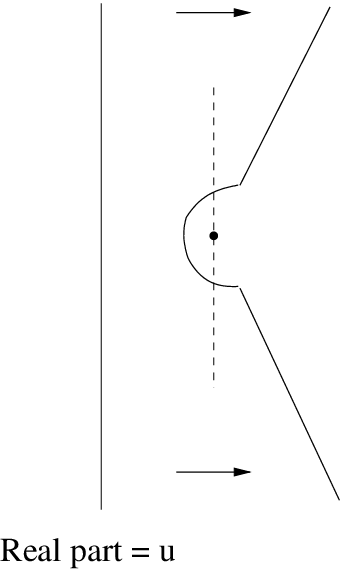}
\caption{the projective chain} 
\label{fig:cone 3}
\end{figure}

Thirdly, they allow us to ``straighten out'' the chain 
of integration.  Figure~\ref{fig:cone 3} shows that the 
chains in figure~\ref{fig:cone 1}, as well as the original
chain, are homotopic near the critical point to a (slightly perturbed) conical
chain.

Combined with the series expansion by homogeneous functions, 
this reduces all necessary integrals to a small class of 
Fourier-type integrals.  Many of these are evaluated
as \Em{generalized functions} in~\cite{ABG,riesz,gelfand-generalized} 
and elsewhere.  In Section~\ref{ss:gen func} and~\ref{ss:fourier}
we summarize the relevant facts about generalized functions.
The above deformations allow us to show, in 
Sections~\ref{ss:cone only}~--~\ref{ss:cone plane}, 
that these generalized functions, defined as integrals over
the straight contour on the left of figure~\ref{fig:cone 3},
do approximate the integrals we are interested in, which we must
evaluate over the chains shown on the right of figure~\ref{fig:cone 1}
and on figure~\ref{fig:cone 3} in order for the localizations to remain 
valid.  Not all of the computations we need are available in
the literature.  In Section~\ref{ss:cone plane} we use a 
construction from~\cite{ABG}, the \Em{Leray cycle}, along
with a residue computation, to reduce the Fourier transform 
of $1/(\qt \cdot \hht)$ to an explicitly computable one-dimensional 
integral.  It is this computation that is responsible for
the explicit asymptotic formula for placement probabilities
in the Aztec Diamond and Cube Grove problems.

To summarize, the organization of the rest of the paper is 
as follows.  Section~\ref{sec:prelim} defines some notation in use 
throughout the paper, and collects preliminary results
on amoebas, convex duals, hyperbolicity, and expansions
by powers of homogeneous polynomials.  Section~\ref{sec:results} 
states the main results.  Section~\ref{sec:applications} has 
five subsections, each discussing one of the five examples.  
The next two sections are concerned with the proofs of the 
main results.  Section~\ref{sec:homotopies} constructs 
homotopies that shift contours of integration, while 
Section~\ref{sec:integrals} evaluates several classes of
integrals via the theory of generalized Fourier transforms.
Finally, Section~\ref{sec:further} concludes with 
a discussion of open problems and further research directions.

\subsection{Comparison with other techniques}

One might ask, in a paper of this length, whether this is 
the best way to obtain these results.  To answer this, we
briefly review comparable published results and a relevant
unpublished failure.  The Arctic Circle Theorem for tilings
of the Aztec diamond was proved in~\cite{aztec} in two steps.
First formulae for the coefficients of the simpler \Em{creation rate} 
generating function $(1 - (X + X^{-1} + Y + Y^{-1}) Z/2 + Z^2)^{-1}$ 
were derived via a relation to Krawtchouk polynomials.  Secondly,
these were summed by means of contour integrals.  The computation
was quite specialized, and did not generalize even to the 
nearly formally equivalent case of cube groves.  In fact, 
for the cube grove model, up to now only the easy half of 
the Arctic Circle theorem was proved (exponential decay outside
of the circle).  

The present paper takes the view that the work is justifiable
if we can then crank out results with relatively little effort.  
The continuous setting clarifies matters considerably.  Our fundamental
result, Theorem~\ref{th:no plane} below, is that the asymptotics of
a generating functions with irreducible quadratic denominator, 
such as in~\eqref{eq:graphpoly},
are its Fourier transform, which is the dual quadratic.  This
is the continuous analogue of the Krawtchouk polynomials that
appear when the computations are done in the discrete setting.
Multiplying the denominator by $1/h$ where $h$ is smooth, 
corresponds to convolving the Fourier transform with a 
heavyside function; this is the analogue of summing and
is made rigorous in Theorem~\ref{th:cone and plane}.  
These two facts (Fourier transform of a cone is the dual cone, 
and Fourier transform of $1/(Qh)$ is the integral of the 
Fourier transform of $1/Q$) are well known, which makes 
the resulting computation predictable although a number
of details need to be addressed.

As far as we know, the Aztec diamond result
is the only one of the five results in Section~\ref{sec:applications}
that was previously known; all five examples are easily handled
once the machine is built.   A sixth example, the fortress tiling 
ensemble, can be analyzed by our methods but explicit formulae
in this case require further computations along the lines of 
Section~\ref{ss:cone plane}.  The limit theorem is known in this 
case, the result of a variational equation which is established 
and explicitly solved in the beautiful paper~\cite{kenyon-okounkov}.

It should be emphasized that evaluating the 
integral~\eqref{eq:new cauchy} involves interplay between 
the form and the chain, and this interplay which is primarily 
responsible for failure of several earlier attempts to analyse 
the asymptotics of the integral.  To be sure, resolution of 
singularities provides one with an efficient toolbox for reducing 
the integrand to a monomial form; see for 
instance~\cite[Chapter~7]{arnold-singularities2}
or~\cite{jeanquartier70}; the resolution is in principle
effective~\cite{bierstone-milman} and the algorithm given
in~\cite{varchenko76} suffices in most cases.  However, it then
becomes difficult to control the chain of integration. A few years ago, 
the second author in collaboration with H. Cohn attempted to
resolve the singularity and compute the integral.  The resolution
was indeed computable.  Unfortunately, as will be true in all such 
cases, the phase function becomes quite degenerate, being constant 
along the exceptional divisor of the resolution.  The resulting 
integral was beyond Cohn and Pemantle's ability to evaluate.

Another important aspect of the problem that is difficult 
to control using resolution of singularities is the real structure. 
Consequences of this include hyperbolicity of the tangent cones to the
pole variety at critical points (see Section \ref{ss:hyp curved}). 
As we will see, this hyperbolicity plays a critical role in the 
constructions of the deformations of the integration chains.

\setcounter{equation}{0}
\section{Notation and preliminaries} \label{sec:prelim}

Several notions arising repeatedly in this paper are
the logarithmic change of variables, duality between
$\rr$ and $\zz$, and the leading homogeneous part of
a function.  We employ some meta-notation designed for
ease of keeping track of these.  We use upper case
letters for variables and functions in the complex torus $\ctorus$, 
and lower case letters in the logarithmic coordinates.  
We will never use the notations without defining
them, but knowing that, for example, $F \circ \exp$ will
always be denoted $f$ and $\zz$ will always be $\log \ZZ$, 
should help the reader quickly recognize the setting.
We will always use $\xx$ and $\yy$ for the real and 
imaginary parts of the vector $\zz \in \C^d$.
Boldface is reserved for vectors.  The leading homogeneous 
part of a function is denoted with a bar.  Rather than
considering the index $\rr$ of $a_\rr$ to be an element of
$\Z^d$, we consider it to be an element of a space $\Rdd$
that is dual to the domain $\C^d$ in which $\zz$ lives, 
with respect to the pairing $\rr \cdot \zz$ (the space
$\Rdd$ is a subset of the full dual space $(\C^d)^*$ but
all our dual vectors will be real).  
Many functions of $\rr$ use in what follows are homogeneous degree $0$;
letting $\rhat$ denote the unit vector $\rr / |\rr|$ we will
often write these as functions of $\rhat$.  The logarithm
and exponential functions are extended to act coordinatewise
on vectors.  Thus
\begin{eqnarray*}
\exp (z_1 , \ldots , z_d) & := & (\exp (z_1) , \ldots , \exp (z_d)) \, ; \\
\log (Z_1 , \ldots , Z_d) & := & (\log (Z_1) , \ldots , \log (Z_d)) \, .
\end{eqnarray*}
We also employ the slightly clunky notation \label{relog}
$$\Log \ZZ \, := \, (\log |Z_1| , \ldots , \log |Z_d|) 
   \, = \, \Real \{ \log \ZZ \}$$
for the coordinatewise log-modulus map, having found that 
the notations in use in~\cite{GKZ} do not allow for quick
visual distinction between $\log$ and $\Log$.  

Our chief concern is with generating functions that are
ordinary power series, convergent on the unit polydisk, 
and whose denominator is the product of smooth and 
quadratically singular factors that intersect the closed 
but not the open unit polydisk.  It costs little, however,
and there is some benefit to work in the greater generality
of Laurent series representing functions with polynomial denominators.  
Indeed, Laurent series arise naturally in the 
examples (though these Laurent series have exponents in proper cones, 
and may therefore be reduced by log-affine changes of coordinates 
to Taylor series).

\begin{defn}[homogeneous part] \label{def:homog}
For analytic germ $f: (\C^d,\zz) \to \C$ at a point
$\zz \in \C^d$, we let $\deg (f , \zz)$ denote the degree of \label{deg}
vanishing of $f$ at $\zz$.  This is zero if $f(\zz) \neq 0$
and in general is the greatest integer $n$ such that
$f(\zz + \ww) = O(|\ww|^n)$ as $\ww \to \zero$.  Also, $\deg (f,\zz)$
is the least degree of any monomial in the ordinary power series expansion
of $f (\zz + \cdot)$ around $\zero$.  We let $\homog (f , \zz)$ denote
the sum of all monomials of minimal degree in the power series
for $f(\zz + \cdot)$ and we call this the \Em{homogeneous part 
of $f$ at $\zz$}.  Thus 
$$f(\zz + \ww) = \homog (f, \zz) (\ww) 
   + O \left ( |\ww|^{\deg (f,\zz) + 1} \right ) \, $$ 
for small $|\ww|$.When $\zz = \zero$, we may omit $\zz$ from the notation: thus,
$\homog(f) := \homog(f , \zero)$. \label{homog}
\end{defn}

A number of connections between zeros of a Laurent polynomial $F$, 
the Laurent series for $1/F$, the Newton polytope for $F$ and 
certain dual cones to this polygon were worked out in the 1990's
by Gelfand, Kapranov and Zelivinsky.  We summarize 
some relevant results from~\cite[Chapter~6]{GKZ}.  

\subsection{The Log map and amoebas} \label{ss:amoeba}

\begin{figure}[ht] \centering
\subfigure[$\amoeba(2-X-Y)$]
{\includegraphics[scale=0.40]{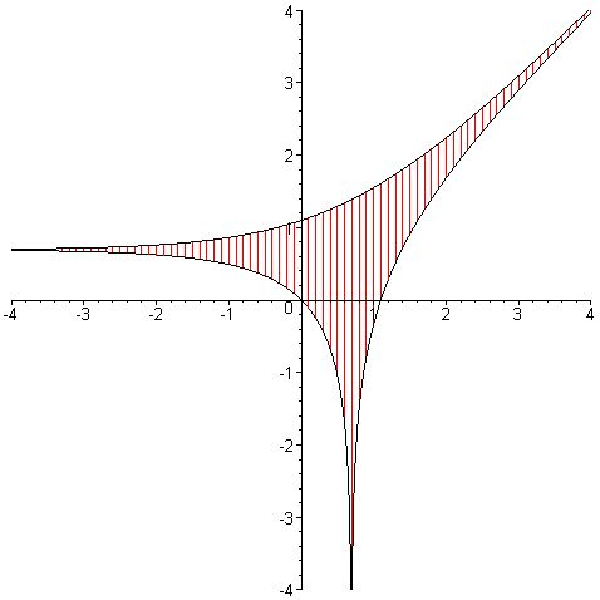} \label{sf:amoeba01}}
\hspace{0.5in}
\subfigure[$\amoeba(3-2X-Y)(3-X-2Y)$]
{\includegraphics[scale=0.40]{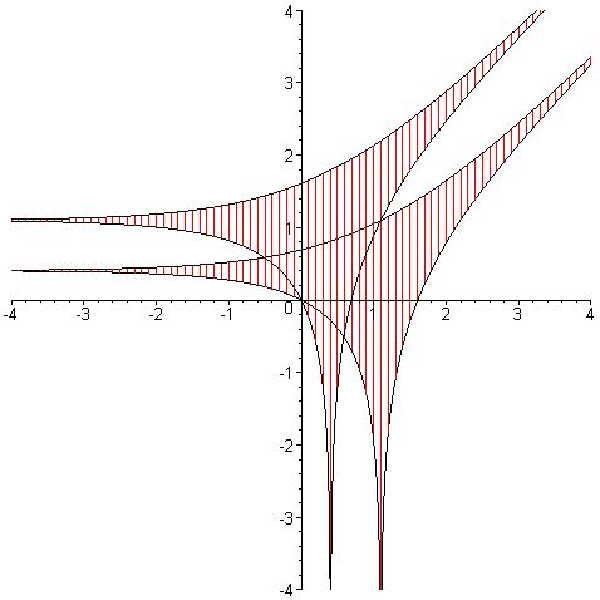} \label{sf:amoeba02}}
\caption{two amoebae}
\label{fig:amoeba01}
\end{figure}

Let $F$ be a Laurent polynomial in $d$ variables.  
Let $\ctorus \subseteq \C^d$ denote the $d$-tuples of
nonvanishing complex numbers and let  $\sing_F$ denote  \label{sing}
the zero set of $F$ in $\ctorus$.
Following~\cite{GKZ} we define the \Em{amoeba} of $F$ to
be the image under $\Log$ of the zero set of $F$: 
$$\amoeba (F) := \{ \Log \zz : \zz \in \sing_F  \cap \ctorus\} 
   \subseteq \R^d \, .$$
The simplest example is the amoeba of a linear function,
such as $F = 2-X-Y$, shown in figure~\ref{sf:amoeba01}.
The amoeba of a product is the union of amoebas, as shown 
in figure~\ref{sf:amoeba02}. \label{amoeba}

The rational function $1/F$ has, in general, a number of
Laurent series expansions, each convergent on a different
subset of $\C^d$.  Combining Corollary~1.6 in Chapter 6 of~\cite{GKZ}
with Cauchy's integral theorem, we have the following result.

\begin{pr} \label{pr:amoeba}
The connected components of $\R^d \setminus \amoeba (F)$ are
convex open sets.  The components are in bijective correspondence 
with Laurent series expansions for $1/F$ as follows.  Given a 
Laurent series expansion of $1/F$, its open domain of convergence 
is precisely $\Log^{-1} B$ where $B$ is a component of $\R^d
\setminus \amoeba (F)$.  Conversely, given such a component $B$,
a Laurent series $1/F = \sum a_\rr \ZZ^\rr$ convergent on $B$ may 
be computed by the formula
$$a_\rr = \frac{1}{(2 \pi i)^d} \int_{\bf T} \ZZ^{-\rr - \one} 
   \frac{1}{F(\ZZ)} \, d\ZZ$$
where ${\bf T}$ is the torus $\Log^{-1} (\xx)$ for any $\xx \in B$. 
Changing variables to $\ZZ = \exp (\zz)$ and
$d\ZZ = \ZZ \, d\zz$ gives
\begin{equation} \label{eq:cauchy torus}
a_\rr = \frac{1}{(2 \pi i)^d} \int_{\xx + i \flattorus}
   e^{- \rr \cdot \zz} \frac{1}{f (\zz)} \, d\zz
\end{equation}
where $f = F \circ \exp$ and $\flattorus$ is the torus 
$R^d / (2 \pi \Z)^d$. \label{flattorus}
$\Cox$
\end{pr}

\subsection{Dual cones, tangent cones and normal cones}
\label{ss:cones}

Let $\Rdd$ denote the dual space to $\R^d$ 
and for $\yy \in \R^d$ and $\rr \in \Rdd$, use the
notation $\rr \cdot \yy$ to denote the pairing.  
\label{dual}
Let $\cL$ be any convex open cone in $\R^d$.  The (closed) 
convex dual cone $\cLd \subseteq \Rdd$ is defined to be the set of 
vectors $\vv \in \Rdd$ such that $\vv \cdot \xx \geq 0$ for all $\xx
\in \cL$.  Familiar properties of the dual cone are
\begin{eqnarray}
L \subseteq M & \Rightarrow & L^* \supseteq M^* \, ; \label{eq:dual contain} \\
(L \cap M)^* & = & \hull (L^* \cup M^*) \, . \label{eq:dual hull}
\end{eqnarray}
Suppose $\xx$ is a point on the boundary of a 
convex set $C$.  Then the intersection of all halfspaces
that contain $C$ and have $\xx$ on their boundary is a
closed convex affine cone with vertex $\xx$ (a translation by $\xx$ 
of a closed convex cone in $\R^d$) that contains $C$.  Translating by $-\xx$ 
and taking the interior gives the (open) \Em{tangent cone} to $C$ at $\xx$, 
denoted by $\tangent_\xx (C)$.  An alternative definition is  
\label{tan}
$$\tangent_\xx (C) = \{ \vv : \xx + \ee \vv \in C
   \mbox{ for all sufficiently small } \ee > 0 \} \, $$
(where $B$ is the unit ball).   
The (closed) \Em{normal cone} to $C$ at $\xx$, denoted 
$\normalx (C)$, is the convex dual cone to the negative of the 
tangent cone:
\label{normalx}
$$\normalx (C) = (- \tangent_\xx (C))^* \, .$$
Equivalently, it corresponds to the set of linear functionals 
on $C$ that are maximized at $\xx$, or to the set of outward 
normals to support hyperplanes to $C$ at $\xx$.
\label{xmax}
\begin{defn}[proper dual direction] \label{def:proper}
Given a convex set $C$, say that $\rr$ is a proper direction
for $C$ if the maximum of $\rr \cdot \xx$ on $L$ is achieved
at a unique $\xmax \in \overline{C}$.  We call $\xmax$
the dual point for $\rr$.  The set of directions $\rr$ for which 
$\rr \cdot \xx$ is bounded on $C$ but $\rr$ is not proper has 
measure zero.
\end{defn}
The term {\em tangent cone} has a different meaning in algebraic 
contexts, which we will require as well.  (The term
{\em normal cone} has an algebraic meaning as well, which
we will not need.)
To avoid confusion, we define the \Em{algebraic tangent cone}
of $f$ at $\xx$ to be $\sing_{\homog (f,\zz)}$.  
An equivalent but more
geometric definition is that the algebraic tangent cone is
the union of lines through $\xx$ that are the limits of secant
lines through $\xx$; thus for a unit vector $\uu$, the line 
$\xx + t \uu$ is in the algebraic tangent cone if there are 
$\xx_n \in \sing_f$ distinct from but converging to $\xx$ 
for which $(\xx_n - \xx) / ||\xx_n - \xx|| \to \pm \uu$.
This equivalence and more is contained in the following
results.  We let $S_1$ denote the unit sphere 
$\{ (z_1 , \ldots , z_d) : |z_1|^2 + \cdots + |z_d|^2 = 1 \}$
and let $S_r := r S_1$ denote the sphere of radius $r$.

\begin{lem}[algebraic tangent cone is the limiting secant cone] 
\label{lem:homog approx}
Let $q$ be a polynomial vanishing to degree $m \geq 1$ at the
origin and let $\qt := \homog (q)$ be its homogeneous part;
in particular,
$$q(\zz) = \qt (\zz) + R(\zz)$$
where $\qt$ is a nonzero homogeneous polynomial of degree $m$
and $R(\zz) = O(|\zz|^{m+1})$.  
Let $q_\ee$ denote the polynomial
$$q_\ee (\zz) := \ee^{-m} q (\ee \zz) = \qt (\zz) + R_\ee (\zz)$$
where $R_\ee (\zz) = \ee^{-m} R(\ee \zz) \to 0$ as $\ee \to 0$. 
Let $\sing_\ee := \sing_{q_\ee} \cap S_1$ denote the intersection 
of $\{ q_\ee = 0 \}$ with the unit sphere.  Then $\sing_\ee$ 
converges in the Hausdorff metric as $\ee \to 0$ to $\sing_{\qt}
\cap S_1$.
\end{lem}

\noindent{\sc Proof:} On any compact set, in particular $S_1$,
$R_\ee \to 0$ uniformly.  If $\zz^{(n)} \to \zz$ and $\zz^{(n)} 
\in \sing_{1/n}$ then for each $n$, 
$$|\qt (\zz^{(n)})| = |q_{1/n} (\zz^{(n)}) + R_{1/n} (\zz^{(n)})| = 
   |R_{1/n} (\zz^{(n)})| \to 0 \, .$$
Hence $\qt (\zz) = 0$ by continuity of $\qt $ and we see that any
limit point of $\sing_\ee$ as $\ee \to 0$ is in $\sing_{\qt} \cap S_1$.   
Conversely, fix a unit vector $\zz \in \sing_{\qt}$.  The homogenous 
polynomial $\qt$ is not identically zero, therefore there is a projective 
line through $\overline{\zz}$ along which $\qt$ has a zero of finite order 
at $\overline{\zz}$.  Back in affine space, there is a complex
curve $\gamma$ in the unit sphere along which $\qt$ is holomorphic 
with a zero of some finite order $k$ at $\zz$.  As $\ee \to 0$,
the holomorphic function $R_\ee$ goes to zero uniformly in a 
neighborhood of $\zz$ in $\gamma$; hence there are $k$ zeros of 
$q_\ee$ converging to $\zz$ as $\ee \to 0$, and therefore $\zz$ 
is a limit point of $\sing_\ee$ as $\ee \to 0$.
$\Cox$

\subsection{Hyperbolicity for homogeneous polynomials} \label{ss:hyp}

The notion of hyperbolic polynomials arose first 
in~\cite{garding-hyperbolic} in connection with solutions 
to wave-like partial differential equations.  The same
property turns out to be very important as well for 
convex programming, cf.~\cite{gulen} from which much of the
next several paragraphs is drawn.

Let $f$ be a complex polynomial in $d$-variables and $f(D)$ denote 
the corresponding linear partial differential operator with 
constant coefficients, obtained by replacing each $x_j$ by
$\partial / \partial x_j$.  For example, if $f$ is the standard
Lorentzian quadratic $\lorentz (\xx) := x_1^2 - \sum_{j=2}^d x_j^2$ 
then $\lorentz (D)$ is the wave operator $(\partial / \partial x_1)^2 - 
\sum_{j=2}^d (\partial / \partial x_j)^2$.  G{\aa}rding's object was 
to determine when the equation 
\begin{equation} \label{eq:LPDE}
f(D) u = g 
\end{equation}
with $g$ supported on a halfspace has a solution supported 
in the same halfspace.  The wave operator has this property, 
and in fact there is a unique such solution for any such $g$.
It turns out that the class of $f$ such that~\eqref{eq:LPDE}
always has a solution supported on the halfspace is precisely
characterized by the property of being \Em{hyperbolic}, as 
defined by G{\aa}rding.  
In this case, it was later shown~(\cite[Theorem~12.4.2]{hormander} 
that the solution is in fact unique.  The theory of hyperbolic 
polynomials serves in the present paper to prove the existence
of deformations of chains of integration past points of the
pole manifold at which the pole polynomial is locally hyperbolic. 

We begin with hyperbolicity for homogeneous polynomials,
which is a simpler and better developed theory.  We use $A$ 
rather than $f$ for a homogeneous polynomial.
\begin{defn}[hyperbolicity]
Say that a homogeneous complex polynomial $A$ of degree $m \geq 1$
is hyperbolic in direction $\vv \in \R^d$ if $A(\vv) \neq 0$
and the polynomial $A(\xx + t \vv)$ has only real roots 
when $\xx$ is real.  In other words, every line in real space 
parallel to $\vv$ intersects $\sing_A$ exactly $m$ times 
(counting multiplicities).  
\end{defn}

While seemingly weaker, the requirement of avoiding purely
imaginary roots is in fact easily seen to be equivalent.
\begin{pr} \label{pr:no imaginary roots}
Hyperbolicity of the homogeneous polynomial $A$ in the direction $\vv$
is equivalent to the condition that $A(\vv + i \yy) \neq 0$
for all $\yy \in \R^d$.  
\end{pr}

\noindent{\sc Proof:} Because $A$ is homogeneous, when $\lambda \neq 0$,
we have $A(\lambda \zz) = 0$ if and only if $A(\zz) = 0$.  With
$\lambda = i \cdot s$, a purely imaginary number not equal to zero, 
we see that $A(\vv + i \yy) \neq 0$ for all $\yy \in \R^d$ is equivalent
to $A(\yy + i s \vv) \neq 0$ for all $\yy \in \R^d$ and all
nonzero real $s$.  This becomes $A(\yy + t\uu + i s \uu) \neq 0$ 
for all $\yy \in \R^d$ and real $s \neq 0$; writing $z = t + i s$, this
is equivalent to $A(\yy + z \uu) \neq 0$ when $z$ is not real, which
is the definition of hyperbolicity.   $\Cox$

The further properties we need are well known and are proved,
among other places, in~\cite[Theorem~3.1]{gulen}.  

\begin{pr} \label{pr:hyp properties}
The set of $\vv$ for which $A$ is hyperbolic in direction $\vv$
is an open set whose components are convex cones.  Denote by 
\label{hypcone}
$\cK^\vv (A)$ the connected component of this cone that contains a 
given $\vv$.  Some multiple of $A$ is positive on $\cK^\vv (A)$ and 
vanishing on $\partial \cK^\vv (A)$, and for $\xx \in \cK^\vv (A)$, 
the roots of $A(\xx + t \vv)$ will all be negative.  
$\Cox$
\end{pr}

Semi-continuity properties for cones of hyperbolicity play a 
large role in the construction of deformations.  A lower
semi-continuous function $f$ satisfies $f(x) \leq \liminf_{y \to x}
f(y)$.  The property is important in elementary analysis
because a lower semi-continuous function on a compact set 
achieves its infimum; generalizing to set-valued functions, 
the conclusion is roughly that the empty set is not a limit
value and therefore that a continuous section can be defined.
In this section, we develop semi-continuity properties for cones 
of hyperbolicity (a topic that occupies many pages of~\cite{ABG}).  

The following proposition and definition define a family
of cones $\disp{\{ \cK^{A,B} (\xx) \}_{\xx \in \R^d}}$ 
which will be used to prove two critical semi-continuity 
results for cones of hyperbolicity for log-Laurent polynomials 
(Theorem~\ref{th:semi} below).
\begin{pr}[first semi-continuity result]
\label{pr:relinearization}
Let $A$ be any hyperbolic homogeneous polynomial, and let $m$ be 
its degree.  Fix $\xx$ with $A(\xx) = 0$ and let $\At := \homog (A , \xx)$
denote the leading homogeneous part of $A$ at $\xx$.  If
$A$ is hyperbolic in direction $\uu$ then $\At$ is also
hyperbolic in direction $\uu$.  Consequently, if $B$ is
any cone of hyperbolicity for $A$ then there is some
cone of hyperbolicity for $\At$ containing $B$.  
\end{pr}
\noindent{\sc Proof:} This follows from the conclusion~(3.45)
of~\cite[Lemma~3.42]{ABG}.
Because the development there is long and complicated, we
give here a short, self-contained proof, provided by J. Borcea 
(personal communication).  If $P$ is a polynomial whose degree
at zero is $k$, we may recover its leading homogeneous part
$\homog(P)$ by
$$\homog (P) (\yy) = \lim_{\lambda \to \infty} 
   \lambda^k P (\lambda^{-1} \yy) \, .$$
The limit is uniform as $\yy$ varies over compact sets.
Indeed, monomials of degree $k$ are invariant under the scaling
on the right-hand side, while monomials of degree $k+j$ scale
by $\lambda^{-j}$, uniformly over compact sets.  

Apply this with $P (\cdot) = A(\xx + \cdot)$ and $\yy + t \uu$
in place of $\yy$ to see that for fixed $\xx, \yy$ and $\uu$,
$$\At (\yy + t \uu) = \lim_{\lambda \to \infty} 
   \lambda^{k} A(\xx + \lambda^{-1}(\yy + t \uu))$$
uniformly as $t$ varies over compact sub-intervals of $\R$.
Because $A$ is hyperbolic in direction $\uu$, for any fixed
$\lambda$, all the zeros of this polynomial in $t$ are real.
Hurwitz' theorem on the continuity of zeros~\cite[Corollary~2.6]{conway} 
says that a limit, uniformly on bounded intervals, of polynomials 
having all real zeros will either have all real zeros or vanish 
identically.  The limit $\At (\yy + t \uu)$ has degree $k \geq 1$; 
it does not vanish identically and therefore it has all real zeros.  
This shows $\At$ to be hyperbolic in direction $\uu$.
$\Cox$

\begin{defn}[family of cones in the homogeneous case]
\label{def:KAB}
Let $A$ be a hyperbolic homogeneous polynomial and let $B$
be a cone of hyperbolicity for $A$.  If $A(\xx) = 0$, define
$$\cK^{A,B} (\xx)$$
\label{hypcone2}
to be the cone of hyperbolicity of $\homog (A,\xx)$ containing
$B$, whose existence we have just proved.  If $A(\xx) \neq 0$
we define $\cK^{A,B} (\xx)$ to be all of $\R^d$.
\end{defn}

As an example of a hyperbolic polynomial, let $S = x_1^2 - x_2^2
- \cdots - x_d^2$ be the standard Lorentzian quadratic.  Then 
$\cK^{e_1}(S)$ is the Lorentz cone $\{ \xx : x_1 \geq \sqrt{x_2^2 + 
\cdots + x_d^2} \}$.  Any quadratic of Lorentizan signature is
obtained from this one by a real linear transformation; we see 
therefore that for any Lorentzian quadratic, the boundary of 
the cone of hyperbolicity is the algebraic tangent cone.

The class of hyperbolic polynomials in a given direction
direction $\vv$ is closed under products, and $\cK^\vv (A A')
= \cK^\vv (A) \cap \cK^\vv (A')$.  The class contains all 
linear polynomials not annihilating $\vv$ and all real 
quadratic polynomials $p$ of Lorentzian signature for which 
$p(\vv) > 0$ ($\vv$ is time-like).  

\subsection{Hyperbolicity and semi-continuity for log-Laurent polynomials 
on the amoeba boundary} 
\label{ss:hyp curved}

For a function that is not locally homogeneous, there are
two natural generalizations of the definition of hyperbolicity.
Both are equivalent to the notion of hyperbolicity already
defined, in the case of a homogeneous polynomial.  Useful
features of the two definitions are revealed in the subsequent
two propositions.
\begin{defn} \label{def:inhomogeneous hyperbolicity}
Let $f: \C^d \to \C$ vanish at $\zz$ and be holomorphic
in a neighborhood of $\zz$.  We say that $f$ is \Em{strongly
hyperbolic} at $\zz$ in direction of the unit vector $\vhat$
if there is an $\ee > 0$ such that $f(\zz + t \vv' + i \uu) \neq 0$ 
for all real $t \in (0, \ee)$, all $\vv'$ with $|\vv' - \vhat| < \ee$,
and all $\uu \in \R^d$ of magnitude at most $\ee$.  In this
case we may say that $f$ is strongly hyperbolic at $\zz$
in direction $\vhat$ \Em{with radius} $\ee$.
Say that $f$ is \Em{weakly hyperbolic} in direction $\vv$ if 
for every $M > 0$, there is an $\ee > 0$ such that
$f(\zz + t \vv + i \uu) \neq 0$ for all real $0 < t |\vv| < \ee$, 
and for all $\uu \in \R^d$ of magnitude at most $\ee$ additionally
satisfying $|\uu| / (t |\vv|) \leq M$.  
\end{defn}

\begin{pr} \label{pr:weak hyperbolicity}
Let $A = \homog (f , \zz)$.  Then $A$ is hyperbolic in direction
$\uu$ if and only if $f$ is weakly hyperbolic in direction $\uu$
at $\zz$.  
\end{pr}

\noindent{\sc Proof:} The homogeneous polynomial $A$ fails to
be hyperbolic at in direction $\uu$ if and only if there
is some real $\yy$ such that $A(\uu + i \yy) = 0$.  By 
Lemma~\ref{lem:homog approx}, this happens if and only if
$f(\zz + \ww_n) = 0$ for some sequence $\{ \ww_n \}$ 
converging to $\zero$ with $\ww_n / |\ww_n|$ converging
to $(\uu + i \yy) / |\uu + i \yy|$.  This is equivalent
to failure of weak hyperbolicity of $f$ at $\zz$ in 
direction $\uu$.   $\Cox$

\begin{unremark} 
It is immediate from the definition that strong hyperbolicity
is a neighborhood property: if $f$ is strongly hyperbolic
at $\xx + i \yy$ in direction $\vhat$ with radius $\ee$,
then for $|\yy' - \yy| < \ee$ and $|\vhat' - \vhat| < \ee$, 
$f$ is strongly hyperbolic at $\xx + i \yy'$ in direction 
$\vhat'$ with direction $\ee - \max \{ |\yy' - \yy| , |\vv' - \vhat| \}$.
Weak hyperbolicity of $f$ at $\zz$ in direction $\vv$ extends to 
a neighborhood of $\vv$ by Propositions~\ref{pr:hyp properties}
and~\ref{pr:weak hyperbolicity}.  Extending weak hyperbolicity
to neighboring $\zz$ is much trickier. 
\end{unremark}

\begin{pr} \label{pr:hyperbolic}
Let $F$ be a Laurent polynomial in $d$ variables.  Suppose 
that $B$ is a component of $\amoeba (F)$ and $\xx \in 
\partial B$, so that $f := F \circ \exp$ vanishes at some point
$\xx + i \yy$.  Let $\ft := \homog (f , \xx + i \yy)$ 
\label{ft}
denote the leading homogeneous part of $f (\xx + i \yy + \cdot)$.  
Then $f$ is strongly hyperbolic at $\xx + i \yy$, some complex 
scalar multiple of $\ft$ is real and hyperbolic, and some cone of 
hyperbolicity $\cK^\uu (\ft)$ contains $\tangent_\xx (B)$.  
\end{pr}

\noindent{\sc Proof:} Strong hyperbolicity of $f$ in 
any direction $\uu \in \tangent_\xx (B)$ follows
from the definition of the amoeba.  Strong hyperbolicity
is stronger than weak hyperbolicity, hence hyperbolicity of $\ft$
in direction $\uu$ follows from Proposition~\ref{pr:weak hyperbolicity}. 
The vector $\uu \in \tangent_\xx (B)$ is arbitrary, whence
$\cK^\uu (\ft) \supseteq \tangent_\xx (B)$.  To see that some 
multiple of $\ft$ is real, 
let $\uu$ be any real vector in $\tangent_\xx (B)$, let 
$m$ denote the degree of $\ft$, and let $\gamma$ denote the 
coefficient of the $z^m$ term of $A (\ft \uu + \yy)$.  Then 
$\gamma$ is the degree $m$ coefficient of $\ft (z \uu)$, 
hence is nonzero and does not depend on $\yy$.  For any 
fixed $\yy$, the fact that $\ft (z \uu + \yy)$ has all real 
roots implies that the monic polynomial $\gamma^{-1} \ft (z \uu + \yy)$ 
has all real coefficients.   $\Cox$

\begin{defn}[hyperbolicity and normal cones at a point of $\sing_f$]
Let $F$ be a Laurent polynomial, $B$ a component of $\R^d \setminus
\amoeba (F)$, and $\ZZ = \exp (\xx + i \yy) \in \sing_f$ with
$\xx \in \partial B$.  We let $f := F \circ \exp$ and let
\begin{equation} \label{eq:hom cone}
\cK^{f,B} (\ZZ) := \cK^\uu (\homog (f , \xx + i \yy)) \, ,
\end{equation}
\label{hypcone3}
denote the (open) cone of hyperbolicity of $\ft := \homog (f , \xx + i \yy)$
that contains $B$, whose existence is guaranteed by 
Proposition~\ref{pr:hyperbolic}.  Although it is a slight
abuse of notation, we also write
$$\cK^{f,B} (\yy) := \cK^{f,B} (\ZZ)$$
when $\ZZ = \exp (\xx + i \yy)$ and the specification of $\xx$ 
is clearly understood.  We also define the normal cone
\label{normal}
\begin{equation} \label{eq:normal cone}
\normal (\ZZ) := (\normal)^{f , B} (\ZZ) = (\cK^{f , B} (\ZZ))^* \, .
\end{equation}
\end{defn}
We see immediately from Proposition~\ref{pr:hyperbolic} that
\begin{equation} \label{eq:contains tan}
\cK^{f , B} (\ZZ) \supseteq \tan_\xx (B)
\end{equation}
and hence 
$$\normal (\ZZ) \subseteq \normalx (B) \, .$$

In order to produce deformations, we will need to know
that the cones $\cK^{f , B} (\ZZ)$ vary semi-continuously
as $\ZZ$ varies over the torus $\exp (\xx + i \R^d)$.  
We have seen already that all of these cones contain
$\tangent_\xx (B)$.  What is needed, therefore, is an
argument showing that $\cK^{f , B} (\ZZ')$ contains
any $\uu \in \cK^{f , B} (\ZZ)$ when $\ZZ'$ is 
sufficiently close to $\ZZ$ and $\uu \notin \tangent_\xx (B)$.
In fact, not every polynomial admits a semi-continuous choice 
of tangent subcone; a counterexample is $xy + z^3$.  However,
in the case where $\xx \in \partial B$, we are able to use
strong hyperbolicity in directions $\vv \in \tangent_\xx (B)$
to prove semi-continuity even outside of $\tangent_\xx (B)$.
The main result of this section is exactly such an analogue of 
Proposition~\ref{pr:relinearization}:
\begin{thm} \label{th:semi}
Suppose that an analytic function $f$ is strongly hyperbolic 
in direction $\vv$ at the point $\zz = \xx + i \yy$.  Let 
$\ft := \homog (\ft , \zz)$.
Let $\uu \in \cK^\vv (\ft)$ be any point in the cone of hyperbolicity 
of $\ft$ containing $\vv$.  Then $f$ is strongly hyperbolic in direction
$\alpha \vv + (1 - \alpha) \uu$ for every $0 \leq \alpha \leq 1$.
\end{thm}

\begin{cor} \label{cor:semi}
~~\\[-2ex] 

$(i)$~If $B$ is a component of $\amoeba (F)^c$ and $\xx \in \partial B$,
then $\cK^{f , B} (\ZZ)$ is semi-continuous in $\ZZ$ as $\yy$ varies 
with $\ZZ = \exp (\xx + i \yy)$, meaning that $\cK^{f , B} (\ZZ) 
\subseteq \liminf_{\ZZ'} \cK^{f , B} (\ZZ')$. 

$(ii)$~If $A$ is a homogeneous polynomial and $B$ is a cone of 
hyperbolicity for $A$, then $\cK^{A,B} (\yy)$ is semi-continuous in $\yy$.  
\end{cor}

\noindent{\sc Proof:}  Pick any $\vv \in \tangent_\xx (B)$.  The 
function $f$ is strongly hyperbolic in direction $\vv$, hence
by Theorem~\ref{th:semi}, it is strongly hyperbolic at $\xx + i \yy$
in every direction $\uu \in \cK^\vv (\ft)$.  Because strong 
hyperbolicity is a neighborhood property,
it follows that for every $\yy'$ in some neighborhood of $\yy$, 
some cone of hyperbolicity of $\ft$ contains $\cK^\vv (f)$.  
All these cones contain $\vv$, hence these are the cones 
$\cK^{f , B} (\ZZ')$ (with $\ZZ' := \exp (\xx + i \yy')$), 
and hence all these cones contain $\cK^{f , B} (\ZZ)$.  
The proof in the homogeneous case is analogous, again because 
each $\At := \homog (A , \yy)$ is strongly hyperbolic in direction 
$\vv$ for any $\vv \in B$.
$\Cox$

\noindent{\sc Proof of Theorem}~\ref{th:semi}:
The proof is based on a technique of 
G{\aa}rding~\cite[Theorem~H~5.4.4]{garding-hyperbolic}
that is used in the proof of~\cite[Lemma~3.22]{ABG}.  Let $f$ be 
strongly hyperbolic at $\xx + i \yy$ in direction $\vv$ with radius 
$\ee$ and choose any $\uu \in \cK^\vv (f)$.  For the remainder of this 
argument, we assume that $\yy'$ and $\vhat'$ satisfy
$$|\yy' - \yy|  , |\vhat' - \vhat| < \frac{\ee}{2} \, ;$$
a consequence is that $f$ is strongly hyperbolic in direction
$\vhat'$ at $\xx + i \yy'$ with radius $\ee / 2$. 
For any $\bb \in \R^d$, if $s$ is purely imaginary 
with $|s| |\bb| < \ee / 2$, then the imaginary vector 
$s \bb + i (\yy' - \yy)$ will have magnitude less than $\ee$. 
By hypothesis, when $0 < t < \ee$, the function 
\begin{equation} \label{eq:bb}
s \mapsto f (\xx + i \yy' + s \bb + t \vhat')
\end{equation}
will therefore be nonzero.  

As the complex argument $s$ tends to zero, there is an expansion
$$f(\xx + i \yy + s (\alpha \vhat' + (1 - \alpha) \uu)) = 
   s^m \ft (\alpha \vhat' + (1 - \alpha) \uu) + s^{m+1} B(\alpha , s)$$
where $B$ is analytic.  The homogeneous function $\ft$
does not vanish on the convex hull of $\uu$ and the $(\ee/2)$-ball 
about $\vhat$, hence $|\ft (\alpha \vhat' + (1 - \alpha) \uu)|$ 
is uniformly bounded away from zero for $\alpha \in [0,1]$ 
and $|\vhat' - \vhat| \leq \ee/2$.  It follows that
for a sufficiently small $\delta$ (which we take also to be
less than $\ee$), the function 
$$s \mapsto f( \xx + i \yy' + s (\alpha \vhat' + (1-\alpha) \uu) 
   + t \vhat')$$ 
has exactly $m$ roots bounded in absolute value by $\delta$,
as long as $|\yy' - \yy|$ and $t$ are both bounded in
magnitude by $\delta$.  Once $2 \delta |\alpha \vhat' + 
(1 - \alpha) \uu| < \ee$ for all $0 \leq \alpha \leq 1$, then, 
taking $\bb = \alpha \vhat' + (1 - \alpha) \uu$ in~\eqref{eq:bb}, 
we see that these $m$ roots cannot be purely imaginary, and their 
real parts must therefore retain the same sign as $\alpha , \beta$
and $\yy'$ vary.  When $\alpha = 1$, these are the $m$ roots in $s$
of $f (\xx + i \yy' + (s+t) \vhat')$, so the real parts are $t$ less
than the real parts of the roots of $f (\xx + i \yy' + s \vhat')$
which are all negative by strong hyperbolicity of $f$ at 
$\xx + i \yy'$ in direction $\vhat'$.  We conclude that for
all positive real $s$ in the interval $0 < s < \delta$,
the function $f (\xx + i \yy' + s (\alpha \vhat' + (1 - \alpha) \uu))$
does not vanish, finishing the proof of strong hyperbolicity with
neighborhood size $\delta$, for any $\alpha \in [0,1]$.    $\Cox$

\begin{cor} \label{cor:section}
Let $F$ be a Laurent polynomial and $f := F \circ \exp$.
Let $\xx \in \partial B$ for some component $B$ of $\amoeba (F)^c$.
Let $\theta$ be a continuous unit section of $\cK^{f , B} 
(\exp ( \xx + i \cdot))$.  In other words, $\theta : 
(\R / (2 \pi \Z))^d \to S^{d-1}$ is continuous and $\theta (\yy) 
\in \cK^{f , B} (\exp (\xx + i \yy ))$ for each $\yy$.
Then there is some $\ee_0 > 0$ such that for all $0 < \ee < \ee_0$,
$f (\xx + i \yy + \ee \theta (\yy)) \neq 0$.
\end{cor}

\noindent{\sc Proof:} For each $\yy$, let $\ee(\yy)$ be a
radius of strong hyperbolicity for $f$ at $\xx + i \yy$
in direction $\theta (\yy)$.  Choosing a neighborhood $\nbd (\yy)$
such that $|\theta (\yy') - \theta (\yy)| < \ee (\yy) / 2$ when 
$\yy' \in \nbd (\yy)$, we see that $\ee (\yy) /2$ is a radius 
of strong hyperbolicity for $f$ at $\xx + i \yy'$ for any 
$\yy' \in \nbd (\yy)$.  Covering the compact set $(\R / (2 \pi \Z))^d$
with finitely many neighborhoods $\nbd (\yy^{(1)}) , \ldots , 
\nbd (\yy^{(n)})$, we may choose $\ee_0 = \min_n \ee (\yy^{(n)})$.
$\Cox$

\subsubsection*{Examples and counterexamples}

It is important to understand how semi-continuity may fall short 
of continuity.  This is illustrated in the following examples. 
To avoid misleading you with the pictures, we note that all of 
the upcoming figures show complex lines in $\C^2$, but that for 
obvious dimensional reasons, only the intersection with the 
$\R \times \R$ subspace is shown.
\begin{example}[cones drop down on a substratum] \label{eg:semi}
Let $F = L_1 L_2 = (3-X-2Y)(3-2X-Y)$.  This differs from
figure~\ref{fig:ghost} in that now $l_2$ also passes through $(1,1)$.
However, since $L_2$ in this example is the inversion of $L_2$
in example~\ref{eg:ghost}, the amoeba is the same as in
example~\ref{eg:ghost}.  We will see that the cone $\cK (\ZZ)$
drops discontinuously as $\ZZ \to (1,1)$, in contrast to 
example~\ref{eg:ghost}.
\begin{figure}[ht]
\hspace{2.2in} \includegraphics[scale=0.50]{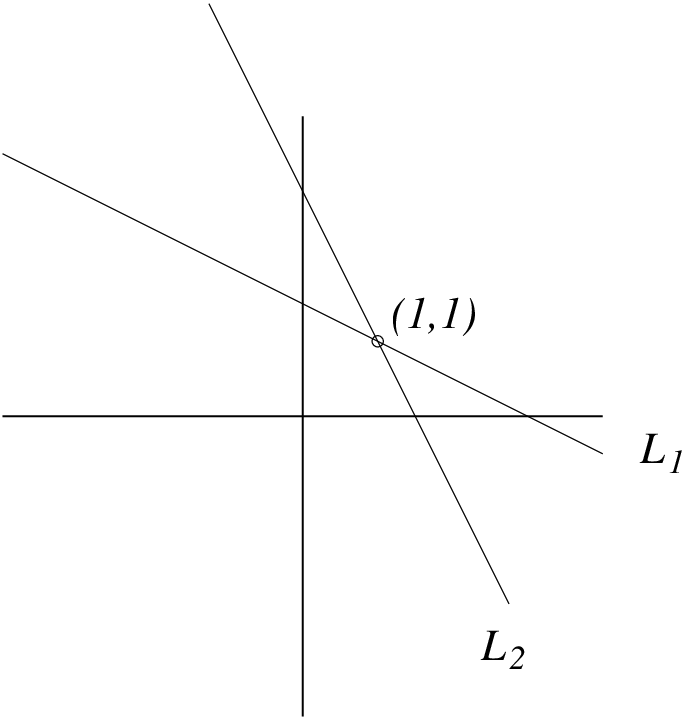}
\caption{the zero set of the function $L_1 L_2 := (3-X-2Y)(3-2X-Y)$} 
\label{fig:non-ghost}
\end{figure}
The subset of $\sing_F$ lying in $\Log^{-1} (B)$ is the union of 
two rays $\{ (Y-1) = -(X-1)/2 : X \leq 0 \} \cup 
\{ (Y-1) = -2 (X-1) : X \geq 0 \}$ with the common 
endpoint $(1,1)$.  For any point $\ZZ$ in this set other than
$(1,1)$, the cone $\cK(\ZZ)$ is equal to $\tan_{\log \ZZ} B$
which is a halfspace.  For $\ZZ = (1,1)$, the cone $\cK (\ZZ)$
is still equal to $\tan_{\log \ZZ} (B)$, but now this is 
a proper cone bounded by rays with slope $-1/2$ and $-2$.  
This cone is the intersection of the two halfspaces that are
possible values of the cone at nearby points, thus $\cK (1,1)$
is equal to the lim inf of $\cK (\ZZ)$ for nearby $\ZZ$, but
there is a discontinuity at $(1,1)$.  

Compare this to example~\ref{eg:ghost}.  Here, 
$\sing_F \cap \Log^{-1} (B)$ is the union of two rays with 
different endpoints $(1,1)$ and $(-1,-1)$ and $\cK (\ZZ)$ 
is continuous, being constant on each ray and
equal to a different halfspace on each ray.
\end{example}

The containment $\tan_\xx (B) \subseteq \cK^{f , B} (\ZZ)$ for
$\ZZ = \exp (\xx + i \yy) \in \sing_f$ may be strict.  
We will see later that this causes a headache, so we
formulate a property allowing us to bypass this trouble
in some cases.
\begin{defn} \label{def:well covered}
Say that $\xx$ is a well covered point of $\partial B$ 
if $\cK^{f , B} (\ZZ) = \tan_\xx (B)$ for some 
$\ZZ = \exp (\xx + i \yy)$.  
\end{defn}
\noindent{We} now give two examples of points that are not well covered.
\begin{example}[two lines with ghost intersection] \label{eg:ghost}
Let $F = L_1 L_2 = (3 - X - 2Y)(3 + 2X + Y)$.  The variety
$\sing_F$ is shown on the left of figure~\ref{fig:ghost}.  Its amoeba is
identical to the amoeba on the right of figure~\ref{fig:amoeba01}.
Indeed, it is the union of $\amoeba (3-X-2Y)$ and $\amoeba (3+2X+Y)$,
the latter of which is identical to $\amoeba(3-2X-Y)$ because
the amoeba of $F(-X,-Y)$ is the same as the amoeba of $F(X,Y)$.
\begin{figure}[ht] \centering
\subfigure[$\sing_F$]
{\includegraphics[scale=0.50]{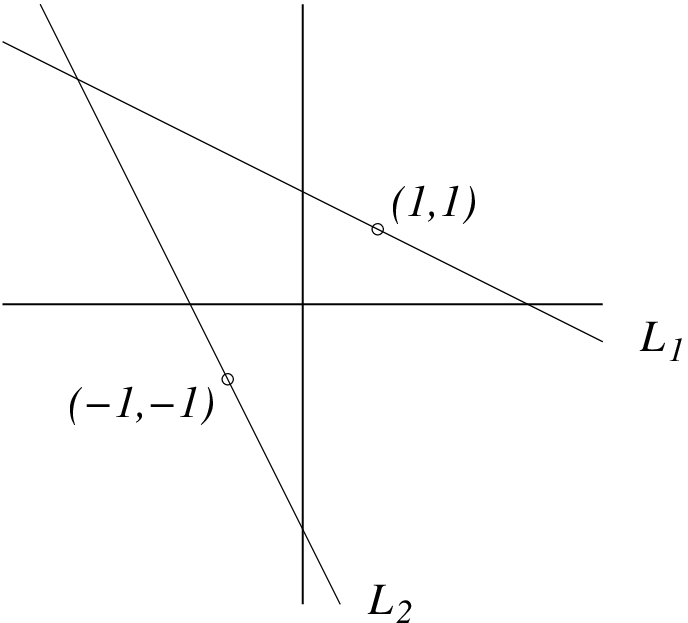}}
\hspace{0.5in}
\subfigure[a component of $\amoeba (F)$]
{\includegraphics[scale=0.50]{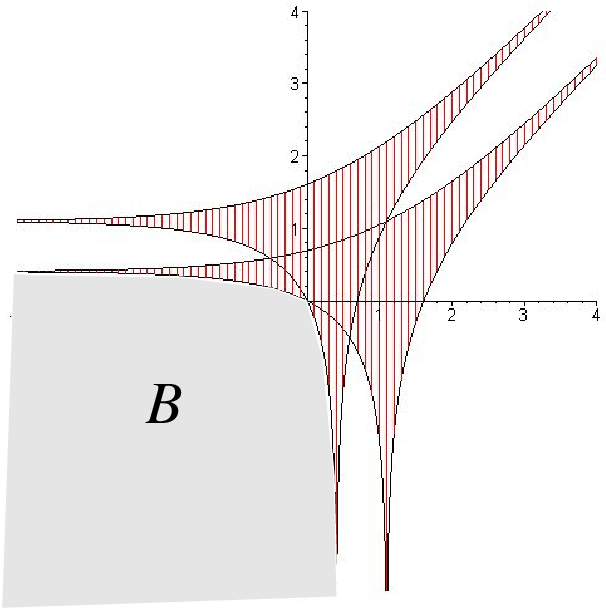}}
\caption{the zero set of $(3-X-2Y)(3+2X+Y)$ 
from example~{\protect{\ref{eg:ghost}}} and 
the OPS component}
\label{fig:ghost}
\end{figure}
The component $B$ of $\R^d \setminus \amoeba (F)$ containing 
the negative quadrant corresponds to the ordinary power series.  
An enlargement of this component is
shown on the right of figure~\ref{fig:ghost}.  
For $\xx \neq (0,0) \in \partial B$, the only point $\ZZ = 
\exp (\xx + i \yy)$ of $\sing_f$ is the real point 
$\ZZ = \pm \exp (\xx)$, the positive point being chosen for
the part of $\partial B$ in the second quadrant and the negative
point for the part of $\partial B$ in the fourth quadrant.
In either case, $\cK (\ZZ)$ is equal to the half space 
$\tan_\xx (B)$.  

On the other hand, when $\xx = (0,0)$, the linearization of 
$f$ at $\xx$ is just $\ell_1 \ell_2 := (X+2Y)(2X+Y)$.  
The zero set of which contains
the two rays forming the boundary of 
$$\tan_\xx (B) = \{ (u , v) \in \R^2 : 2 u + v < 0 \mbox{ and } 
   u + 2v < 0 \} \, . $$  
There are two points $\ZZ \in \sing_F$ in $\Log^{-1} (0,0)$,
namely $(1,1)$ and $(-1 , -1)$.  The first is in $\sing_{L_1}$
and the second is in $\sing_{L_2}$.  The cone $\cK (1,1)$
is the halfspace $\{ (u , v) \in \R^2 : u + 2 v < 0 \}$,
while the cone $\cK (-1 , -1)$ is the halfspace 
$\{ (u , v) \in \R^2 : 2 u + v < 0 \}$.  Both of these
cones strictly contain the cone $\tan_\xx (B)$.
The term ``ghost intersection'' refers to the fact that the
two curves $\Log \sing_{L_1}$ and $\Log \sing_{L_2}$ intersect
at $(0,0)$ but the lines $\sing_{L_1}$ and $\sing_{L_2}$ 
have different imaginary parts and have no intersection 
on the unit torus (though they do intersect at $(-3,3)$).  
\end{example}

Next we include an example which is the closest we can get in
two dimensions to the amoeba of a quadratic point (which can occur
only in dimensions three and higher). 

\begin{example}[critical set has large intersection with a torus] 
\label{eg:toral}
Let $F = 1 - \sqrt{1/2} (1-X) Y - X Y^2$ be the denominator for the
generating function for a one-dimensional Hadamard quantum random 
walk (see~\cite{bressler-pemantle}).  The component $B$ of
$\R^d \setminus \amoeba (F)$ corresponding to the ordinary power
series is that component of the complement of the shaded
region in figure~\ref{fig:amoeba04} which contains the 
negative quadrant.
\begin{figure}[ht] 
\hspace{1.2in} \includegraphics[scale=0.70]{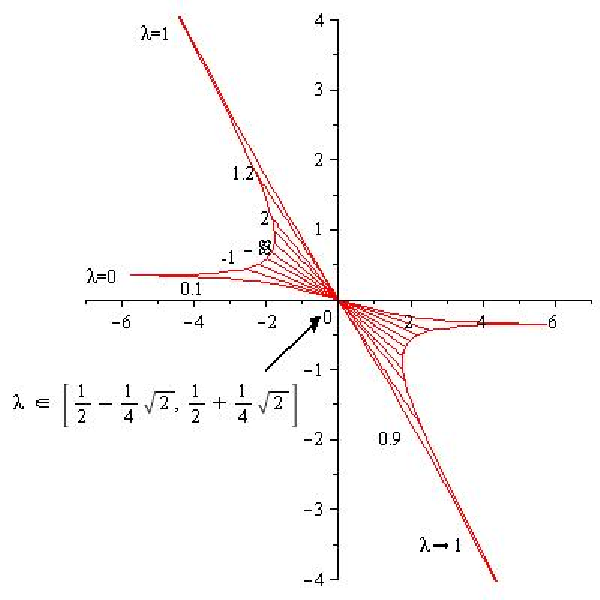}
\caption{the amoeba for $F=1-(1-X) Y/\sqrt{2}-XY^2$}
\label{fig:amoeba04}
\end{figure}

To illuminate this example a little more, observe that as we go around
the boundary of the amoeba, starting at the origin and leaving
to the northwest, the dual cone is a single projective direction
$\lambda \in \RP^1$ at every point other than the origin.  Parametrizing
$\RP^1$ by $\lambda = y/x$, we see $\lambda$ decreasing from 
$c := (1-\sqrt{1/2})/2$ to zero as the tentacle goes to infinity,
then from 0 to~$-\infty$ coming back down the other side of tentacle
and from $+\infty$ to~1 going up and out the northwest tentacle, and
so forth.  For each point of $\partial (\amoeba (F))$ other than the origin, 
there is a unique $\xx \in \R^2$ and $\yy \in T^2$ with 
$f (\xx + i \yy) = 0$; the cone $\cK (\exp (\xx + i \yy))$ 
is equal to the halfspace $\tan_\xx (B)$.  

On the other hand, when $\xx = (0,0)$, the cone $\tan_\xx (B)$ 
is bounded by the two rays $\lambda = c$ and $\lambda = 1-c$.  
This is noted in figure~\ref{fig:amoeba04} by the arrow matching
the interval $[c,1-c]$ to the single point at the origin.
It is easy to check that if $(X,Y) \in \sing_F$ then $|X| = 1$ 
if and only if $|Y| = 1$.  Thus the intersection of $\sing_F$ 
with the unit torus is the smooth topological circle parametrized
by $\{ (\phi (e^{iy}) , e^{iy}) : y \in \R / (2 \pi \Z) \}$.  

As $(x,y)$ varies over this curve, the cone $\cK(e^{ix} , e^{iy})$ 
remains a halfspace, the slope of whose normal varies smoothly 
between $(1 + \sqrt{1/2}) / 2$ and $(1 - \sqrt{1/2}) / 2$
and back.  All of these cones strictly contain $\tan_\xx (B)$.  
Thus the cone $\tan_\xx (B)$ is the intersection of the cones 
$\{ \cK (\ZZ) : \ZZ \in \sing_F \cap T^2 \}$ but these all strictly 
contain $\tan_\xx (B)$.  
\end{example}

\subsection{Critical points} \label{ss:dir}

It is time to give further examination to the role of $\xmax$.
The modulus of the term $\ZZ^{-\rr}$ in the Cauchy integral
is constant over tori, and among all tori in $\Log^{-1} (B)$ 
the infimum of $|\ZZ^{-\rr}|$ occurs on the torus 
$\Log^{-1} (\xmax)$.  This already indicates that
this torus is a good choice, but we may get some more
intuition from Morse theory.  The space $\sing$ is a
\Em{Whitney stratified space}: a disjoint union of smooth
real manifolds, called strata, that fit together nicely.  The
axioms for this may be found in Section~1.2 of part I of~\cite{GM},
along with some consequences.  We will use Morse theory only
as a guide, quoting precisely one well known result, namely
local product structure:
\begin{equation} \label{eq:locally trivial}
\parbox{4.8in}
{A point $\pp$ in a $k$-dimensional stratum $S$ of a stratified
space $\sing$ has a neighborhood in which $\sing$ is homeomorphic
to some product $S \times X$.}
\end{equation}
This is needed only for the proof of second part of
Proposition~\ref{pr:dir equiv} below, which in turn
is used only for classifying critical points when
computing examples.  According to~\cite{GM}, a proof
may be found in mimeographed notes of Mather from 1970;
it is based on Thom's Isotopy Lemma which takes up fifty 
pages of the same mimeographed notes.

A point $\ZZ \in \sing$
is a \Em{critical point} for the smooth function $h$
if $dh |_S$ vanishes at $\ZZ$, where $S$ is the stratum
containing $\ZZ$.  Goresky and MacPherson show that
in fact such points are the only possible topological
obstacles to lowering the value of $h$.  Taking 
$h = - \rhat \cdot \Log \ZZ$, we see that $(i)$ if 
there is no critical point in $\Log^{-1} (\xmax)$ 
then this torus is in fact not the best chain of integration, 
and~$(ii)$ if there is a critical point in this torus then 
we may use this fact to help us compute $\xmax$.
Because we do not give a rigorous development of stratified
Morse theory here, we give a definition of the critical set
in terms of cones of hyperbolicity, then indicate the relation
to Morse theory.
\begin{defn}[minimal critical points in direction $\rr$]
\label{def:crit}
Fix a Laurent polynomial $F$ in $d$ variables and a 
component $B$ of of $\R^d \setminus \amoeba (F)$.  
For a proper direction $\rr$, let $\xmax (\rhat)$ 
denote the unique point on $\partial B$ maximizing
$\rhat \cdot \xx$ and let $\sing_1 = \sing_1 (\rhat) =
\sing_1 (\xmax)$ denote the intersection of $\sing$
with $\Log^{-1} (\xmax)$.   Recall the notation 
\label{sing1}
$\normal (\ZZ)$ for the dual cone to the cone $\cK^{f,B} (\ZZ)$
and define the set of minimal critical points by
$$\critical (\rr) := \{ \ZZ \in \sing_1 (\rr) : \rr \in \normal (\ZZ) \} \, .$$
\label{crit}
A logarithmic version of $\critical$ is
\label{logcrit}
$$\logcrit (\rr) := \{ \yy \in \flattorus : 
   \exp (\xmax + i \yy) \in \critical (\rr) \} \, .$$
\end{defn}
\noindent{The} term ``minimal'' refers to the fact that $\Log \ZZ \in 
\partial B$ and follows the terminology of~\cite{PW1, PW2}.
\begin{pr} \label{pr:dir equiv}
Fix a Laurent polynomial $F$ in $d$ variables, let
$f := F \circ \exp$, and let $B$ be a component of 
$\R^d \setminus \amoeba (F)$.  If $\ZZ \in \sing_1 (\rr)$
is not in $\critical (\rr)$ then there is some 
$\vv \in \cK^{f,B} (\ZZ)$ with $\rhat \cdot \vv = 1$.
Conversely, if $\ZZ \in \critical (\rr)$ then $\ZZ$
is a critical point for the function $\phi := \rhat 
\cdot \log \ZZ$ on the stratified space $\sing$.
\end{pr}

\noindent{\sc Proof:} If $\ZZ \notin \critical (\rr)$ then by 
definition of the dual cone, the maximum of $\rhat \cdot \xx$ 
on $\tangent_\xx (B)$ is strictly positive.  Letting $\vv'$ 
denote a vector in $\tangent_\xx (B)$ for which $\rr \cdot \xx > 0$, 
we may take $\vv$ to be the appropriate multiple of $\vv'$.  

For the converse, suppose that $\ZZ$ is not a critical point of 
the function $\phi$ on $\sing$.  Then $\zz := \log \ZZ$ is not 
a critical point for $f := F \circ \exp$ on $\log \sing$; 
denoting $d(\phi \circ \exp)$ by $\rr$, we see, by definition 
of criticality in the stratified sense, that $\rr|_S$ is not 
identically zero, where $S$ is the stratum of $\log \sing$ 
in which $\zz$ lies.  

We claim that the linear space $T_\zz (S)$ is what~\cite{ABG}
call a \Em{lineality} for the function $\ft := \homog (f , \zz)$, 
meaning that $\ft (\ww + \ww') = \ft (\ww)$ for any
$\ww' \in T_\zz (S)$ and any $\ww \in \C^d$.  To see this,
for any $\ww \in \C^d$, let $\ww = \ww_{\parallel} + \ww_\perp$
denote the decomposition into an element $\ww_\parallel \in T_\zz (S)$
and an element in the complementary space $T_\zz (S)^\perp$.
Write $f$ as a power series $\sum c_\rr \ww_\parallel^\rr$
in $\ww_\parallel$ with coefficients that are power series in
$\ww_\perp$.  The coefficients $c_\rr (\zero)$ vanish for 
$|\rr| < m := \deg (f,\zz)$.  By~\eqref{eq:locally trivial}, the 
degree of vanishing of $f$ at any point of $S$ is the same, hence 
$c_\rr (\ww_\perp)$ vanish identically for $\rr < m$.  This implies 
that the only degree $m$ terms in the power series for $f$ near $\zz$
are those of degree $m$ in $\ww_\parallel$, which implies that
$\ft (\ww)$ depends only on $\ww_\parallel$, proving the claim. 

By Proposition~\ref{pr:hyperbolic} we know that $\ft$ is 
hyperbolic.  By~\cite[Lemma~3.52]{ABG}, the real part of the
linear space $T_\zz (S)$ is in the \Em{edge} of $\cK^{f,B} (\ZZ)$,
is invariant under translation by vectors in $T_\zz (S)$ 
meaning that such translations map  $\cK^{f,B} (\ZZ)$ into
itself.  Any real hyperplane not containing 
the edge of a cone intersects the interior of the cone.  
Applying this to the real hyperplane $\{ \xx : \rr \cdot \xx = 0 \}$,
(recall by assumption of noncriticality that this hyperplane
does not contain $T_\zz (S)$), we conclude that there is some point 
$\pp \in \cK^{f,B} (\ZZ)$ with $\rr \cdot \pp = 0$.  This implies
$Z \notin \critical (\rr)$.
$\Cox$

Showing that $\critical (\rr)$ is contained in the set 
of critical points of the logarithmic gradient enables us
to use algebraic computational methods, cf.\ the Aztec Diamond 
computations in Section~\ref{ss:aztec-appl}.  Some of this algebraic
apparatus is detailed further in~\cite{PW2,BP-residues}; for the 
present purpose, the following observations will suffice.  
When $\ZZ$ is a smooth point of $\sing_F$, the homogeneous part 
$\ft$ of $f := F \circ \exp$ is a linear map vanishing on the 
tangent space to $f$ at $\log \ZZ$.  Hence the cone of hyperbolicity 
of $\ft$ is an open halfspace, and the dual is the normal vector to this 
halfspace, which is the logarithmic normal to $\sing_F$ at $\ZZ$.
(Thus in some sense, the dual cone $\normal$ is a set-valued 
generalization of the logarithmic gradient map.)  To compute the
smooth points of $\critical (\rr)$, we observe that the gradient
of $\rr \cdot \log \ZZ$ is $(r_1 / Z_1) , \ldots , (r_d / Z_d)$.
Thus, for $\ZZ$ to be a smooth critical point, on the divisor
$\{ H_j = 0 \}$ we must have 
\begin{eqnarray} 
H_j & = & 0 \nonumber \, ; \\[2ex]
\left ( Z_1 \frac{\partial H_j}{\partial Z_1} , \ldots , 
   Z_d \frac{\partial H_j}{\partial Z_d} \right ) & \parallel & \rr \, .
\label{eq:crit eqns}
\end{eqnarray}
Similarly, for a stratum which is the transverse intersection
of $k$ smooth divisors $\{ H_j : 1 \leq j \leq k \}$ with 
\label{loggrad}
logarithmic normals $\loggrad H_j$, the equations for critical points 
in direction $\rr$ are $H_1 = \cdots = H_k = 0$ and
\begin{equation} \label{eq:span}
\rr \in \langle \loggrad H_1, \ldots , \loggrad H_k \rangle \, , 
\end{equation}
the linear span of the $k$ logarithmic gradients.
Generically, this defines a zero-dimensional variety, meaning
that the number of solutions is finite and nonzero.

For functions $f$ and $g$, define the notation
\label{oexp}
$$f = o_{\exp} (g) \Leftrightarrow |f(x)| \leq e^{-\beta x} g(x)$$
for some $\beta > 0$ and all sufficiently large $x$.
\begin{pr} \label{pr:dir}
Let $F$ be a Laurent polynomial and $\sum_\rr a_\rr \ZZ^\rr$ be a 
Laurent series for $1/F$, convergent on a domain $\Log^{-1} (B)$ 
where $B$ is a component of $\amoeba (F)$.  Let $-B^*$ denote
the negative convex dual of the set $B$.  
\begin{enumerate} \romenumi
\item If $\rr \notin -B^*$ then $a_\rr = O(e^{-\beta |\rr|})$
   for any $\beta$.
\item If $\xx \in B$ then $a_\rr = o_{\exp} (e^{-\rr \cdot \xx})$
   for all $\rr$.  
\item If $\xx \in \partial B$ but the dot product with $\rr$
is not maximized over $\overline{B}$ at $\xx$.
   then $a_\rr = o_{\exp} (e^{-\rr \cdot \xx})$.
\item If $\rr$ is proper and $\critical (\rr)$ is empty, then
   $a_\rr = o_{\exp} (e^{-\rr \cdot \xx})$.
\end{enumerate}
\end{pr} 

\noindent{\sc Proof:} The first three statements follow
directly from the integral formula~\eqref{eq:cauchy torus}
by taking $\xx \cdot \rr$ to $+\infty$ in~$(i)$ and taking
$\xx' \cdot \rr > \xx + \rr$ in~$(ii)$ and~$(iii)$.  
The fourth conclusion is an immediate consequence of 
something we will prove in Section~\ref{sec:homotopies}:
under the hypotheses, the contour of integration in~\eqref{eq:cauchy torus} 
may be deformed so that $\Real \{ - \rr \cdot \yy \} < - \rr \cdot \xx$
for every $\yy$ on the contour. 
$\Cox$

\subsection{Quadratic forms and their duals} \label{ss:quadratic}

\label{lorentz}
Let $\lorentz$ denote the standard Lorentzian quadratic
$x_1^2 - x_2^2 - \cdots - x_d^2$.  Any real quadratic form 
\label{form}
$\form$ with signature $(1,d)$ may be written as $\lorentz 
\circ M^{-1}$ for some invertible linear map $M$.  We now define 
the dual quadratic form $\dual$ in two ways.  The classical definition
is that $\dual (\rr)$ is the reciprocal of the unique critical value 
of $\form$ on the set $\rr^{(1)} := \{ \xx : \rr \cdot \xx = 1 \}$.  
It is easy to compute the dual $\lorentzd$ to $\lorentz$.  
The point $\xx$ is critical for $S|_{\rr^{(1)}}$ if and only if 
$\grad S || \rr$, that is, if and only if $\xx || 
(r_1 , -r_2 , \ldots , -r_d)$.  Thus the unique critical point 
of $S|_{\rr^{(1)}}$ is $(r_1 , -r_2 , \ldots , -r_d) / (r_1^2
- r_2^2 - \cdots - r_d^2)$ and the reciprocal of $S$ there
is $\lorentzd (\rr) := r_1^2 - r_2^2 - \cdots - r_d^2$.
In other words, $\lorentzd$ in the dual basis $\{ r_1 , \ldots , r_d \}$
looks exactly like $\lorentz$ in the original basis
$\{ x_1 , \ldots , x_d \}$.  For the second definition, note
any real quadratic form $\form$ with signature $(1,d)$ may be 
written as $\lorentz \circ M^{-1}$ for some invertible real linear 
map $M$.  Let $M^*$ denote the adjoint to $M$, that is, 
$\langle M^* \rr , \xx \rangle = \langle \rr , M \xx \rangle$; 
in our coordinates, this is 
just the transpose.  We see from the diagram below that 
$M \xx$ is a critical point for $\form |_{\rr^{(1)}}$ if and only 
if $\xx$ is a critical point for $\lorentz |_{(\lorentzd \rr)^*}$,
leading to the alternative definition $\dual (\rr) = \lorentz^* (M^* \rr)$.
 
For computation, it is helpful to compute the matrix for
the quadratic form $\dual$.  We have 
$$\form (\xx) = \lorentz (M^{-1} \xx) = \xx^T (M^{-1})^T D M^{-1} \xx$$
where $D$ is the diagonal matrix with entries $(1 , -1 , \ldots , -1)$.
Thus the matrix for $\form$ is $(M^{-1})^T D M^{-1}$.  On the
other hand, since $\dual (\rr) = \lorentzd (M^T \rr) = 
\rr M D M^T \rr^T$, we see that the matrix for $\dual$ is
$M D M^T$.  In other words, the matrices for the quadratic forms
$\form$ and $\dual$ are inverse to each other. 

Our definition of the dual quadratic is coordinate free in 
the following sense.  Let $A = S \circ M^{-1}$ as above, and 
let $\vv = (v_1 , \ldots , v_d)$ denote coordinates in which
$A$ is represented by the standard form; in other words,
$\vv = M^{-1} \xx$ and
$$A = S (M^{-1} (\xx)) = v_1^2 - v_2^2 - \cdots - v_d^2 \, .$$
Suppose that an element $L \in (\R^d)^*$ is represented
by $(\ell_1 , \ldots , \ell_d)$ in $\vv$-coordinates, that is,
$L$ maps $\sum a_j v_j$ to $\sum a_j \ell_j$.  Then
$L \xx = (\ell_1 , \ldots \ell_d) M^{-1} \xx$, that is,
$L$ is represented by the row vector $(\ell_1 , \ldots , \ell_d) M^{-1}$
with respect to the $\xx$-basis.  Computing in the $\xx$-basis,
using this row vector for $L$ and the representation $M D M^T$ 
for $A$ computed above, we have
\begin{eqnarray*}
A^* (L,L) & = & (\ell_1 , \ldots \ell_d) M^{-1} 
   \left ( M \, D \, M^T \right ) (M^{-1})^T (\ell_1 , \ldots , \ell_d)^T \\
& = & (\ell_1 , \ldots \ell_d) \, D (\ell_1 , \ldots \ell_d)^T \, .
\end{eqnarray*} 
In the $\vv$ coordinates, $A=S$ and $A^* = S^*$, whence
$A^* (L,L) = \ell_1^2 - \ell_2^2 - \cdots - \ell_d^2$
and we see that dualization indeed commutes with linear
coordinate changes.

Dual quadratics are important because they and their
partial derivatives appear in the asymptotic formulae 
for $a_\rr$ given in Theorems~\ref{th:no plane}, 
\ref{th:cone and plane} and~\ref{th:general}.  In order 
to interpret such asymptotic estimates and series, it is 
good to know the size of $A^*$ and its partial derivatives.  
It is easy to see that if $F$ is homogeneous of degree $n$
then $\partial F / \partial r_j$ is homogeneous of
degree $n-1$.  It follows that for any multi-index $\mm \in (\Z^+)^d$,
the $\mm$-partial derivative of $(A^*)^\alpha$ is homogeneous of 
degree $2 \alpha - |\mm|$ and hence that
\begin{equation} \label{eq:size}
\left ( \frac{\partial}{\partial \rr} \right )^\mm 
   [\lorentz^* (\rr)^\alpha ]
   = O \left ( |\rr|^{2 \alpha - |\mm|} \right )
\end{equation}
The upper estimate is sharp, in the sense that the left-hand side is
$\Theta (|\rr|^{2 \alpha - |\mm|})$ except on a subset of positive 
codimension where the $\mm$-partial derivative may vanish.  

\subsection{Linearizations} \label{ss:linearization}

The Fourier integral in~\eqref{eq:cauchy torus} turns
out to be much easier to evaluate if the function
$f$ in the denominator is replaced by its leading
homogenous part.  Unfortunately, if $q$ is a polynomial
with homogenous part $\qt$, then the fact that $q - \qt$ 
is of smaller order at the origin than $\qt$ does not imply 
that $q \sim \qt$, which would be necessary for a straightforward 
estimate of $q^{-1}$ by $\qt^{-1}$.  However, on any cone where 
$\qt$ does not vanish, we do have such an estimate, and in fact 
a complete asymptotic expansion of $q^{-s}$ in decreasing powers 
of $\qt$.
\begin{lem}[expansion in decreasing powers of one function] 
\label{lem:linearization}
Suppose that  $q(\xx) = \qt (\xx) + R(\xx)$, 
where $\qt$ is homogeneous of degree $h$, and $R$ is analytic
in a neighborhood of the origin with $R(\xx) = O(|\xx|^{h+1})$.
Let $K$ be any closed cone on which $\qt$ does not vanish.  Then
on some neighborhood of the origin in $K$, $q$ does not vanish
and there is an expansion
\begin{equation} \label{eq:series 1}
q(\xx)^{-s} = \sum_{n=0}^\infty \qt (\xx)^{-s-n} \left [ 
   \sum_{|\mm| \geq n(h+1)} c (\mm , n) \xx^\mm \right ] \, .
\end{equation}
Furthermore, 
\begin{equation} \label{eq:est 1}
q(\xx)^{-s} - \sum_{|\mm| - hn < N} c(\mm , n) \xx^\mm \qt (\xx)^{-s-n} 
   = O(|\xx|^{-hs+N})
\end{equation}
on $K$ as $\xx \to \zero$.  An expansion of the same type is
possible for $p(\xx) q(\xx)^{-s}$ whenever $p$ is analytic in a 
neighborhood of the origin.
\end{lem}

\noindent{\sc Proof:} Let $R(\xx) = \sum_{|\mm| \geq h+1} b(\mm) \xx^\mm$
be a power series for $R$ absolutely convergent in some ball $B_\ee$ 
centered at the origin.  Let 
$$M := \frac{\sup_{|\xx| \in B_\ee} \sum |b(\mm)| |\xx|^\mm}
   {\inf_{|\xx| \in \partial B_\ee \cap K} \qt (\xx)} \, .$$ 
Then by homogeneity,
$$\sum_\mm \frac{|b (\mm) \xx^\mm|}{|\qt (\xx)|} \leq 1/2$$
on the $\ee / (2M)$ ball.  The binomial expansion $(1+u)^{-s}
= \sum_{n \geq 0} {-s \choose n} u^n$ converges for $|u| < 1$
and in particular for $|U| = 1/2$.  Therefore, plugging in
$\sum_\mm b(\mm) \xx^\mm / \qt (\xx)$ in for $u$ yields a series
$$\left ( 1 + \frac{R(\xx)}{\qt (\xx)} \right )^{-s} = 
   \sum_{n \geq 0} {-s \choose n} 
   \left ( \sum_\mm b(\mm) \frac{\xx^\mm}{\qt (\xx)} \right )^n$$
that converges on $B_{\ee/(2M)} \cap K$.  Multiply through by
$\qt^{-s}$ to get~\eqref{eq:series 1}.   Convergence on any
neighborhood of the origin implies the estimate~\eqref{eq:est 1}. 
$\Cox$

\setcounter{equation}{0}
\section{Results} \label{sec:results}

\subsection{Cone point hypotheses and preliminary results}

We are interested in the asymptotics of the power series 
coefficients $a_\rr$ of a rational generating function $F_0$, 
in cases where there is a cone singularity and previous known
results do not apply.  Among the properties of $F_0$ discussed in
Section~\ref{sec:prelim} there are a number of hypotheses and
notations that will arise repeatedly.  So as to be able to refer
to these {\em en masse}, we state them here.
\begin{hyps}[quadratic point hypotheses] \label{hyps:cone}
~~
\begin{enumerate}
\item Let $F$ be the product $P_0 F_1^{s_1} \cdots F_p^{s_\eta}$ 
of an analytic function $P_0$ with nonzero real powers of 
Laurent polynomials $F_j$ with no common factor.  Assume without
loss of generality that $s_j \notin \Z^+$ (since otherwise we may
absorb $F_j^{s_j}$ into $P_0$).
\item Let $B$ be a component of the complement of 
$\amoeba (\prod_{j=1}^\eta F_j)$ so that $F$ has a Laurent
series expansion on $B$.
\item Let $\rr$ be a dual vector in the dual cone $-B^*$ and
assume $\rr$ is proper with $-\rr \cdot \xx$ minimized at $\xmax$.
\item Assume that $\logcrit (\rr)$ is finite and nonempty.  Let 
$\ww \in \flattorus$ be an element of $\logcrit (\rr)$ and denote
$$\zz := \xmax + i \ww \; , \, \ZZ := \exp (\zz) \, .$$
The remaining assumptions enforce a particular set of degrees
for the denominator, namely a real power of a quadratic together
with positive integer powers of smooth divisors. 
\item With $\ZZ$ fixed, we let $P$ be the product of $P_0$ with 
all $F_j$ such that $F_j (\ZZ) \neq 0$ and collect terms, writing 
$$F = \frac{P}{Q^s \prod_{j=1}^k H_j^{n_j}} \, .$$
Denote $q := Q \circ \exp$, $h_j := H_j \circ \exp$, $p := P \circ \exp$, 
and denote the homogeneous parts of $q$ and $h_j$ by $\qt := 
\homog (q , \zz)$ and $\hht_j := \homog (h_j , \zz)$. 
\item Assume that $\qt$ is an irreducible quadratic with signature
$(1,-1, \ldots -1)$ and let $M$ be a linear map such that 
$\qt = S \circ M^{-1}$.  We allow $s=0$, in which case there 
is no quadratic factor vanishing at $\zz$.
\item Assume that $\hht_j$ are linear and that $n_j$ are a positive 
integers.
\end{enumerate}
\end{hyps}
\begin{unremark} \label{rem:exp}
We lose little generality in assuming $\logcrit (\rr)$ is non-empty 
in clause~4 above, for the following reason.  If $\logcrit (\rr)$ is
empty, then part~$(iv)$ of Proposition~\ref{pr:dir} guarantees
that $|a_\rr|$ is less than $\xx^{-\rr}$ by a factor that grows
exponentially with $|\rr|$.  
\end{unremark}

In the Aztec and cube grove examples, at the point $\ZZ$ of interest,
$s=1, k=1, n_1 = 1$, in other words, the denominator of $F$ is
a product of (the first power of) a quadratic and a smooth factor.
In the QRW example, $\eta=2$ but $k=1$ (at each of the two quadratic 
points, only one of the other factors vanishes).  There are contributions 
at the quadratic points (where $s = k = n_1 = 1$) but they turn out to be 
dominated by the contributions at smooth points ($s=0$).  
In the superballot example, $\eta = 2$ with $F_1 = 1-4xz$, 
$F_2 = 1 - x - y - z + 4xyz$, $s_1 = -1/2$ and $s_2 = -1$.  
At the quadratic point, $\ZZ = (1/2 , 1/2 , 1/2)$, $F_2$ is
quadratic, $s=1$, and $n_1 = 1/2$.  In the graph polynomial
example, $\eta = 1$ and $s = \beta$.  

We extend the expansion in Lemma~\ref{lem:linearization} 
to the generality of the quadratic point hypotheses as follows.
\begin{lem}[general quadratic point expansion] \label{lem:exp}
Assume the quadratic point hypotheses.  Let $K$ be any closed cone 
on which $\qt \prod_{j=1}^k \hht_j$ is nonvanishing.  
Then there is some neighborhood of $\zero$ in $K$ 
such that for all $N \geq 1$ the following estimate holds uniformly:
\begin{equation} \label{eq:remainder}
f (\xmax + i \ww + \yy) =
   \sum_{\mm , \ell , n : |\mm| - 2 \ell - kn < N}
   c(\mm , \ell , n) \yy^\mm \qt (\yy)^{-s-\ell} \prod_{j=1}^k
   \hht_j (\yy)^{-n_j - n} + O \left ( |\yy|^{2 \ell + |\nn| + N}
   \right ) \, .
\end{equation}
The sum is finite because $c(|\mm| , \ell , n)$ vanishes
unless $|\mm| \geq 3 \ell + (k+1) n$. 
\end{lem}

\noindent{\sc Proof:} Apply Lemma~\ref{lem:linearization} once with 
$q (\xx + \cdot)$ in place of $q$, for the given value of $s$, 
and once with $\prod_{j=1}^k h_j^{n_j} (\xx + \cdot)$ in
place of $q$, setting $s=1$.  This yields two convergent 
power series.  Multiply the two series together and multiply 
as well by the power series for $p(\xx + \cdot)$.  
$\Cox$

The results in this paper can be summarized as follows.  
First, $a_\rr$ is well approximated by a sum of contributions 
indexed by $\logcrit (\rr)$, these contributions being integrals
localized near points $\xmax + i \ww$, for $\ww \in \logcrit (\rr)$.
Secondly, depending on the geometry at $\xmax + i \ww$, this
contribution is well approximated by a certain explicit function
of $\rr$.  The result giving the decomposition as a sum
is stated as Theorem~\ref{th:localize} below, with the
remaining theorems in this section giving the contributions
in various special cases.  It should be noted that 
Theorem~\ref{th:localize} is like a trade for the proverbial
``player to be named later'', in that it allows us to state 
a complete set of results even though the meaning will not be
clear until the other results have been stated.

\begin{thm}[localization] \label{th:localize}
Assume the quadratic point hypotheses and notations.  Then there 
is a conical neighborhood $\nbd$ of $\rr$ in $(\R^d)^*$ and
there are chains $\{ \generic (\ww) : \ww \in \logcrit (\rr_0) \}$ 
defined in the text surrounding Theorem~\ref{th:avoid f} in 
Section~\ref{sec:homotopies} below, such that 
\label{contrib}
\begin{equation} \label{eq:localize}
a_\rr = \sum_{\ww \in \logcrit (\rr_0)} \contrib (\ww) 
   + o_{\exp} \left ( \xmax^{-\rr} \right ) \, .
\end{equation}
The estimate is uniform when $|\rr| \to \infty$ while
remaining within $\nbd$.  The summand is defined by 
\begin{equation} \label{eq:contrib-chain}
\contrib (\ww) := \left ( \frac{1}{2 \pi i} \right )^d
\int_{\generic (\ww)} e^{-\rr \cdot \zz} \frac{p(\zz)}{q(\zz)^s 
   \prod_{j=1}^k h_j (\zz)^{n_j} } 
   \, d\zz  \, .
\end{equation}
\end{thm}
\noindent{\sc Proof:} This is an immediate consequence of
Corollary~\ref{cor:local chains} below.   $\Cox$

\subsection{Asymptotic contributions from quadratic points}

Next, we identify the contributions $\contrib (\ww)$.  In the case 
where $\ww$ is a smooth point ($s=0, k=1, n_1 = 1$), these are 
already known.  A formula involving the Hessian determinant
for a parametrization of the singular variety $\sing$ of $F$ 
was given in~\cite[Theorem~3.5]{PW1}, which was then given in 
more canonical terms in~\cite{BBBP}.
\begin{thm}[\protect{\cite{PW1,BBBP}}]
\label{th:smooth}
Assume the quadratic point hypotheses and suppose that $s=0, k=1$
and $n_1 = 1$, so $\ZZ$ is a simple pole for $F$.  Let 
$$\loggrad := \grad f (\zz) = \left ( x_1 \frac{\partial f}{\partial x_1} 
   , \ldots , x_d \frac{\partial f}{\partial x_d} \right )$$
denote the gradient of $H := H_1$ in logarithmic coordinates and let
\label{gauss}
$\gauss = \gauss (\zz)$ denote the (possibly complex) Gaussian 
curvature of $\sing_{f}$ at $\zz$.  Suppose that $\gauss \neq 0$.  
Letting $|\cdot|$ denote the euclidean norm, we have:
$$\contrib (\ww) \sim \left ( 2 \pi \, |\rr| \right )^{(1-d)/2} 
   \, \frac{p(\zz)}{\sqrt{\gauss (\zz)} |\loggrad|}
   \zz^{-\rr}$$
The estimate holds uniformly over a sufficiently small
neighborhood of $\rr$ such that: $(i)$ the quadratic point hypotheses 
are satisfied, $(ii)$ $\kappa \neq 0$, and~$(iii)$ the point
$\ZZ = \ZZ (\rr)$ varies smoothly.  The square root should be
taken as the product of the principal square roots of the eigenvalues
of the Gauss map.
$\Cox$
\end{thm}
In the case where $\ww$ is on the transverse intersection of 
smooth (local) divisors, formulae are also already
known.  There are a number of special cases, depending on the
dimension of the space, the dimension of the intersection, and
the number of intersecting divisors.  We will not need these 
results in this paper (we need only the upper bound in
Lemma~\ref{lem:big-O}) but statements may be found 
in~\cite[Theorems~3.1,~3.3,~3.6,~3.9,~3.11]{PW2} and
in~\cite{BP-residues}.
The novel results in this paper concern the case at a quadratic point, 
that is, where $s \neq 0$.  Let $\normalx (B)$ denote the dual 
to the tangent cone $\tan_\xx (B)$.  The cone $\normalx (B)$ will
have nonempty interior.  By contrast, in example~\ref{eg:ghost} 
the cone $\tan_\xx (B)$ is always a half space and $\normalx (B)$ 
is always a single ray.  
\begin{defn}[obstruction] \label{def:obstructed}
Assume the quadratic point hypotheses and notations.  
Say that $\rr$ is \Em{non-obstructed}
if $\rr$ is in the interior of $\normalx (B)$ and if
for any $\xx$ in the boundary of the cone of
hyperbolicity of $\qt$, the cone $\cK^{\qt , B} (\xx)$
contains a vector $\vv$ with $\rr \cdot \vv > 0$.  
\end{defn}

This condition is not transparent, so we pause to discuss it.  
First, note that the non-obstruction condition will turn out to be 
satisfied for all $\rr$ in the interior of $\normalx (B)$ when 
$k=0$ (locally, the denominator of $F$ is an irreducible quadratic).
To see this, recall from 
Proposition~\ref{pr:hyp properties} that the cone of hyperbolicity 
of $\qt$ is a component of its cone of positivity.  At any point
$\vv$ on the boundary of this cone, other than the origin, 
$\qt$ is smooth and hence $\homog (\qt , \vv)$ is linear, vanishing
on the tangent plane at $\vv$ to $\{ \qt = 0 \}$.  The normals to these
planes are precisely the extreme points of the cone $\normalx (B)$.
Therefore, for any $\rr$ in the interior of $\normalx (B)$,
$\rr$ is not perpendicular to the tangent plane at to $\{ \qt = 0 \}$
at any point $\vv \neq 0$, which implies that $\rr$ is non-obstructed.
An example where there are obstructed directions interior
to $\normalx$ is as follows.
\begin{example}[obstruction] \label{eg:obstructed}
Suppose the denominator of $F$ is $H := H_1 H_2 H_3 := 
(1-X)(1-Y)(1-XY)$.  Then $\hht := C x y (x+y)$.  
The cone $\tan_{(0,0)} (B)$ is the negative orthant.
The dual cone $\normalx (B)$ is the positive orthant.
The vector $\rr = (1,1)$ lies in the interior of the dual cone.
Let $\xx = (t , -t)$ for some $t \neq 0$.  Then $\cK^{\qt,B} (\xx)$
is the halfspace $\{ (x,y) : x + y < 0 \}$ and $-\rr \cdot (x,y)$
is minimized at zero on this cone.    
\end{example}

Secondly, we see that the condition of non-obstruction is not 
merely technical, but is necessary for the conclusions we wish 
to draw.  To elaborate, we would like our asymptotics to be 
uniform as $\rr$ varies over the interior of $\normalx (B)$.
Unfortunately, this is not always possible.  In the previous 
example, if $F = 1/H = 1 / [(1-X)(1-Y)(1-XY)]$, then 
$a_\rr = \min \{ r_1 , r_2 \}$.  Analytic expressions for
$a_\rr$ will not be uniform as $\rr$ approaches the diagonal.  
This is in fact because movement of the contour of integration 
in~\eqref{eq:cauchy torus} will be obstructed, requiring different 
deformations for $\rr$ in the positive quadrant on different 
sides of the diagonal.

\begin{thm}[quadratic, no other factors] \label{th:no plane}
Assume the quadratic point hypotheses~\ref{hyps:cone}, and suppose
that $k=0$, in other words, $F = P/Q^s$ with no further factors
in the denominator and $s \neq 0, -1, -2, \ldots$.  Let $c (\mm , n)$ 
be the coefficient of 
$\xx^\mm \qt (\xx)^{-1-n}$ in the expansion~\eqref{eq:series 1}.
Let $\cKd$ be any compact subcone of the interior of $\normal$.
Then, uniformly over $\rr \in \cKd$, when the Gamma functions
in the denominator are finite, there is an expansion
\begin{equation} \label{eq:s power}
\contrib (\ww) \sim
   \frac{|M|}{2^{2s-1} \pi^{d/2-1} \Gamma (s) \Gamma (s+1-d/2)}
   \ZZ^{-\rr} \sum_n \sum_{|\mm| \geq 3n} c (\mm , n) 
   (-1)^{|\mm|} \frac{\partial^\mm}{\partial \rr^\mm}
   \left ( \qt^* (\rr)^{s+n-d/2} \right ) \, .
\end{equation}
The series is asymptotic in the following sense.  For any $N$,
restricting the series to terms with $|\mm| - 2n < N$ yields
an approximation whose remainder term is $O(|\rr|^{2s-d-N})$,
all of whose terms are generically of order $|\rr|^{2s-d-N+1}$.  
If $P(\ZZ) \neq 0$ then
\begin{equation} \label{eq:leading s power}
\contrib (\ww) \sim
\frac{P(\ZZ) \, |M|}
   {2^{2s-1} \pi^{d/2-1} \Gamma (s) \Gamma (s+1-d/2)}
   \ZZ^{-\rr} \left [ \qt^* (\rr)^{s-d/2} \right ] \, .
\end{equation}
When $s+1-d/2$ is a nonpositive integer, and thus the denominator
of~\eqref{eq:s power} is infinite, the conclusion should
be understood to say that
$$\contrib (\ww) = o \left ( |\ZZ^{-\rr}| |\rr|^{-N} \right )$$
for all $N > 0$.
\end{thm}
\begin{unremark}
Comparing to equation~\ref{eq:size}, we see that the remainder 
terms are no larger than the first omitted term of~\eqref{eq:size}.  
For a true asymptotic expansion, this should be smaller than the 
last term that was not omitted, but in general there may be 
directions $\rr$ in which $(\partial^\mm / \partial \rr^\mm)
\qt^* (\rr)^{s - d/2}$ is of smaller order than $|\rr|^{2s-d-|\mm|}$.
This may occur after the first term in the expansion~\eqref{eq:s power},
though not in the leading term~\eqref{eq:leading s power}.
Also, by Theorem~\ref{th:localize}, we may be adding up
several of these formulae, thereby obtaining some cancellation.
For example in the case of the Aztec diamond, $a_\rr = 0$ 
when $\sum r_j$ is odd.  This manifests itself in the symmetry 
$F(\ZZ) = F(-\ZZ)$, and in two quadratic points at $(1,1,1)$ and 
$(-1,-1,-1)$.  Contributions from the two quadratic points will sum 
or cancel according to the parity of~$\rr$.  
\end{unremark}
As a corollary, for ease of application, we state the asymptotics 
in the three variable case for a single power of $Q$ in the denominator.
Theorem~\ref{th:no plane} is proved in Section~\ref{ss:cone only},
while Corollary~\ref{cor:simple} follows immediately.
\begin{cor} \label{cor:simple}
Assume the quadratic point hypotheses with $d=3, k=0$ and $s=1$.  Let 
$c(\mm , n)$ be the coefficients in the expansion~\ref{eq:series 1}.
Let $\cKd$ be any compact subcone of the interior of $\normal$, 
the dual cone to $\tan_\xx (B)$.  Then, uniformly over $\rr \in \cKd$, 
there is an expansion
\begin{equation} \label{eq:simplest}
\contrib (\ww) \sim \frac{|M|}{2\pi} \ZZ^{-\rr}
   \sum_{n=0}^\infty \sum_{|\mm| = n}
   c(\mm , n) \frac{\partial^\mm}{\partial \rr^\mm}
   \left [ \qt^* (\rr)^{-1/2} \right ] \, .
\end{equation}
Here, asymptotic development means that if one stops at the
term $n = N-1$, the remainder term will be $O(|\rr|^{-1-N})$,
while the last term of the summation will be of order $|\rr|^{-N}$.
In particular, if $P(\ZZ) \neq 0$ then
\begin{equation} \label{eq:leading simplest}
\contrib (\ww) \sim \frac{P(\ZZ) \, |M|}{2\pi} \ZZ^{-\rr}
   \left [ \qt^* (\rr)^{-1/2} \right ]
\end{equation}
uniformly on $\cKd$.
$\Cox$
\end{cor}
\begin{unremark}
Again, the leading term estimate~\eqref{eq:leading simplest}
is a true asymptotic estimate, while the right-hand side
of~\eqref{eq:simplest} may vanish for certain $\mm$ and $\rr$.
\end{unremark}

\subsection{The special case of a cone and a plane}
\label{ss:3.3}

Our last main result addresses the simplest case where the are
both a quadratic and a linear factor.  The case of a quadratic
along with multiple linear factors is also interesting.  We address
this in Section~\ref{ss:many planes}.  Because there are a great
number of subcases and we have no motivating examples, we do not
state here any theorems about that case, and instead describe
in Section~\ref{ss:many planes} a number of results that pertain
to this case.  In the case of a single factor of each type, in
three variables, significant simplification of the general computation 
is possible.  The remaining results concern this special case.

\begin{figure}[ht] \centering
\includegraphics[scale=0.46]{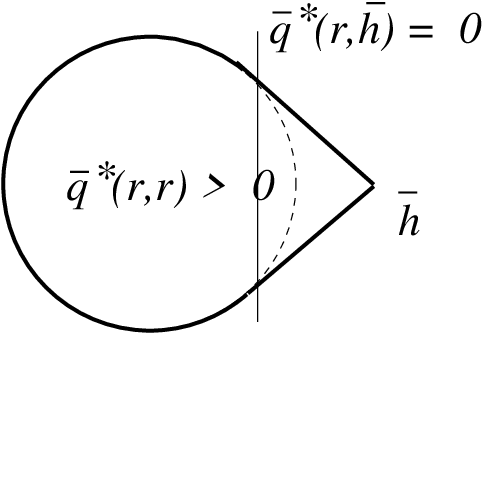} 
\caption{the cone $\normal$ depicted by its slice at $t=1$}
\label{fig:teardrop}
\end{figure}

Assume the cone-point hypotheses with $d=3, s=1$ and $k=1$.  
Because $k=1$, we drop the subscript and denote $H := H_1$.  
We assume also that the linear factor $\ell := \hht_1$ of
the homogeneous part of $(Q H) \circ \exp$ shares two real, 
distinct projective zeros with the quadratic factor $\qt$,   
and we denote these by $\alpha_1$ and $\alpha_2$.  The given
component $B$ on which the Laurent series $\sum_\rr a_\rr \ZZ^\rr$
converges is the intersection $B_1 \cap B_2$ of a component
of $\amoeba (Q)^c$ and $\amoeba (H)^c$.  By hyperbolicity,
we know that the quadratic $\qt$ is a scalar multiple of a real 
hyperbolic quadratic; multiplying by $-1$ if necessary, we may 
assume the signature to be $(1,2)$; in particular, we may write
$$\qt (\lambda \vv + \ww) = \lambda^2 - |\ww|^2$$
for some $\vv \in \tan_\xx (B_1)$ and all $\ww \in 
\tan_\xx (B_1)^\perp$.  The set $B_1$ is a cone over an ellipse $\ellipse$
and its dual $-B_1^*$ is a component of positivity of the cone
$\qt^* = 0$.  The linear function $\hht$ may be viewed as a point
of $(\R^3)^*$.  Figure~\ref{fig:teardrop} shows a plot of $\qt^* = 0$ 
and of the point $\hht$ in $(\R^3)^*$.  Also shown is the line
of points $\rr$ for which $\qt^* (\rr , \hht) = 0$.  These shapes
in the projective $(r|s|t)$-space $(\RP^2)^*$ are shown via their 
slices at $t=1$.
%
%

The assumption that $\alpha_1, \alpha_2$ are real implies 
that the point $\hht$ lies outside $B_1^*$.  The normal cone $\normalx (B)$
is the convex hull of the normal cone $B_1^*$ of $Q$ and the normal cone
$\{ \hht \}$ of $H$.  This teardrop-shaped is the entire shape shown
in Figure~\ref{fig:teardrop}.  The tangent lines to $B^*$ from $\hht$ 
intersect $B^*$ in two projective points, namely those $\rr$ for which
$\qt^* (\rr , \hht) = \qt^* (\rr , \rr) = 0$.  The non-obstructed set 
is a disjoint union $B_1^* \cup E$, where the cone $E$ is the
non-convex region $\normalx (B) \setminus B_1^*$.  Observe that
the dotted arc in figure~\ref{fig:teardrop} is obstructed and thus is
in neither $B_1^*$ nor $E$, these being the two components of the
non-obstructed set.

\label{dblres}
To state the final theorem, we must define one more quantity.
If $A$ is a homogeneous quadratic and $L$ is a linear function,
define a quantity $\dblres$ as follows.  Let $\theta$ 
denote the form $(z\, dx dy - y \, dz dx + x \, dy dz) / (A \cdot L)$.
The second iterated residue of $\theta$ is a 0-form, defined on
the two lines $\alpha_1, \alpha_2$ where $A=L=0$.  
Because $\theta$ is homogeneous of degree zero, its second
residue has a common value at any affine point in the line
$\alpha_j$.  We let $\dblres = \dblres_{A,L} (\alpha_j)$ denote
this value.  In coordinates, we have a number of formulae for
$\dblres$, one being
\begin{equation} \label{eq:dblres}
\dblres (\alpha) = \left. \frac{\disp{z}} {\disp{
   \frac{\partial A}{\partial x} \frac{\partial L}{\partial y} -
   \frac{\partial A}{\partial y} \frac{\partial L}{\partial x} }}
   \right |_{(x,y,z) \in \alpha}  \, .
\end{equation}
\begin{thm}[quadratic and one smooth factor] 
\label{th:cone and plane}
Assume the quadratic point hypotheses with $d=3, s=1$ and $k=1$ and 
let $\ell$ denote the linear factor, $\hht_1$ at the point
$\zz$.  Assume $p(\zz) \neq 0$ and assume that the two projective 
solutions $\alpha_1, \alpha_2$ to $\ell = \qt = 0$ are real and distinct,
so that the non-obstructed set $\normal$ is the union of an 
elliptic cone $B_1^*$ and a nonconvex cone $E$ as described above.

Let $\qt^*$ denote the dual to the quadratic $\qt$.  Let $\arctan$ denote 
the branch of the arctangent mapping $(0,\infty)$ to $(0,\pi/2)$, while 
mapping $(-\infty , 0)$ to $(\pi/2 , \pi)$ rather than to $(-\pi/2 , 0)$.  
Then
\begin{equation} \label{eq:atn}
\contrib (\ww) = \ZZ^{-\rr} P(\ZZ) 
   \left [ \frac{\dblres}{\pi} 
   \arctan \left ( \frac{\sqrt{\qt^* (\rr , \rr)} 
   \sqrt{- \qt^* (\ell , \ell)}} {\qt^*(\rr , \ell)} \right ) + R \right ]
\end{equation}
where the remainder term satisfies $R = O(|\rr|^{-1})$ uniformly 
as $\rr$ ranges over compact subcones of $B_1^*$.  On the other hand, 
we have the estimate
$$\contrib_\ww = \dblres P(\ZZ) \ZZ^{-\rr} + R$$
where $R = O(|\rr|^{-1})$ uniformly as $\rr$ ranges over
compact subcones of $E$.
\end{thm}

\setcounter{equation}{0}
\section{Five motivating applications} \label{sec:applications}

One feature is common to all but one of our applications, namely that 
$\zero$ is on the boundary of the amoeba of the denominator
of the generating function.  In this case, by part~$(iii)$ of
Proposition~\ref{pr:dir}, the coefficients $a_\rr$ decay 
exponentially as $|\rr| \to \infty$ in directions $\rhat$ 
for which $\sup_{\yy \in \tan_\zero (B)} \rhat \cdot \yy > 0$,
in other words for $\rr \notin \normal$, the dual cone to 
$\tan_\zero (B)$.  In such a case, the only significant
(not exponentially decaying) asymptotics are in directions
in the dual cone $(\tan_\zero (B))^*$.  We therefore restrict 
our attention in every case but the superballot example to 
$\rr \in (\tan_\zero (B))^*$, and consequently, to $\xmax = \zero$.  

\subsection{Tilings of the Aztec diamond} \label{ss:aztec-appl}

\subsubsection*{The model}

The Aztec diamond of order $t$ is a union of lattice squares in
$\Z^2$.  Its boundary is the polygon whose vertices are the pairs 
$(\pm r , \pm s)$ with $r , s \geq 1$ and $r+s = t$ or $t+1$.  
Thus the order~1 Aztec diamond consists of the four squares 
adjacent to the origin and the order~2 diamond consists of these 
together with the square centered at $(3/2 , 1/2)$ and its seven 
images under the symmetries of the lattice rooted at the origin.  
\begin{figure}[ht] \centering
\includegraphics[scale=0.78]{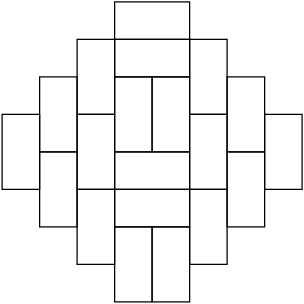}
\caption{the Aztec diamond of order~4, tiled by dominoes}
\label{fig:aztec}
\end{figure}
This was defined in~\cite{EKLP}, where questions were considered 
regarding the statistical ensemble of domino tilings of the 
Aztec diamonds.  A domino tiling of a union of lattice squares
is a representation of the region as the union of $1 \times 2$
or $2 \times 1$ lattice rectangles with disjoint interiors.
Figure~\ref{fig:aztec} shows an example of a domino tiling of
an order~4 Aztec diamond.  The set of domino tilings of the 
order~$n$ Aztec diamond has cardinality $2^{{n \choose 2}}$~\cite{EKLP}.
Let $\mu_n$ be the uniform measure on this set, that is, the 
probability measure giving each tiling a probability of 
$2^{-{n \choose 2}}$.  Limit theorems for characteristics of
$\mu_n$ have been proved, the most notable of which is the
Arctic Circle Theorem which states that outside a $(1 + \ee)$ 
enlargement of the inscribed circle the orientations of the
dominoes are converging in probability to a deterministic 
brick wall pattern, while inside a $(1-\ee)$ reduction of the
inscribed circle the measure has positive entropy~\cite{JPS}.
A new proof and a distributional limit at rescaled locations 
inside the circle were given in~\cite[Theorem~1]{aztec}.

Via an algorithm called domino shuffling~\cite{propp-shuffling},
the following generating function was found.  Color the square
centered at $(r-\half , s-\half)$ in the Aztec diamond of order $t$
black if $r+s+t$ is odd and white if $r+s+t$ is even.  A domino 
is said to be \Em{northgoing} if its white square is the 
$(0,1)$-translate of its black square.
For $r+s+t$ odd, let $p(r,s,t)$ denote the probability
that the domino covering the square centered at $(r-\half , s-\half)$
is northgoing.  The generating function~\eqref{eq:aztec},  
which we recall here, known in the 1990's to users of the 
Domino Forum and is proved, for example, in~\cite{DGIP}:
$$F := \sum p(r,s,t) X^r Y^s Z^t = \frac{Z/2}{(1 - YZ) 
   [1 - (X + X^{-1} + Y + Y^{-1}) Z/2 + Z^2]} \, . 
   \eqno(\protect\ref{eq:aztec})$$
The sum is taken over $t \geq 1$ and $-t < r , s \leq t$ 
with $|r - \half | + |s - \half | \leq t$ and $r+s+t-1 \equiv 0$ 
modulo~2.  The first results on these probabilities
were derived using bijections and other algebraic combinatorial
methods~\cite{EKLP}.  We will show that Theorem~\ref{th:cone and plane} 
implies the following asymptotic formula for $a_{rst}$.
\begin{thm}[\protect{\cite[Theorem~1]{aztec}}] \label{th:new aztec}
Let $\rt := r/t , \st := s/t$.  Let $U$ be the union of the 
two sets $\{ \rt^2 + \st^2 < 1/2 \}$ and $\{ \rt^2 + \st^2 > 1/2 \} 
\cap \{ 0 < \st < 1 - |\rt| \}$ (see figure~\ref{fig:disk-teardrop}). 
Then
\begin{equation} \label{eq:parity}
a_{rst} \sim \frac{1}{\pi} \arctan \left ( 
   \frac{\sqrt{1 - 2 \rt^2 - 2 \st^2}}{1 - 2 \st} \right ) 
\end{equation}
when $r+s+t$ is odd and zero when $r+s+t$ is even.  Here, the 
arctangent is taken to lie in $[0,\pi]$ so that it varies continuously
as $\st$ increases through $1/2$.
The asymptotic is uniform as $t \to \infty$ as long as $(\rt , \st)$
remains in a compact subset of $U$.
\end{thm}

\subsubsection*{The amoeba and normal cone}

We apply the results of Section~\ref{sec:results}.  An outline 
is as follows.  After verifying the quadratic point hypotheses,
the localization Theorem~\ref{th:localize} computes $a_\rr$
asymptotically as a finite sum
$$\sum_{\ww \in \logcrit (\rr)} \contrib (\ww) \, .$$
The point $(0,0,0)$ is on the boundary of the component $B$
and is in fact a quadratic point.  We will compute its
normal cone $\normal$ which is the teardrop shaped region 
shown in figure~\ref{fig:teardrop}.
Outside of $\normal$, the probabilities decay exponentially.  
When $\rr \in \partial \normal$ we cannot say anything, 
but for $\rr$ interior to $\normal$ we will obtain, via
Theorem~\ref{th:cone and plane}, a 2-periodic contribution at 
the critical points $\pm (1,1,1)$.  The leading term 
asymptotics~\eqref{eq:parity} will follow once we show all 
other contributions to be negligible.

Corresponding to the notation in the quadratic point hypotheses,
we write $F = P / (QH)$ where $Q := 1 - (X+X^{-1}+Y+Y^{-1}) Z/2 + Z^2$,
$H := 1-YZ$ and $P := Z/2$.  Using a computer algebra system to
compute a Gr\"obner basis for $\{ Q, Q_X, Q_Y , Q_Z \}$,
we find that $\sing_Q$ is singular precisely at $\ZZ = \pm (1,1,1)$.
Letting $q := Q \circ \exp$ and $\qt := \hom (q , \zero)$, we find 
at the point $(1,1,1)$ that $\qt (x,y,z) = z^2 - \half x^2 - \half y^2$;
the computations for the point $(-1,-1,-1)$ are analogous and
are done at the end of the discussion.  We see that near $(1,1,1)$,
$Q$ is an irreducible quadratic, while $\hht$ is linear, with 
linearization $\hht (x,y,z) = y+z$.  To specify $B$, 
observe that the components of $\amoeba (f)^c$ are intersections 
of complements of $\amoeba (Q)$ with components of the complement 
of $\amoeba (H)$.  A glance at the series $\sum p(r,s,t) X^r Y^s Z^t$ 
shows that the series is convergent for any fixed $X$ and $Y$ as 
long as $Z$ is sufficiently small.  Hence the component $B$ of 
the complement of $\amoeba (QH)$ corresponding to this series 
is the one containing $(0,0,-\lambda)$ for sufficiently large 
$\lambda$.  The amoeba of $1-YZ$ is just the line $y = - z$ in 
log space, and the component of $\amoeba (H)$ containing the ray 
$(0,0,-\lambda)$ is the halfspace $B_1 := \{y + z < 0 \}$.  
Turning to $Q$, we recall from~\cite[Chapter~6]{GKZ} that the 
components of the complement of the $\amoeba (Q)$ 
correspond to vertices of the Newton polytope $\newton (Q)$.
The Newton polytope is an octahedron with vertices $(\pm 1 , 0 , 1),
(0 , \pm 1 , 1)$ and and $(0 , 0 , 1 \pm 1)$.  There is one vertex, 
namely $(0,0,0)$, for which $(0,0,-\lambda)$ is in the interior of 
the normal cone.  Let $B_2$ to be the component of $\amoeba (Q)^c$
containing a translate of this cone.  Let $B = B_1 \cap B_2$.
This completes (1)--(2) of the quadratic point hypotheses.  

\begin{figure}[ht] \centering
\includegraphics[scale=0.38]{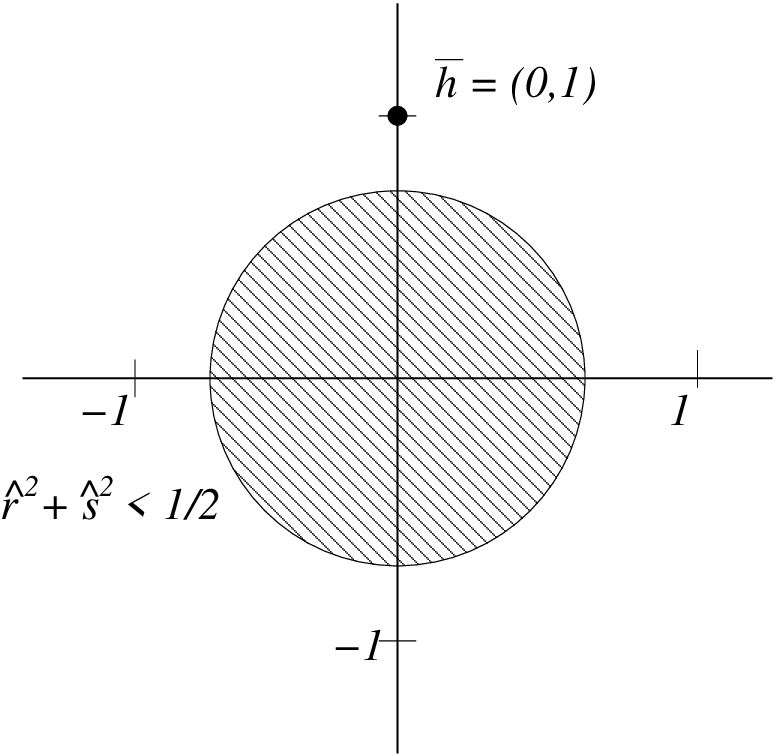} \label{sf:disk}
\caption{the disk $B_2^*$ and the point $\hht$; the region $U$
   is the interior of the convex hull of $B_2^* \cup \hht$, 
   minus the boundary of $B_2^*$}
\label{fig:disk-teardrop}
\end{figure}

As discussed at the beginning of Section~\ref{sec:applications},
in the case where $\zero \in \partial B$, we will be chiefly
interested in asymptotics in directions $\rr$ for which
$\rr \cdot \xx \leq 0$ for $\xx \in B$.  
Let us verify that $\zero \in \partial B$.  Let $Z = UXY$.  
Observe that if $X , Y , U < 1$ then the series~\eqref{eq:aztec}
is absolutely convergent.  Sending $U,X,Y$ to~1 sends $(X,Y,Z)$ 
to $(1,1,1)$ which is therefore on the boundary of the domain 
of convergence of~\eqref{eq:aztec}; hence $\zero = \log (1,1,1)$ 
is on the boundary of $\amoeba (QH)$.  We now compute 
$\normal := -(K^{\qt \cdot \hht , B})^*$.  This was done in
general in Section~\ref{ss:3.3}, so to complete the description,
we need merely to identify the dual quadratic $\qt^*$ and the
dual projective point $\hht$.  The quadratic $\qt$ is already
diagonal: $\qt = z^2 - (x^2 + y^2)/2$; hence 
$\qt^* = (1/2) t^2 - r^2 - s^2$.  Letting 
$(\rt , \st , \tat)$ be the unit vector $\rr / |\rr|$, 
we obtain the plot in figure~\ref{fig:disk-teardrop}.
The projective point $\hht$ is the point $(0|1|1)$, which in
the $t=1$ slice is given by $(0,1)$; this is  
outside the dual cone $B_1^*$, reflecting the fact that
$\qt$ and $\hht$ have two common real solutions.  

\subsubsection*{Classifying critical points}

When $\rr$ is in the interior of $\normal$, $\xmax = \zero$
and $\rr$ is obstructed only when $\rr \in \partial B_2^*$.
To finish verifying quadratic point hypotheses~(4)--(5), we need
to identify $\logcrit (\rr)$ and check that it is finite. 
As noted before, it will turn out that 
$\sum_{\ww \in \logcrit (\rr)} \contrib (\ww)$
is dominated by the contributions from $\ww = \zero$ and
$\ww = (\pi , \pi , \pi)$.  We may therefore identify the 
remaining critical points somewhat less explicitly.

Finding the critical points requires an explicit stratification
of $\sing_F$.  The coarsest Whitney stratification is as follows.
\begin{eqnarray*}
\sing_1 & := & \{ (1,1,1) , (-1,-1,-1) \} \\
\sing_2 & := & \sing_Q \cap \sing_H \setminus \sing_1 \\
\sing_3 & := & \sing_H \setminus (\sing_1 \cup \sing_2) \\
\sing_4 & := & \sing_Q \setminus (\sing_1 \cup \sing_2) 
\end{eqnarray*}
defines a Whitney stratification of $\sing_F$.  The points of
$\sing_1$ are isolated (quadratic) singularities of $\sing_Q$, 
while the remaining strata are $\sing_Q$, $\sing_H$ and their
intersection, which may be parametrized by $\{ (z^{\pm 1} , 
z^{-1} , z) : z \in \C^* \}$.  
By definition, any function is critical on a zero-dimensional
stratum, whence both points of $\sing_1$ are critical for
all $\rr \in \normal$.  Below, we will show that in fact 
$\contrib (\ww) = \Theta (1)$ for $\exp (i \ww) \in \sing_1$.  
When $\rr$ is in the interior of $\normal$, we will show
that the remaining critical points break down as follows.
\begin{equation} \label{eq:breakdown}
\begin{array}{cl}
\sing_2: & \mbox{No critical points} \\
\sing_3: & \mbox{No critical points} \\
\sing_4: & \mbox{Finitely many critical points}
\end{array}
\end{equation}
By Theorem~\ref{th:smooth}, the critical points in $\sing_4$,
which are smooth, each contribute $o(1)$ to the asymptotics,
so we will be done once we evaluate the contributions from
$\pm (1,1,1)$ and prove~\eqref{eq:breakdown}. 

Turning to the issue of counting critical points, we begin with
the easiest stratum $\sing_3$.  Recall from~\eqref{eq:crit eqns}
that on the smooth stratum $\sing_H$, the point $\ZZ$ is critical 
if and only if $\loggrad \ZZ$ is parallel to $\rr$.  The logarithmic
gradient of $H$ the constant vector $(0,1,1)$, which is on the
boundary of $\normal$, whence $\sing_3$ contains no critical points 
interior to $\normal$.  To compute critical points on $\sing_2$, 
we evaluate $\loggrad Q (z^{\pm 1} , z^{-1} , z)$ and find that
independent of $z$, we always obtain the projective point 
$(\mp \half , \half , 1)$.  This shows that $\sing_Q$ intersects
$\sing_H$ transversely, and by~\eqref{eq:span}, that $\sing_2$ 
produces critical points only when $\rr$ is in the union of two
projective lines, one joining $(0,1)$ to $(\half , \half)$ and
the other joining $(0,1)$ to $(-1/2,1/2)$.  This union does not
intersect the interior of $\normal$.

To solve for critical points in $\sing_4$, fix $\rr = (r,s,t)$ 
and solve the equations~\eqref{eq:crit eqns}:
$Q = 0$, $t X Q_X - r Z Q_Z = 0$ and $t Y Q_Y - s Z Q_Z = 0$. 
Multiplying each of these by $2xy$ clears denominators and
allows us to use a computer algebra system to compute a
Gr\"obner basis for the solution.  With lexicographic
term order ${\tt plex} (x,y,z)$, the almost-elimination polynomial 
for $z$ is $yz(1+z)^2(1-z)^2$ times a quadratic polynomial in $r,s,t$ 
and $z^2$:
\begin{eqnarray*}
&& ({r}^{4}-2\,{r}^{2}{s}^{2}-2\,{t}^{2}{r}^{2}+{s}^{4}-2\,{t}^{2}{s}^{2}
+{t}^{4})+ \left( -2\,{s}^{4}+4\,{r}^{2}{s}^{2}-2\,{r}^{4}-4\,{t}^{2}{s
}^{2}-4\,{t}^{2}{r}^{2}+2\,{t}^{4} \right) {z}^{2} \\
& + & \left( {r}^{4}-2\,
{r}^{2}{s}^{2}-2\,{t}^{2}{r}^{2}+{s}^{4}-2\,{t}^{2}{s}^{2}+{t}^{4}
 \right) {z}^{4}
\end{eqnarray*}
It is easy to check that for every $r,s,t$, this polynomial is 
not identically zero, hence there are only finitely many solutions.
The basis contains a polynomial in $y$ and $z$ (over $\C (r,s,t)$)
that is linear and non-constant in $y$, implying that for each $z$ 
there is at most one $y$.  The same is true for $x$ if we use the 
term order ${\tt plex} (y,x,z)$.  It follows that there are finitely 
many critical points in $\sing_4$ for each $\rr$.  Summing up, 
have verified (4)--(5) of the quadratic point hypotheses for $\rr$ 
interior to $\normal$.

\subsubsection*{Computing the estimate}

The computations for $\ww = (\pi , \pi , \pi)$ (hence $\ZZ = (-1, -1, -1)$)
are almost identical to those for $\ww = \zero$ and $\ZZ = (1,1,1)$.
We do the latter computation and indicate changes needed to
do the former.  We observe also that $Q$ and $H$ are invariant
under $(X,Y,Z) \mapsto (-X,-Y,-Z)$, while the numerator, $Z/2$, 
is odd; this corresponds to the parity constraint of $p(r,s,t)$ 
vanishing when $r+s+t$ is even.

The quadratic point hypotheses have now been spelled out and verified.
Let $\ww := \zero$ and $\ZZ := (1,1,1)$.  To check that we are 
in the case covered by Theorem~\ref{th:cone and plane}, we need 
to check that the two projective solutions to $\qt = \hht = 0$ 
are real and distinct.  This is easy: plugging in $y = -z$, we 
get $z^2 - \half z^2 - \half x^2 = 0$ which has the two real
solutions $y = -z = \pm x$.  

\begin{figure}[ht] \centering
\includegraphics[scale=0.60]{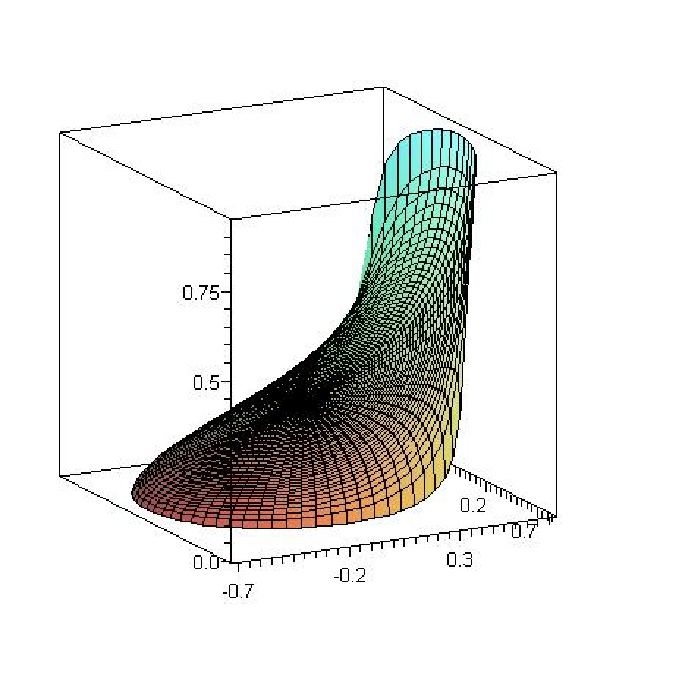} 
\caption{asymptotics for northgoing probabilities in the Aztec Diamond}
\label{fig:aztec plot}
\end{figure}

The quantity $\qt^* (\rr , \rr)$ is the quadratic that is positive 
on the interior of the disk, reaching a maximum of~1 at $(0,0,1)$ 
and vanishing on the boundary of the disk.  In coordinates, it
is given by 
$$\qt^* (\rr , \rr) = t^2 - 2s^2 - 2r^2 \, .$$
The quantity $\qt^*(\rr , \hht)$ is equal to $t-2s$.  
This vanishes on the line shown in figure~\ref{fig:teardrop}.  
The branch of the arctangent chosen in the conclusion of 
Theorem~\ref{th:cone and plane} varies continuously through 
$\pi / 2$ as $\qt^* (-\rr , \hht)$ varies through zero and the argument 
of the arctangent passes through $\pm \infty$.  The arctangent goes to 
zero where $\qt^* (\rr , \rr) = 0$ and $\qt^* (\rr , \hht) > 0$ (the part
of the boundary of the disk to the left of the vertical line) and
to $\pi$ where $\qt^* (\rr , \rr) = 0$ and $\qt^* (\rr , \hht) < 0$ (the
part of the boundary of the disk to the right of the line).  The
residue $\dblres$ is immediately computed from~\eqref{eq:dblres}
and is equal to~1.  Finally, we have $P(\zz) = 1/2$ and $\qt^*(\hht,\hht)
= -1$.  Thus, as $\rr$ varies over the interior of
the projective disk $B_2^*$ we have 
$$\contrib (\zero) \sim \frac{1}{2 \pi} \arctan \left ( 
   \frac{\sqrt{\qt^*(\rr , \rr)}}{\qt^*(\rr , \hht)} \right )
   =  \frac{1}{2 \pi} \arctan \left ( 
   \frac{\sqrt{t^2 - 2r^2 - 2s^2}}{t - 2s} \right ) \, .$$
The computation for $\ww = (\pi , \pi , \pi)$ is entirely 
analogous, leading to the same contribution but with an extra
sign factor of $(-1)^{i+j+n+1}$.  We have already shown that all
other contributions are of order $O(|\rr|^{-1})$.  Therefore, we may
sum these results to finish the proof of Theorem~\ref{th:new aztec}.
A plot of this function is shown in figure~\ref{fig:aztec plot}; 
see~\cite[Figure~2]{aztec} for a contour plot of the same function.
$\Cox$

\subsection{Cube groves} \label{ss:cube-grove-appl}

\subsubsection*{The model}

After~\cite{carroll-speyer,cube-grove}, we define a collection
of lattice subgraphs known as \Em{cube groves}.  Let $L_n$ be the
triangular lattice of order $n \geq 0$, by which we mean the 
set of all triples of nonnegative integers $(r,s,t) \in (\Z^+)^3$
such that $r+s+t = n$ with edges between nearest neighbors (thus the
degree of an interior vertex is 6).  We depict this in the plane
as a triangle with $n+1$ vertices in the top (zeroth) row, and so 
on down to 1 vertex in the $n^{th}$ row.  

The cube groves of order $n$ are a subset $C_n$ of the 
subgraphs of $L_n$.  The set $C_n$ has a description 
where one begins with the unique cube grove of order 
zero, then produces sequentially groves of orders 
$1, 2, \ldots , n$, each produced from the previous by
a ``shuffle'' which injects some information in a manner similar
to the domino shuffling used by~\cite{aztec} in studying
and enumerating domino tilings of the Aztec diamond.  The
set $C_n$ has other, static definitions in terms of 
graphs that look like stacks of cubes and in terms of 
graphical realization of certain terms of generating functions
(see~\cite{cube-grove}), but here we will take the
shuffling procedure to define the set $C_n$ of order $n$ cube groves.

Define $C_0$ to be the singleton whose element is the one-point 
graph.  If $T$ is a downward-pointing triangular face 
of $L_n$, let $T'$ be the rotation of $T$ by $180^\circ$ 
about its center.  The union of the vertices of the triangles
$T'$ is a translation of the graph $L_{n+1}$, provided that
one adds in the three corner vertices of $L_{n+1}$.  The
edge sets of the triangles $T'$ are disjoint and their
union is the edge set of $L_{n+1}$, provided that one adds
in the six edges adjacent to corner vertices.

\begin{figure}[ht]
\hspace{1.2in} \includegraphics[scale=0.58]{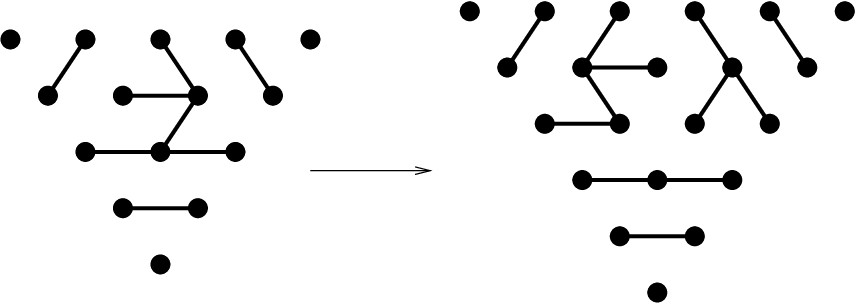}
\caption{an order 4 cube grove, shuffled to become an order 5 grove}
\label{fig:grove-shuffle}
\end{figure}
Given a cube grove $G \in L_n$ and a downward-pointing 
triangular face $T$ of $L_n$, let $G(T)$ be the collection 
of graphs on $T'$ that have: no edges if $G$ has
two edges in $T$; one edge if $G$ has one edge $e \in T$,
in which case the edge of $T'$ must be the edge of $T'$
parallel to $e$; two edges if $G$ has no edges in $T$, in
which case any two of the three edges of $T'$ will do.
Let $C(G)$ be the direct sum of $G(T)$ as $T$ varies over
downward-pointing triangular faces of $L_n$.  That is, 
choose an element of $G(T)$ for each $T$ and take the union 
of these.  Figure~\ref{fig:grove-shuffle} shows an order~4
grove, $G$ and one of the 27 elements of $C(G)$.  Finally, 
let $C_{n+1}$ be the (disjoint) union of the collections 
$C(G)$ as $G$ runs over $L_n$.  

Looking at a picture of a uniformly chosen random cube grove of order
100, one sees regions of order and disorder similar to those of
the Aztec Diamond.
\begin{figure}[ht]
\hspace{1.2in} \includegraphics[scale=0.58]{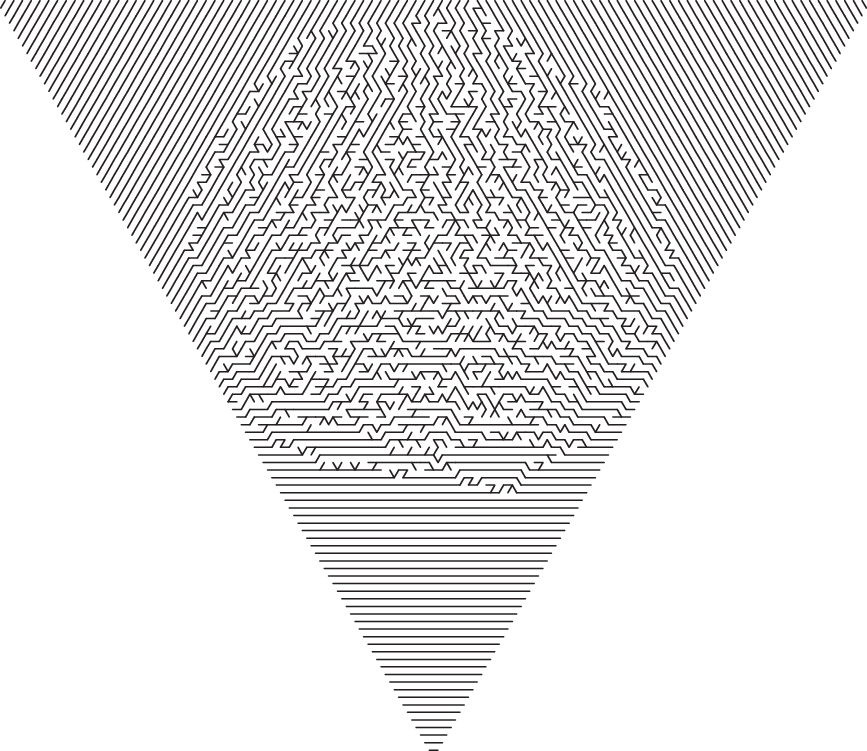}
\caption{a random cube grove of size 100}
\label{fig:grove100}
\end{figure}
Let $p_n(i,j)$ be the probability that the horizontal
edge with barycentric coordinates $(i,j,n-i-j)$ is present
in a uniformly chosen cube grove of order $n$.  The creation
rates $E_n (i,j)$ may be defined in terms of the shuffling procedure
but in this case they satisfy the simple relation $E_{n-1} (i,j)
= \frac{3}{2} (p_n(i,j) - p_{n-1} (i,j))$~\cite[Theorem~2]{cube-grove}.
We recall here the explicit generating function~\eqref{eq:cube-grove},
which is derived in~\cite[Section~2.2]{cube-grove}:  
$$F (X,Y,Z) = \frac{2 Z^2}{(1-Z) 
   (3 + 3XYZ - (X+Y+Z+XY+XZ+YZ))}
   := \frac{2Z^2}{H Q} \, .$$
It is quick to verify that longer factor, $Q$, in the denominator has 
a quadratic cone singularity and that $F$ therefore is singular 
on the union of a quadratic cone with a smooth surface.  The real part 
of this is pictured in figure~\ref{fig:cubeplot}.  Application
of Theorem~\ref{th:cone and plane} will yield the following result.
\begin{thm} \label{th:new cubegrove}
The quantity $p_t (r,s)$, which is the coefficient $a_{rst}$ of 
$X^r Y^s Z^t$ in~\eqref{eq:cube-grove}, is given asymptotically by
\begin{equation} \label{eq:new cubegrove}
\frac{1}{\pi} \arctan \left ( 
   \frac{\sqrt{2(rs+rt+st) - (r^2+s^2+t^2)}}{r+s-t} \right ) 
\end{equation}
where the arctangent is taken in $(0,\pi)$ so that as we cross
the line $t=r+s$ the arctangent varies continuously across $\pi/2$.
\end{thm}

\begin{figure}[ht]
\hspace{1.2in} \includegraphics[scale=0.58]{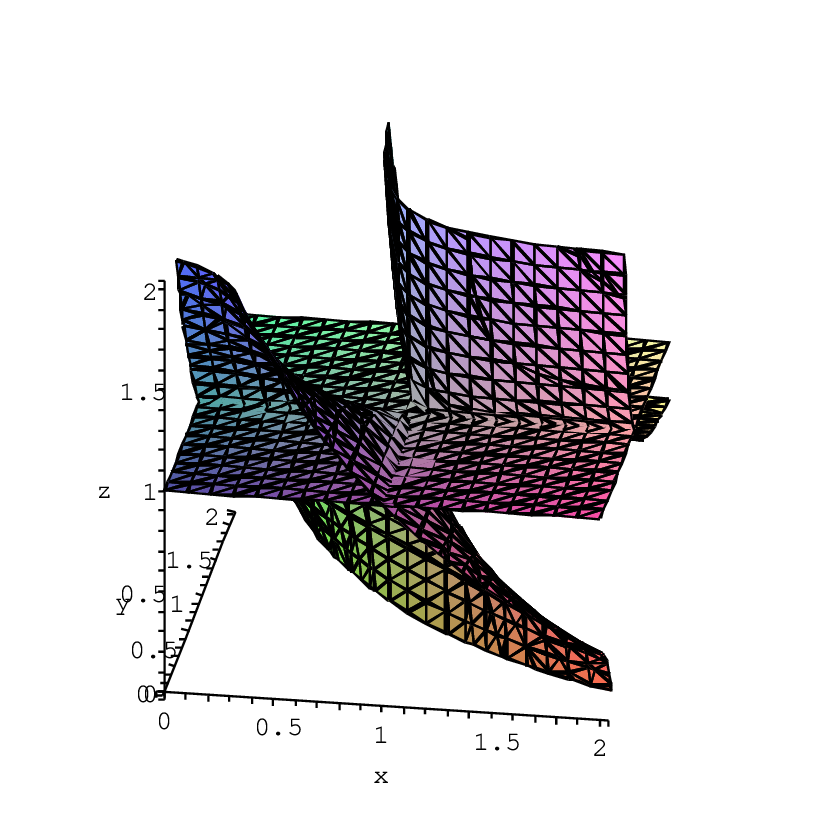}
\caption{the pole variety of the cube grove generating function}
\label{fig:cubeplot}
\end{figure}

\subsubsection*{The amoeba and the normal cone}

\begin{figure}[ht]
\centering
\includegraphics[scale=0.58]{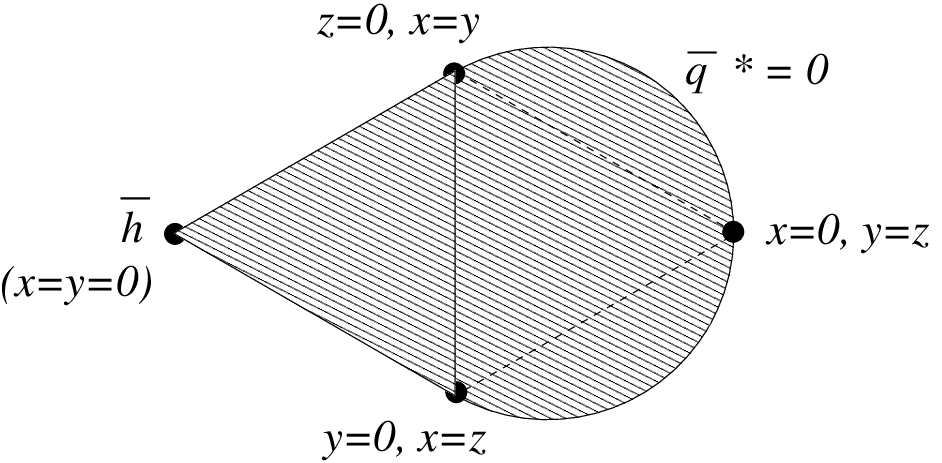}
\caption{the dual cone in symmetrized coordinates}
\label{fig:120}
\end{figure}

All multi-indices in the generating function are nonnegative, so it
is an ordinary generating function and $B$ will be the component
of $\amoeba (F)^c$ containing the negative orthant.  Again, we 
are chiefly interested in directions $\rr$ for which $\xmax (\rr)
= \zero$, these being the directions of non-exponential decay.  
The polynomial $Q$ has a single quadratic point at $(1,1,1)$.  To compute
$\tan_{\xmax} (B)$, we intersect $B_1 := \{ (x,y,z) : z < 0 \}$
with the cone $B_2$ of hyperbolicity of $Q$ at $(1,1,1)$.  Changing
to exponential coordinates via $q := Q \circ \exp$ and computing the
leading homogeneous term gives
\begin{eqnarray*}
q(x,y,z) & = & \qt (x,y,z) + O(|\zz|^3) \\[1ex]
\mbox{ where } \;\;\; \qt (x,y,z) & := & 2xy + 2xz + 2yz \, .
\end{eqnarray*}
It follows that $B_2$ is the cone containing the negative orthant
and bounded by $\{ \qt = 0 \}$.  The dual quadratic is represented
by the matrix
$$\left [ \begin{array}{ccc} 0 & 1 & 1 \\ 1 & 0 & 1 \\ 1 & 1 & 0 
   \end{array} \right ]^{-1} = \half 
   \left [ \begin{array}{ccc} -1 & 1 & 1 \\ 1 & -1 & 1 \\ 1 & 1 & -1
   \end{array} \right ] \, ,$$
hence 
$$\qt^*(r,s,t) = rs + rt + st - \half (r^2+s^2+t^2) \, .$$
The dual cone $-B_1^*$ is the subcone of the positive orthant
bounded by $\qt^* = 0$.  Again, the point $\hht$, which is equal
to $(0,0,1)$ in the $(r,s,t)$ coordinates, lies outside this cone,
and again the solutions to $\qt = \hht = 0$ are the solutions
to $z = 0 = xy$ which are two distinct projective points,
namely the $x$-axis and the $y$-axis.  We could again depict 
this by the slice through $t=1$, viewing $B_2^*$ as the interior
of a parabola in the first quadrant, opening in the northeast direction 
and tangent to the axes at $(1,0)$ and $(0,1)$, with $\hht$ at $(0,0)$.
It is easier to see what is going on if we change coordinates to
$\uu := (1,1,1)/\sqrt{3}$, letting $U^\perp$ denote the complementary
space.  For $\rr = (\rr \cdot \uu) \uu + \rr^\perp$ with $\rr^\perp 
\in U^\perp$, we then have $|\rr|^2 = (\rr \cdot \uu)^2 + |\rr^\perp|^2$, 
whence 
$$|\rr^\perp|^2 = r^2 + s^2 + t^2 - \frac{1}{3} (r+s+t)^2  = 
   \frac{1}{3} (r^2 + s^2 + t^2) - \frac{2}{3} \qt^*(r,s,t) \, .$$
Thus $\qt^* (r,s,t) = 0$ when $|\rr|^2 = 3 |\rr^\perp|^2$,
or equivalently, $|\rr \cdot \uu|^2 = 2 |\rr^\perp|^2$.
Viewing projective space via the slice $\rr \cdot \uu = 1$, we
see that $\qt^*$ vanishes on the circle centered at the origin
of radius $\sqrt{1/2}$.  The projective point $\hht = (0,0,1)$ 
intersects the slice $|\rr \cdot \uu| = 1$ at $(0,0,\sqrt{3})$,
whose projection to $U^\perp$ has squared norm $2$.  This is 
pictured in figure~\ref{fig:120}.  In these coordinates, the 
only difference between this figure and that for the Aztec Diamond 
is that the distance from the origin to the point $\hht$ is 
twice the radius of the circle, rather than $\sqrt{2}$ times
the radius, and the tangents subtend an arc of $120^\circ$
rather than $90^\circ$.  

\subsubsection*{Classifying the critical points}

The stratification is similar to that for the Aztec Diamond 
generating function.  There is just one singular point of $Q$, 
namely $(1,1,1)$.  This is on $\sing_H$ as well.  The surfaces 
$\sing_H$ and $\sing_Q$ intersect in the set $\{ x = z = 1 \}
\cup \{ y = z = 1 \}$, which is smooth away from $(1,1,1)$,  
leading to the following stratification.
\begin{eqnarray*}
\sing_1 & := & \{ (1,1,1) \} \\
\sing_2 & := & \sing_Q \cap \sing_H \setminus \sing_1 \\
\sing_3 & := & \sing_H \setminus (\sing_1 \cup \sing_2) \\
\sing_4 & := & \sing_Q \setminus (\sing_1 \cup \sing_2) 
\end{eqnarray*}
The logarithmic gradient of $H$ is parallel to $(0,0,1)$,
which is not in $\normal$, so for $\rr \in \normal$ there 
are never any critical points on $\sing_3$.  There are 
critical points on $\sing_4$, but we verify as before that
there are only finitel many.  They are smooth, so by
Theorem~\ref{th:smooth}, their contributions are $o(1)$.
On $\sing_2$, the logarithmic gradient of $Q$ is parallel
to either $(1,0,1)$ or $(0,1,1)$.  The logarithmic gradient
of $H$ is in the $t$ direction, so the span of the two
logarithmic gradients is either the $r$-$t$ plane or the 
$s$-$t$ plane.  Neither of these planes intersects the
interior of $\normal$ (the planes are tangent to $\normal$ 
at the projective points $(1,0,1)$ and $(0,1,1)$ respectively).
The for $\rr$ interior to $\normal$, there are no contributions
from $\sing_2$; it remains to compute the contribution from $\sing_1$.

\subsubsection*{Computing the estimate}

Completing the computation as in the Aztec case, we evaluate
$\dblres$ using~\eqref{eq:dblres} but switching the roles of
$x$ and $z$ because only the $z$-derivative of $\hht$ is non-vanishing.
This gives $\dblres = \half$.  With $P = 2Z^2$, and only one contributing
point $\ww = (0,0,0)$, we have $\qt^* (\rr , \hht) = (r+s-t)/2$
and $\qt^* (\hht,\hht) = -1/2$, whence Theorem~\ref{th:cone and plane} 
gives
\begin{eqnarray*}
a_\rr & \sim & \frac{1}{\pi} \arctan \left ( 
   \frac{\sqrt{\half \qt^*(\rr , \rr)}}{(r+s-t)/2} \right ) \\[1ex]
& = & \frac{1}{\pi} \arctan \left ( 
   \frac{\sqrt{2(rs+rt+st) - (r^2+s^2+t^2)}}{r+s-t} \right ) \, ,
\end{eqnarray*}
finishing the proof of Theorem~\ref{th:new cubegrove}.
$\Cox$

\subsection{Two-dimensional quantum random walk} \label{ss:quantum-appl}

\subsubsection*{The model}

We begin with a brief review on one-dimensional quantum random 
walk (QRW).  In the classical simple random walk, the law at 
time $n$ is a probability measure on $\Z$ and the evolution operator 
on this law is $(1/2) \sigma_+ + (1/2) \sigma_-$, where $\sigma_+$
is the right-shift operator $\sigma_+ \mu (n) = \mu (n-1)$ and
$\sigma_-$ is the left-shift operator $\sigma_- \mu (n) = \mu (n+1)$.
In the quantum world, the law at time $n$ is given by the
values of $|\psi (n)|^2$ where the wave function $\psi$ is
not a positive unit vector in $L^1 (\R)$ but rather a unit
vector in $L^2 (\C)$.  Evolution operators must be unitary.
While the shifts $\sigma_\pm$ are unitary, linear combinations
of these such as $(1/2) \sigma_+ + (1/2) \sigma_-$ are not.

An idea for constructing a quantum simple random walk, 
apparently due to~\cite{QRW-meyer}, is to enlarge the
space to $\space := \Z \times \{ U,D \}$, adding a hidden 
``spin'' variable.  To take a step of the random walk, 
first the spin is randomized, then all particles with
spin up move one step right and all particles with spin
down move one step left.  A number of choices are available
for the operator that executes the evolution of spins. 
One common choice is the Hadamard coin-flip, $B := \disp{ \frac{1}{\sqrt{2}}
\left ( \begin{array}{cc} 1 & 1 \\ 1 & -1 \end{array} \right )}$.
The terminology reflects the fact that the matrix is a multiple
of an orthogonal matrix with $\pm 1$ entries, these being known 
as Hadamard matrices.
Under this operator, either state $(1,0)$ or $(0,1)$ becomes
an equal mix of U and D states.  Let $A$ be the operator which
maps state $(n,U)$ to $(n+1,U)$ and $(n,D)$ to $(n-1,D)$.
If we begin in state $(0,U)$, then do $\op := A \circ B$, the 
particle is in an equal mix of states $(1,U)$ and $(-1,D)$.
If we measure the position, we will have executed a step of QRW.

The $n$-step simple random walk is defined to be the operator $\op^n 
:= (AB)^n$.  If the position if this is measured at time $n$, 
the probability of being at position $k$ is $|\op^n (k,U)|^2
+ |\op^n (k,D)|^2$.   
Since no measurement is made until time $n$, the various possible
coin-flips and movements interfere, both positively and negatively,
and the result is somewhat complicated.  The analyses 
in~\cite{QRW-line,QRW-one-dim} show that unlike classical
simple random walk, QRW spreads out linearly,
with location distributed over the interval $[-n/\sqrt{2} , n/\sqrt{2}]$;
see also the review article~\cite{QRW-kempe}.  

To define a two-dimensional QRW, we need a four-fold auxiliary
state.  Denote these four states by $\{ N , S , E , W \}$.  Any
$4 \times 4$ unitary matrix may be used for the quantum coin-flip.
The Hadamard matrix 
$$U := \frac{1}{2} \left ( \begin{array}{cccc} 1 & -1 & -1 & -1 
   \\ -1 & 1 & -1 & -1 \\ -1 & -1 & 1 & -1 \\ -1 & -1 & -1 & 1 
   \end{array} \right )$$
is known (http://www.santafe.edu/~moore/gallery.html) as the Hadamard 
quantum coin-flip.  N.B.: This is different from the Hadamard QRW 
in~\cite{grimmett-QRW}, which also uses a Hadamard matrix, namely
the tensor product of two copies of the one-dimensional Hadamard matrix.
A step of the two-dimensional Hadamard QRW is the product $AU$ where 
$A$ maps $((r,s) , N)$ to $((r,s+1) , N)$, and so forth.  The
following generating function for the probability amplitudes 
of a QRW in any dimension with any quantum coin-flip matrix 
is given in~\cite[Proposition~3.1]{bressler-pemantle}.  Let
$M$ be obtained from $U$ by multiplying the first row by $XZ$, 
the second by $YZ$, the third by $X^{-1}Z$ and the fourth 
by $Y^{-1}Z$.  We consider the rows and columns of $M$ as indexed
by the ordered quadruple $(E,N,W,S)$.  Then an entry of $M^t$ 
such as $(M^t)_{N,E}$ counts the number of $t$-step paths from 
$N$ to $E$, weighted by $M$: the $X^r Y^s Z^t$ coefficient of 
this is the wave function at position $((r,s) , N)$ and time $t$ 
starting from $((0,0) , N)$.  Summing in $t$ shows that the 
components of $(I - M)^{-1} = \sum M(r,s,t) X^r Y^s Z^t$
are the generating functions for the wave function at all 
positions and times: each $M(r,s,t)$ is a matrix, whose
$(\xi,\eta)$-entry is the generating function
$\sum_{r,s,t} c(\xi , \eta ; r,s,t) X^r Y^s Z^t$ where $c(r,s,t)$ 
is the probability amplitude, starting from state $((0,0) , \xi)$ 
at time zero, of being in state $((r,s) , \eta)$ at time $t$.

The entries of $(I-M)^{-1}$ have denominator $(1-Z^2) Q$ where
$$Q := 1 - 2 \frac{X + X^{-1} + Y + Y^{-1}}{4} Z + Z^2 \, .$$
We recognize the same polynomial factor that occurred in the 
Aztec denominator.  We know of no reason for this coincidence.
This particular QRW is somewhat special, both because of 
the occurrence of a quadratic point and because the denominator
is reducible.   The numerators in the first row are half of 
the following.  
\begin{eqnarray*}
P_1 & = & 2 - (Y + Y^{-1} + X^{-1}) Z + Z^3 \\
P_2 & = & XZ - (1 + Y^{-1}X) Z^2 + Y^{-1} Z^3 \\
P_3 & = & XZ - (YX + Y^{-1} X) Z^2 + X Z^3 \\
P_4 & = & XZ - (1 + XY) Z^2 + Y Z^3
\end{eqnarray*}
The chiralities $\{ N,E,S,W \}$ and the location $(i,j)$ are 
simultaneously measurable, so the probability of a QRW started at 
$((0,0),N)$ to be found at $(r,s)$ at time $t$ is the sum over
$1 \leq k \leq 4$ of $|\CC_{r,s,t} P_k|^2$, the squared moduli
of the probability amplitudes of going from $((0,0),N)$ to 
$((r,s),\xi)$ in time $t$ for $\xi \in \{ N,E,S,W \}$.

\subsubsection*{Results}

\begin{figure}[ht]
\centering
\includegraphics[scale=0.68]{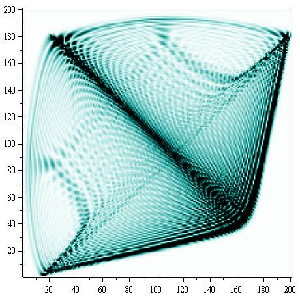}
\caption{Intensity plot of amplitudes at time 200 for a typical quantum walk}
\label{fig:qwalk}
\end{figure}

A number of different two-dimensional QRW's are analyzed
in~\cite{BBBP}.  In these examples, 
the variety defined by the common denominator $\det (I - M)$ 
of the entries of $(I-M)^{-1}$ turns out to be smooth, and
amplitudes may be computed from Theorem~\ref{th:smooth}.
Figure~\ref{fig:qwalk} shows an intensity plot for a 
quantum walk whose unitary matrix $U$ was generated at
random without any symmetries.  
The feasible region, where the probability amplitudes do not decay
exponentially, is a well defined region of irregular shape.  It is 
the image of the torus under the logarithmic Gauss map.  The 
region is determined by the denominator $\det (I-M)$ of the the 
spacetime generating function.  

\begin{figure}[ht]
\centering
\includegraphics[scale=0.38]{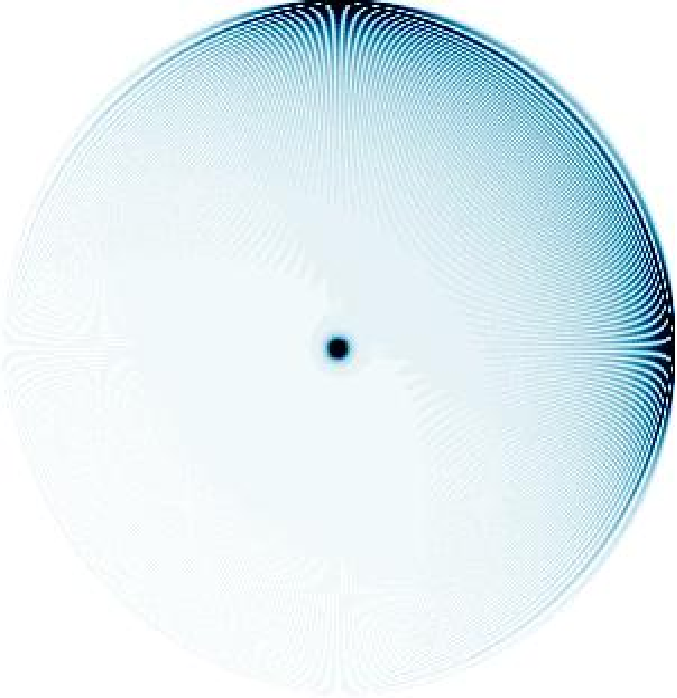}
\caption{two-dimensional Hadamard QRW probability amplitudes 
at time $t=200$}
\label{fig:QRW}
\end{figure}

We now restrict our attention to the Hadamard QRW.  The methods
of~\cite{BBBP} did not suffice to analyze this QRW because 
of the quadratic point and the fact that the denominator $(1-Z^2) Q$ 
is not irreducible.  Two of the factors are binomials $(1 \pm Z)$
and the third factor, $Q$, has quadratic points at $\pm (1,1,1)$,
each of which is on one of the binomial varieties.  The time-200 
intensity plot for this QRW is shown in Figure~\ref{fig:QRW},
where the $x$ and $y$ axes are rescaled spatial variables 
$(r/t, s/t)$, or in other words, velocities.
The feasible region is evidently highly symmetric.  Within the 
feasible region, however, there is asymmetric variation in the
intensities.  This is due to the somewhat arbitrary choice of
starting and ending states, which produce the numerators 
$P_1 , P_2, P_3$ and $P_4$, none of which possesses rotational
symmetry in the $X$-$Y$ plane.  The following result rigorously
establishes the picture in Figure~\ref{fig:QRW}.  

\begin{thm}[feasible region for the two-dimensional Hadamard QRW] 
\label{th:hadamard}
When the velocity is outside the closed disk of radius $\sqrt{1/2}$ 
the amplitudes decay exponentially.  Inside this disk, except possibly 
at the center, the amplitudes do not decay exponentially and
are instead $\Theta (t^{-1})$.
\end{thm}

\begin{unremark}
It will turn out that the 
dominant contribution to the asymptotics of the probability 
amplitudes everywhere except possibly at the center and boundary 
are controlled by smooth points rather than by the two quadratic 
points.  The reason we include this example in the present paper
is that controlling the estimate at the quadratic points will
require the big-O lemma.  
\end{unremark}

\subsubsection*{The amoeba and the normal cone}

Denoting $H_1 = 1-Z, H_{-1} = 1+Z$, we have $F_j = P_j / (Q H_1 H_{-1})$, 
which is in the format of the quadratic point hypotheses with $\eta = 2$.  
The origin is on the boundary of a component $B$ of the complement
of $\amoeba (Q H_1 H_{-1})$ containing the negative $z$-axis.
As before, the cone $\tan_\zero (B)$ is the intersection of
components of $\amoeba (Q)^c$, $\amoeba (H_1)^c$, and 
$\amoeba (H_{-1})^c$.  The latter two are just the halfspaces
$\{ (x,y,z) : z < 0 \}$, which are equal and contain the
component $B_0 := \{ (x,y,z) : z < 0 , \; z^2 > (x^2 + y^2)/ 2$, 
which we recognize from the Aztec Diamond example.  Therefore,
$B = B_0$ and $B^* = \{ t^2 > 2 (r^2 + s^2) \}$ as in 
Section~\ref{ss:aztec-appl}.  An immediate consequence is that
the amplitudes decay exponentially for velocities in directions
outside the disk of radius $\sqrt{1/2}$.  This proves the first
part of Theorem~\ref{th:hadamard}.  

To understand the amplitudes inside the disk, we will first
use the big-O lemma to show that the contribution from the 
quadratic points is $O(t^{-2})$ everywhere except possibly 
at the center (zero velocity).  We will then deduce from
Theorem~\ref{th:smooth} that the contribution from the smooth
points anywhere in the range of the logarithmic Gauss map
is $\Theta (t^{-1})$.  Finally, we will show that the 
range of the logarithmic Gauss map is precisely the disk
of radius $\sqrt{1/2}$.

\subsubsection*{Classification of critical points}

The intersection of $\sing_Q$ with $\sing_{H_j}$ is the set
$$\left \{ Z = j = \frac{X + X^{-1} + Y + Y^{-1}}{4} \right \} \, .$$
The varieties $\sing_{H_1}$ and $\sing_{H_{-1}}$ do not intersect.
We therefore stratify by
\begin{eqnarray*}
\sing_1 & := & \{ (1,1,1)  , (-1,-1,-1) \} \\
\sing_{2+} & := & \sing_Q \cap \sing_{H_1} \setminus \sing_1 \\
\sing_{2-} & := & \sing_Q \cap \sing_{H_{-1}} \setminus \sing_1 \\
\sing_{3+} & := & \sing_{H_1} \setminus (\sing_1 \cup \sing_{2+}) \\
\sing_{3-} & := & \sing_{H_{-1}} \setminus (\sing_1 \cup \sing_{2-}) \\
\sing_4 & := & \sing_Q \setminus (\sing_1 \cup \sing_{2\pm} 
   \cup \sing_{3\pm}) 
\end{eqnarray*}

Again, we are interested in the region of non-exponential decay,
where $\xmax = \zero$; checking where the strata intersect
the unit torus, we find that the strata $\sing_{2\pm}$ do not
intersect the unit torus.  The strata $\sing_{3\pm}$ intersect
the unit torus on the set $\{ |X| = |Y| = Z = 1 \}$, and 
the logarithmic gradient is always in the $t$ direction. 
The factors $H_1$ and $H_{-1}$ cause this direction to be
obstructed and we are therefore not able to say anything 
about asymptotics in the $(0,0,t)$ direction.  This direction
corresponds to the bright spot in the middle of the amplitude 
intensity plot in figure~\ref{fig:QRW}, where there appears
to be a bound state (probability amplitude for being precisely
at the origin does not decay with time).  

Contributions from $\sing_1$ occur for $\rr$ interior 
to $\normal$.  The big-O estimate, Lemma~\ref{lem:big-O} below,
allows us to bound the magnitude of these contributions.
We first compute the homogeneous degree of $F$ at the point
$(1,1,1)$.  The factor $1/Q$ has degree $-2$ here, the factor
$1/H_1$ has degree $-1$, and the factor $1/H_2$ has degree zero.
The numerators $P_j$ vanish to order two at $(1,1,1)$ for 
all $1 \leq j \leq 4$.  We therefore have $\deg (F , (1,1,1))
= \deg (f , (0,0,0)) = 2 - 2 - 1 = -2$.  Applying the lemma,
we find that for the quadratic point $\ww = \zero$, we have
$\contrib (\ww) = O(t^{-2})$, proving the second part of
Theorem~\ref{th:hadamard}.

It is shown in~\cite{brady-thesis}, that the projective range 
of the logarithmic Gauss map is the disk of radius $\sqrt{1/2}$,  
but because this is an unpublished Masters thesis, we will give 
an alternative derivation.  Assuming this for the moment,
it is evident that the contribution from the smooth points 
is $\Theta (t^{-1})$.  If more detail is desired, one can
follow Brady~\cite{brady-thesis} to verify periodicity and
the visually obvious Moire patterns as follows.  Brady shows that
for each $\rr$ interior to $\normal$, there are precisely 
four smooth critical points on the unit torus, a conjugate pair 
and its negative: $\ZZ (\rhat), \overline{\ZZ} (\rhat), 
- \ZZ(\rhat) , -\overline{\ZZ} (\rhat)$.  Denoting $\log \ZZ (\rhat) 
= \zz = i \ww$, Theorem~\ref{th:smooth} tells us that
$$\contrib (\ww) \sim C(\rhat) |\rr|^{-1} \exp ( - i \rr \cdot \ww)$$
where the magnitude of $C(\rhat)$ is proportional to the $-1/2$-power
of the complex curvature of $\log \sing_Q$ at $\zz$.  Adding this
to the contribution from $\overline{\ZZ}$, namely $\contrib (-\ww)$,
we obtain a quantity whose magnitude is $2 \cos \theta (\rr)$
times the magnitude of $\contrib (\ww)$, where $\theta$ is
the argument of $\contrib (\ww)$; note that $\theta$ differs 
from $\rr \cdot \ww$ by $\pi/4$ because the curvature is complex
and its $-1/2$ power has argument $-\sigma \pi / 4$, where
$\sigma$ is the signature of $\qt$, which in our case is $1-2 = -1$;
see~\cite[Section~2.3]{BBBP} for details on 
the phase of the curvature.  Adding the contribution from the
negatives of these two points kills the terms for which $r+s+t$
is odd and doubles the even terms.  This corresponds to periodicity
of the walk.  For fixed $t$, the phase term $\cos \theta (\rr)$ 
varies rapidly (with period of order~1).  Ignoring the Moir\'e pattern
resulting from this term, the probabilities are of order $t^{-2}$
(amplitudes are of order $t^{-1}$) and are spread over the disk 
$r^2 + s^2 = t^2/2$, which is the slice of the normal cone at 
the fixed value of $t$.  We now finish the proof of
Theorem~\ref{th:hadamard} by showing that the range of the
logarithmic Gauss map on $\sing_4$ is the disk of radius
$\sqrt{1/2}$.

\noindent{\sc Proof of the remainder of Theorem}~\ref{th:hadamard}:
Let $(x,y,z)$ be a point on $\sing_4 \cap \tor3$, where
$\tor3$ denote the 3-torus $\{ |x| = |y| = |z| = 1 \}$.
The projective logarithmic Gauss map is the map 
$\gamma_\Proj$ that takes $(x,y,z)$ to the projective
point $(a | b | c)$ where $(a,b,c)$ is the gradient 
of $Q \circ \exp$ at $\log (x,y,z)$.  We represent the
class $(a | b | c)$ by the normalized vector $(a/c , b/c) 
\in \R^2$.  Because the domain is a subset of $\tor3$ we
may use polar coordinates $x=\exp(i\phi), y=\exp(i\psi),
z=\exp(i\theta)$.  In these coordinates, $\sing\cap \tor3$ 
is given by
$$\Sigma=\left\{ 2\cos(\theta)=\cos(\phi)+\cos(\psi)\right\} \, ,$$
where $(\phi, \psi, \theta) \in \tor3$.
The Gauss map takes $(\phi, \psi, \theta)\in\Sigma$ to $-(\sin(\phi),
\sin(\psi), -2\sin(\theta))=:(a,b,c)$. A simple computation shows that
$c^2/2=a^2+b^2+(\cos(\phi)-\cos(\psi))^2/2$, whence the image of the
Gauss map is contained in the two space-like cones $\{c^2/2\geq
a^2+b^2\}$, and its projectivization is contained in the disk
$D:=\{r^2+s^2\leq 1/2\}$.

We first show that the boundary of the disk $D$ belongs to 
the closure of the range of $\gamma_\Proj$.  Near
$(\phi,\psi,\theta)=(0,0,0)$, the surface $\Sigma$ is given by
$$\Sigma=\{\theta^2 - \frac{\phi^2}{2} - \frac{\psi^2}{2}
   + R_3 (\theta , \phi , \psi) = 0\} \, ,$$
where $R_3$ denotes the terms of order three and higher in local
coordinates.  Therefore the projectivization of Gauss map
takes vicinity of the point $(0,0,0)$ to vicinity of the image of
the projectivization of the Gauss map of the conic
$\{\theta^2 - \phi^2/2 - \psi^2/2 = 0\}$, which is the 
boundary of $D$.

To prove, finally, that all of the interior of $D$ is in the range
of the projectivization of the Gauss map, we notice that were that
not the case, this mapping would have some critical values in the
interior of $D$. Critical points of the projectivization of the
Gauss map correspond to the parabolic points of $\Sigma$ (points
where at least one of the principal curvatures of the second
quadratic form of $\Sigma$ vanishes). To find those, we use the
standard trick characterizing the locus of the
parabolic points of the surface
given by $\{f=0\}$ by the equation
$$\langle df, \det H\cdot H^{-1} df\rangle \, ,$$
where $H$ is the Hessian of $f$; $df$ is the gradient (and the matrix
$\det H\cdot H^{-1}$ is the adjunct matrix to $H$). In our case, quick
computation implies that the parabolic points of $\Sigma$ lie on the
intersection of $\Sigma$ with
$$ \left\{\cos(\phi)\cos(\psi)\cos(\theta)\left(
  \frac{2\sin(\theta)^2}{\cos(\theta)}-\frac{\sin(\phi)^2}{\cos(\phi)}
  +\frac{\sin(\psi)^2}{\cos(\psi)}\right)=0\right\} \, ,$$
which, as again a short computation shows, yields critical values only
on the boundary of $D$.  Hence the interior of $D$ is contained in the
range of $\gamma_\Proj$.
$\Cox$

\subsection{Friedrichs-Lewy-Szeg\"o graph polynomials}
\label{ss:4.4}

In the the study of a discretized time-dependent wave equation
in two spatial dimensions, Friedrichs and Lewy required a 
nonnegativity result for the coefficients of $Q^{-1}$ where
$Q(X,Y,Z) := (1-X)(1-Y) + (1-X)(1-Z) + (1-Y)(1-Z)$.
To solve this problem, Szeg\"o~\cite{szego} showed that
the coefficients of $Q^{-\beta}$ are nonnegative for all 
$\beta \geq 1/2$.  Scott and Sokal~\cite{scott-sokal} later
observed that $Q$ is a special case of a spanning tree 
polynomial of a graph.  They proved a generalization of 
Szeg\"o's result to all series-parallel graphs.  Their
results are proved via the stronger property of complete 
monotonicity and are related to the half-plane property.
In order to investigate whether these results might hold
for the polynomials of a larger class of graphs, Scott and
Sokal needed a means of checking the asymptotics  of the 
coefficients: asymptotic nonnegativity is a necessary 
condition for term-by-term nonnegativity.  We will apply
Theorem~\ref{th:no plane} to obtain:
\begin{thm} \label{th:FLS}
Fix $\beta > 1/2$ and let $\sum_\rr a_\rr \ZZ^\rr$ be a 
Taylor expansion for $Q^{-\beta}$.  Then
$$a_\rr \sim \frac{4^{1-\beta}}{\sqrt{\pi} \Gamma (\beta) 
   \Gamma (\beta - 1/2)} (2rs+2rt+2st-r^2-s^2-t^2)^{-1/2}$$
as $\rr$ varies over compact subsets of the cone
$2(rs+rt+st) > r^2+s^2+t^2$.  
\end{thm}

The simplest nontrivial case in which asymptotics may be 
worked out is the one above.  Szeg\"o's 1933 proof of 
nonnegativity was, according to Scott and Sokal, ``surprisingly
indirect, exploiting Sonine-type integrals for products of
Bessel functions.''  It is evident that asymptotics in this case 
may be derived directly from Theorem~\ref{th:no plane}.  We remark 
that the connection between these coefficients and harmonic analysis
of symmetric cones is known to Scott and Sokal, who exploit
the connection and cite several results on the subject
from the sources~\cite{faraut-koranyi,ishi}.

\noindent{\sc Proof of Theorem}~\ref{th:FLS}:
We first check that $\zero$ is on the boundary of the component
$B$ of $\amoeba (Q)^c$ corresponding to the ordinary power series
$1/Q = \sum_\rr a_\rr \ZZ^\rr$.  This follows if we show that 
$Q (X,Y,Z) \neq 0$ for $X, Y, Z$ in the open unit disk.  
To see this, let $D_1$ denote
the open unit disk, let $D_2$ denote the open disk $\{ |z-1| < 1 \}$,
and let $D_3$ denote the halfspace $\{ z : \Real \{ z \} > 1/2 \}$.
The set $D_3$ is the image under $z \mapsto 1/z$ of $D_2$ and
$D_2 = 1 - D_1$.  Therefore, $Q$ has a zero on the open unit polydisk 
$D_1^3$ if and only if $XY+YZ+ZX$ has a zero on $D_2^3$; this
is equivalent to $1/X+1/Y+1/Z$ having a zero on $D_2^3$ which is
equivalent to $X+Y+Z$ having a zero on $D_3^3$.  This is impossible 
because $D_3$ is contained in the open right half-plane.  

Composing with the exponential, then taking the leading homogeneous
part, we obtain
$$\qt = \homog(Q \circ \exp , \zero) = xy+xz+yz \, .$$
We recognize this as half the quadratic factor in 
Section~\ref{ss:cube-grove-appl}.  Therefore $\qt^*$ is twice 
what is was there:
$$ \qt^* (r,s,t) = 2(rs+rt+st) - (r^2+s^2+t^2) \, .$$
Let $P(Z) \equiv 1$.
Recalling that $M$ is chosen so that the matrix for the quadratic 
form is $(M^{-1})^T D M^{-1}$, we see that the determinant of $M$
is $\det (q)^{-1/2}$; plugging in the matrix 
$\disp{\left [ \begin{array}{ccc} 0 & \frac{1}{2} & \frac{1}{2} \\ 
\frac{1}{2} & 0 & \frac{1}{2} \\ \frac{1}{2} & \frac{1}{2} & 0 
\end{array} \right ]}$ for $q$ we obtain $|M| = 2$.  Thus,
for $\beta > 1/2$, equation~\eqref{eq:leading s power} gives
$$a_\rr \sim \frac{4^{1-\beta}}{\sqrt{\pi} \Gamma (\beta) 
   \Gamma (\beta - 1/2)} (2rs+2rt+2st-r^2-s^2-t^2)^{-1/2} \, ,$$
finishing the proof.   $\Cox$

To check this estimate, let $r=s=t=50$ and compute 
$$a_\rr \approx \frac{4^{1-\beta} 7500^{\beta - 3/2}}
   {\sqrt{\pi} \, \Gamma (\beta) \, \Gamma (\beta - 1/2)}
   \approx 0.000222832 \ldots$$
when $\beta = 3/4$.  We then use Maple to crank out the true value 
of $a_{50,50,50}$ which is a rational number near $0.000223464$, for 
a relative error of around $1/400$.

\subsection{Super ballot numbers and multiset permutations} 
\label{ss:super-ballot-appl}

Gessel~\cite{gessel-superballot} defines the 
\Em{super ballot numbers} by
$$g(n,k,r) := \frac{(k+2r)! \, (2n+k-1)!}
   {(k-1)! \, r! \, n! \, (n+k+r)!} \, .$$
These are a generalization of the \Em{ballot numbers} 
$\disp{\frac{k}{2n+k} {2n+k \choose n}}$ (obtained by setting $r=0$),
which are in turn a generalization of the Catalan numbers
(set $k=1$).  The Catalan number and the ballot numbers
are integral and have combinatorial interpretations.  
Gessel shows that the super ballot numbers are integers
as well and sets as a goal to find a combinatorial interpretation.

After re-indexing via $B(a,b,c) := g(a , b-a-c , c)$ for $b > a+c$,
one may extend this definition to all nonnegative $(a,b,c)$ and
obtain a generating function 
$$F(X,Y,Z) = \sum_{a,b,c \geq 0} B(a,b,c) X^a Y^b Z^c = 
   \frac{1-2X}{\sqrt{1-4XZ}} \, G(X,Y,Z)$$
where $G$ is the generating function from equation~\eqref{eq:superballot}.  
The coefficients $N(a,b,c)$ of $G$ are of independent interest.  
They satisfy a similar recurrence to the super ballot numbers
$B(a,b,c)$ but with different boundary conditions.  The numbers
$\{ N(a,b,c) \}$ were shown in~\cite{askey-monograph} to have 
nonnegative coefficients.  They count a difference of cardinalities 
of multi-set permutations~\cite{askey-ismail}.  Gessel goes on to find 
several more identities involving these numbers and their 
generating functions, but no asymptotics are derived.  We will
use Corollary~\ref{cor:simple} to obtain the following asymptotic 
estimate.
\begin{thm} \label{th:superballot}
The numbers $N(a,b,c)$ are estimated asymptotically by
$$N(a,b,c) \sim 2^{a+b+c} \frac{4}{2 \pi} (2ab+2ac+2bc-a^2-b^2-c^2)^{-1/2}$$
uniformly on compact subcones of $\normal = \{ (a,b,c) : 2(ab+ac+bc)
> a^2 + b^2 + c^2 \}$.  
\end{thm}

Our motivation for analyzing the numbers $\{ N(a,b,c) \}$ is
admittedly ``because we can''.
The coefficients of $F$ are of greater interest than the coefficients
of $G$, but the fractional power on the non-quadratic term takes
this problem beyond the main results of this paper.  The deformations
in Section~\ref{sec:homotopies} still apply, but the further analysis
in Section~\ref{ss:cone plane} via Leray cycles does not work when this
factor is algebraic rather than a simple pole.  We therefore do not
state detailed asymptotics for $F$, reserving this for future work.

Let $Q(X,Y,Z) := 2 - X - Y - Z + XYZ$, so that $2/Q$ is the ordinary
power series generating function for $2^{-a-b-c} N(a,b,c)$.  It is
easy to check (e.g., via Gr\"obner bases) that $Q$ and its gradient
vanish simultaneously exactly at the point $\one$.  We have $q := 
Q \circ \exp = xy + yz + xz + O(|(x,y,z)|^3$, whose homogeneous part 
(at $\zero$) is given by $\qt = xy + yz + xz$.  Again, 
$$\qt^* (r,s,t) = 2(rs+rt+st) - (r^2+s^2+t^2)$$
and $|M| = 2$.  

There is a component $B$ of $\amoeba (Q)^c$ containing 
a translate of the negative orthant (corresponding to the 
ordinary power series expansion); let us check that 
$\zero \in \partial B$.  Proceeding as in Section~\ref{ss:4.4},
it suffices to verify that $Q$ has no zero in the open unit 
polydisk $D_1$, which is equivalent to checking that
$Q(\one + \ZZ)$ has no zero in $-D_2^3$ where $D_2 = \{
z : |z+1| < 1 \}$.  We have
$$Q(\one+\ZZ) = XYZ + XY + XZ + YZ 
   = XYZ \left ( 1 + \frac{1}{X} + \frac{1}{Y} + \frac{1}{Z} \right ) \, ,$$
whence this is further equivalent to $1 + X + Y + Z$ having
no zero in $-D_3^3$, where $-D_3$ is the halfplane $\{ z : \Real \{ z \}
< -1/2 \}$.  This is obvious, because the real part of $X + Y + Z$
is bounded above by $-3/2$ on $D_3^3$.  

We now apply Corollary~\ref{cor:simple} to obtain the asymptotics
of $2^{-a-b-c} N(a,b,c)$ inside the cone $\normal$, these 
asymptotics being exponentially small outside $\normal$.
Letting $P(Z) \equiv 2$, we plug $P$, $|M|$ and $\qt^*$ 
into~\eqref{eq:leading simplest} to obtain the leading term 
asymptotics for the coefficients of $2/Q$, which gives the 
expression in Theorem~\ref{th:superballot} and finishes the proof.
As an example, if $a=1, b=20,c=30$, then the approximation yields 
$N(a,b,c) \approx 2.595 \times 10^{16}$ while the actual value of
$N(10,20,30)$ to three decimal places is $2.547 \times 10^{16}$.
$\Cox$

\setcounter{equation}{0}
\section{Homotopy constructions} \label{sec:homotopies}

Recall from the heuristic discussion following~\eqref{eq:new cauchy}
that moving the chain of integration in~\eqref{eq:mv cauchy} to the 
torus $\Log^{-1} (\xmax)$ is not enough.  Our goal in this 
section is to construct homotopies moving this chain of integration, 
within the domain of holomorphy of the integrand (but not within
the domain of convergence of the Laurent series) to a different
chain on which the maximum modulus of the integrand is small
except in a neighborhood of $\critical (\rr)$.  While the Morse 
theoretic methods of~\cite{GM} are in principle constructive, 
we follow~\cite{ABG}, taking advantage of hyperbolicity in order 
to produce vector fields along which chains may be shifted.

\subsection{Vector fields} \label{ss:VF}

Let $\flattorus := (\R / 2 \pi \Z)^d$ denote the $d$-dimensional
flat torus.  Given a Laurent polynomial $F$ and a component $B$
of the complement of $\amoeba (F)$, pick a unit vector $\rhat$ in the 
interior of the convex dual $-B^*$.  We suppose that $\rr$ is 
a proper direction for $B$.  Let
\label{U}
\begin{equation} \label{eq:U}
U = \bigcup_{\ww \in \logcrit (\rr)} U_\ww
\end{equation}
be the disjoint union of neighborhoods of each $\ww \in \logcrit (\rr)$,
where $\logcrit (\rr)$ are the logarithmic critical sets from 
Definition~\ref{def:crit}.
\begin{lem}[Vector field away from the critical set]
\label{lem:eta hyp}
Let $F$ be a Laurent polynomial with $f := F \circ \exp$, let $B$ 
be a component of $\R^d \setminus \amoeba (f)$, and let $\rhat$ 
be a unit vector in the interior of the convex dual $-B^*$.  
Suppose $- \rhat \cdot \xx$ is minimized at a unique $\xmax$ 
in $\partial B$ and define the
neighborhood $U$ of $\logcrit (\rr)$ as in~\eqref{eq:U}.
Then there is a smooth vector-valued function $\eta_{U^c} : 
\label{eta}
\flattorus \setminus U \to \R^d$ such that:
\begin{enumerate} \romenumi
\item $\eta_{U^c} (\yy) \in \cK^{f,B} (\exp (\xmax + i \yy))$; 
\item $\rhat \cdot \eta_{U^c} (\yy) = 1$ for all $\yy \in \flattorus
   \setminus U$.
\end{enumerate}
\end{lem}

\noindent{\sc Proof:}  First, for each $\yy \notin U$, we will 
find a neighborhood $\nbd_\yy$ and a vector $\vv_\yy$ such that 
$\eta_{U^c} \equiv \vv_\yy$ fulfills $(i)$--$(ii)$ on $\nbd_\yy$.

Fix $\yy \notin U$.  If $f (\xmax + i \yy) \neq 0$
then choose a neighborhood $\nbd_\yy$ of $\yy$ in $\R^d$ 
such that for $\vv \in \nbd_\yy$, the quantity $f (\xmax + i \vv)$ 
does not vanish.  Choose $\vv_\yy$ with $\rhat \cdot \vv_\yy = 1$.

Alternatively, suppose that $f (\xmax + i \yy) = 0$.  By 
Proposition~\ref{pr:hyperbolic}, the homogeneous part, 
call it $A_\yy$, of the function $\vv \mapsto f(\xmax + i \yy + \vv)$ 
is real and hyperbolic, and by Proposition~\ref{pr:relinearization},
there is a cone $K$ of hyperbolicity containing $\tan_{\xmax} (B)$.  
Also by the first part of Proposition~\ref{pr:dir equiv}, there is some 
$\vv_\yy \in K$ with $\rhat \cdot \vv_\yy = 1$.  
By semi-continuity (part~$(i)$ of Corollary~\ref{cor:semi}), 
$\vv_\yy \in \cK (\exp(\xmax + i \uu))$ for every $\uu$ in some 
neighborhood $\nbd_\yy$ of $\yy$.  

The collection $\{ \nbd_\ww : \ww \notin U \}$ covers $U^c$;
shrink it slightly if necessary so that the closure of its union 
does not intersect $\logcrit (\rr)$.  We may (shrinking some
of the $\eta_\ww$ slightly if necessary) choose a finite subcover 
$\{ \nbd_\ww : \ww \in \Xi \}$ whose union has closure disjoint 
from $\logcrit (\rr)$.  Choose a partition of unity $\{ \psi_\ww : 
\ww \in \Xi \}$ subordinate to the subcover.  Define
$$\eta_{U^c} (\yy) := \sum_{\ww \in \Xi} \psi_\ww (\yy)
   \vv_\ww (\yy) \, .$$
Now~$(i)$ is satisfied by convexity and~$(ii)$ is satisfied
by linearity.   $\Cox$

\begin{cor}[Vector field defined everywhere] \label{cor:eta global}
Under the conditions of Lemma~\ref{lem:eta hyp} there is a smooth
vector field $\eta$ satisfying 
\label{eta2}
\begin{eqnarray}
\eta (\yy) & \in & \cK^{f , B} (\exp(\xx + i \yy)) \, ;
   \label{eq:eta 1} \\
\rhat \cdot \eta (\yy) & = & 1 \;\; \mbox{ on } U^c \, .
   \label{eq:eta 2}
\end{eqnarray}
\end{cor}

\noindent{\sc Proof:} Let $\eta_\ww : U_\ww \to \C^d$ be any map
for which $\eta_\ww (\yy) \in \cK^{\ft , B} (\exp(\xx + i \yy))$.  
To see that we may choose such a map smoothly, note that the
constant map $\eta_\ww (\yy) \equiv \vv$ is such a map
whenever $\vv \in \tan_\xx (B)$.  The reason for allowing
a general function $\eta_\ww$ in place of a constant is 
that later we will use~\eqref{eq:eta} with functions $\eta_\ww$
tailored to more specific needs.  The collection $\{ \nbd_ \ww : 
\ww \in \Xi \} \cup \{ U_\ww : \ww \in \logcrit (\rhat) \}$ 
covers $\flattorus$.  Choose a partition 
of unity $\{ \psi_\ww \}$ subordinate to this and define
\begin{equation} \label{eq:eta}
\eta (\yy) := \sum_{\ww \in \Xi} \psi_\ww (\yy) \vv_\yy (\ww)
   + \sum_{\ww \in \logcrit (\rhat)} \eta_\ww (\yy) \, .
\end{equation}
This proves the corollary.  
$\Cox$

We remark for later use that if $\rhat$ is replaced by a non-unit
vector $\rr$, then applying the above constructions to $\rhat$
replaces~\eqref{eq:eta 2} by $\rr \cdot \eta (\yy) = |\rr|$ on $U^c$.
Next, we give a projective version of the above construction.
We say that a 1-homogeneous function $\phi$ is smooth if 
it is smooth away from the origin.  

\begin{lem}[projective vector field] \label{lem:eta proj}
Let $A$ be a real homogeneous polynomial in $d$ variables
of degree $m \geq 1$ and let $B$ be a cone of hyperbolicity
for $A$ whose dual $-B^*$ has nonempty interior.  For each 
$\yy \in \R^d$, recall the cone $\cK^{A,B} (\yy)$ defined in 
Proposition~\ref{pr:relinearization}.  Let $\rr$ be a non-obstructed
vector in the interior of $-B^*$.  Then there is a 1-homogeneous, smooth 
vector field $\eta$ on $\R^d$ such that for all $\yy \in \R^d$
and all $\rr'$ in a neighborhood of $\rr$,
\begin{enumerate} \romenumi
\item $\eta (\yy) \in \cK^{A , B} (\yy)$;
\item $\rr' \cdot \eta (\yy) \geq |\rr'| |\yy|$. 
\end{enumerate}
\end{lem}

\noindent{\sc Proof:} This is a homogeneous version of the
proof of Lemma~\ref{lem:eta hyp}.  Assume first that $|\yy| = 1$.
We define $\eta$ locally and then piece these
together via a partition of unity.  When $A(\yy) \neq 0$ we
can find neighborhoods $\nbd_\yy$ of $\yy$ and $\nbd_\yy^*$ of 
$\rr$ such that $A$ vanishes nowhere on $\nbd_\yy$ and 
there is a $\vv$ for which $\rr' \cdot \vv > |\rr'|$ on $\nbd_\yy^*$.
By the trivial part of the definition, $\vv_\yy \in \cK^{A,B} (\yy)$.

When $A(\yy) = 0$, because $\rr$ is non-obstructed, there is 
a vector $\vv_\yy \in \cK^{A,B} (\yy)$ with $\rr \cdot \vv_\yy > 0$.
By semi-continuity (part~$(ii)$ of Corollary~\ref{cor:semi}),
$\vv_\yy \in \cK^{A,B} (\uu)$ for every $\uu$ in some neighborhood
$\nbd_\yy$ of $\yy$.  By continuity, $\rr' \cdot \vv_\yy > 0$ for 
every $\rr'$ in some neighborhood $\nbd_\yy^*$ of $\rr$.  We may 
then replace $\vv_\yy$ by some positive multiple so that 
$\rr' \cdot \vv_\yy > |\rr'|$ for $\rr \in \nbd_\yy^*$.

To define the 1-homogeneous function $\eta$, it suffices to
define it on the set $S_1$ of vectors $\yy$ of norm~1.  Cover
$S_1$ by finitely many neighborhoods $\{ \nbd_\ww : \ww \in \Xi \}$
and use a partition of unity subordinate to the cover to define
$\eta$ via~\eqref{eq:eta} on $S_1$.  Extending 1-homogeneously 
via $\eta (\lambda \yy) := \lambda \eta (\yy)$ finishes the 
construction.    $\Cox$

\subsection{Homotopies} \label{ss:hom}

Any piecewise differentiable map from a compact manifold to 
another manifold defines a chain of integration.  
Let $\eta$ be any continuous vector field on $\flattorus$ and
fix any $\ee > 0$.  Define the homotopy $\Phi = \Phi^{\ee , \eta} : 
\label{Phi}
\flattorus \times [0,1] \to \C^d$ by
\begin{equation} \label{eq:Phi}
\Phi_t (\yy) := i \yy + \xx 
   + \ee \left [ (1-t) \uu + t \eta (\yy) \right ] 
\end{equation}
where $\uu$ is fixed vector in $\tan_\xx (B)$.
We specialize now to $\eta$ given by Corollary~\ref{cor:eta global}.
Setting $t=0$ gives a cycle (thinking of the map as a chain 
of integration) whose range is the torus
$T := \xx + \ee \uu + i \flattorus$.  Setting $t=1$ gives another cycle, 
which we call $\generic (\eta)$, shown to be homotopic to $T$ in $\C^d$ 
via the homotopy $\{ \Phi_t : 0 \leq t \leq 1 \}$.  
\label{Ceta}

\begin{thm}[The homotopy defined by $\eta$ avoids $\sing_{f}$]
\label{th:avoid f}
Let $\eta$ satisfy~\eqref{eq:eta 1}--\eqref{eq:eta 2} and define
$\Phi_t$ by~\eqref{eq:Phi}.  Then for $\ee > 0$ sufficiently small 
and all $0 \leq t \leq 1$, $f (\Phi_t (\yy)) \neq 0$.
\end{thm}

\noindent{\sc Proof:}  For fixed $t$ this follows from 
Corollary~\ref{cor:section}.  To find $\ee$ that works
for all $t$ simultaneously, use compactness of the
interval in $t$ and continuity of the homotopy.
$\Cox$

Let $\{ U_\ww : \ww \in \logcrit (\rr) \}$ be disjoint open 
neighborhoods of the points of $\logcrit (\rr)$ as before, and 
\label{Cw}
let $\generic (\ww)$ denote the restriction of $\generic (\eta)$ 
to the closure of $U_\ww$.  Let $\generic (U^c)$ denote the 
restriction of $\generic (\eta)$ to the closure of $U^c$.  
The chain $\generic (\eta)$ is representable as a sum
$$\generic (\eta) = \generic (U^c) 
   + \sum_{\ww \in \logcrit (\rr)} \generic ( \ww) \, .$$
An immediate consequence of the previous constructions is:
\begin{cor}[localization of the Cauchy integral] \label{cor:local chains}
The chain $T$ is homotopic in $(\C / (2 \pi \Z))^d \setminus 
\sing_{f}$ to a sum of chains
$$\generic (U^c) + \sum_{\ww \in \logcrit (\rr)} \generic (\ww)$$
where $\generic (U^c)$ and each $\partial \generic (\ww)$ are 
supported on the set of $\zz$ such that $\rr \cdot \Real \{ \zz \} 
- \rr \cdot \xx \geq \ee |\rr|$.  Consequently, the 
integral~\eqref{eq:cauchy torus} decomposes as
\begin{eqnarray*}
a_\rr & = & \left ( \frac{1}{2 \pi i} \right )^d 
   \int_{\xx + \ee \uu + i \flattorus}
   e^{- \rr \cdot \zz} \frac{1}{f (\zz)} \, d\zz \\
& = & R + \sum_{\ww \in \logcrit (\rr)}\left ( \frac{1}{2 \pi i} \right )^d
   \int_{\generic (\ww)} e^{- \rr \cdot \zz} \frac{1}{f (\zz)} \, d\zz 
\end{eqnarray*}
where $R = O(e^{-\ee |\rr|})$, and Theorem~\ref{th:localize} follows.
$\Cox$
\end{cor}

While the above general construction, e.g., with $\eta_\ww \equiv \uu$,
suffices to localize the Cauchy integral, explicit computations
for a cone and a plane, done in Section~\ref{ss:cone plane}, 
will require specific choices of $\eta_\ww$, resulting in 
specific chains $\dchain (\ww)$ which we will use in the 
decomposition from Corollary~\ref{cor:local chains}.  Say 
that $\eta_\ww$ is \Em{projective} if $\eta_\ww (\ww + \cdot)$ 
is homogeneous of degree~1 and smooth away from $\zero$.  
The first of these two results follows immediately from 
Lemma~\ref{lem:eta proj}.

\begin{thm} \label{th:projective}
Let $A$ be a hyperbolic homogeneous polynomial with cone 
of hyperbolicity $B$.  Let $\rr$ in the interior of $-B^*$ be
non-obstructed and let $\eta$ be the projective vector field
of Lemma~\ref{lem:eta proj}.  Let $\{ \Phi_t \}$ be the homotopy
on $\R^d$ defined by~\eqref{eq:Phi} with $\xx = \zero$.  
Then $A(\Phi_t (\yy)) \neq 0$ 
for all $0 \leq t \leq 1$ except when $t=1$ and $\yy = 0$, and
$\rr \cdot \Phi_1 (\yy) \geq c |\yy|$.  Consequently, for
$\uu \in \tan_\xx (B)$, the chain $\uu + i \R^d$ is homotopic
through the complement of $\sing_A$ to the projective
\label{proj}
chain $\proj := \Phi_1 [\R^d]$ on which $\rr \cdot \yy$
grows linearly in $|\yy|$.  This construction is uniform as
$\rr$ varies over some neighborhood $\nbd$.
$\Cox$
\end{thm}

The projective chain $\proj$ provides the concept we need,
but this idealized chain has two problems: it is infinite,
and it touches $\sing_A$ at the origin.  To take care of
the second problem, we stop the homotopy early in a 
small ball about the origin.
\begin{defn} \label{def:conechain}
\label{projd}
Let $\proj^{(\delta)}$ denote the chain parametrized by
$\flattorus$ defined by 
$$\yy \mapsto i \yy + \ee [ (1 - t(1-(\delta - |\yy|)^+))) \uu
   + t(1 - (\delta - |\yy|)^+) \eta (\yy) ] \, .$$
This is obtained by replacing $t$ in the definition of
the homotopy~\eqref{eq:Phi} by $t (1 - (\delta - |\yy|)^+)$
where $\eta = \eta_\ww$ is the projective vector field 
near $\ww$.  The definition does not change the homotopy 
outside the ball of radius $\delta$, but inside this ball 
the homotopy stops early, stopping at time $1-\delta$ at 
the origin and interpolating linearly in $|\yy|$.  
\end{defn}

Finally, we glue together pieces looking like $\proj^{(\delta)}$
near each point $\ww \in \logcrit$ to produce the chains
we will use to prove all the remaining results.

Let $F$ be a Laurent polynomial and let $B$ be a component of 
$\R^d \setminus \amoeba (F)$.  Suppose $\rr$ is proper with
dual point $\xmax$, that $\logcrit (\rr)$ is finite, and that
$\rr$ is non-obstructed.  For each $\ww \in \logcrit (\rr)$,
let $\ft_\ww = \homog (f , \ww)$ and let $\eta_\ww$ be the
vector field constructed in Lemma~\ref{lem:eta proj} with
$\ft_\ww$ in place of $A$.  

We piece these together into one locally projective vector
field on the torus via a partition of unity as before.
Let $\{ U_\ww : \ww \in \logcrit (\rr) \}$ and $U$ be 
defined as in~\eqref{eq:U} and define $\eta$ by~\eqref{eq:eta}.
Define $\{ \Phi_t \} = \{ \Phi_t^{\ee , \delta} \}$ by~\eqref{eq:Phi}
with $t (1 - [\delta - d(\yy))^+]$ replacing $t$, where
$d(\yy) := \min_{\ww \in \logcrit} |\yy - \ww|$ is the
minimum distance from $\yy$ to a point of $\logcrit$.
\label{dchainw}
We let $\generic$ denote the chain $\Phi_1$
and for each $\ww$ we let $\dchain (\ww)$ denote
the intersection of $\generic$ with the 
radius-$\delta$ neighborhood of $\ww$.

\begin{thm}[locally projective homotopy] \label{th:eta cone}
If $\ee , \delta > 0$ and the neighborhoods $\{ U_\ww \}$ 
are taken to be sufficiently small, then the homotopy 
$\{ \Phi_t^{\ee , \delta} \}$ will avoid $\sing_{f}$.  
In particular, the chain $\xx + \ee \uu + i \flattorus$
is homotopic in the complement of $\sing_f$ to the chain 
\label{dchain}
$$\dchain = \generic_{U^c} 
   + \sum_{\ww \in \logcrit} \dchain (\ww)$$
for which the inequality
\begin{equation} \label{eq:min}
\rr' \cdot \Phi_1 (\yy) - \rr' \cdot \xx \geq 
   c \min_{\ww \in \logcrit (\rr)} |\yy - \ww|
\end{equation}
will be satisfied for some $c > 0$, for every $\yy \in \flattorus$
and every $\rr'$ in some neighborhood of $\rr$.
\end{thm}

\noindent{\sc Proof:}  We have constructed deformations 
using vector fields $\eta_\ww$ defined in terms of local 
homogenizations $\ft_\ww$, so the main content of the proof
is to ensure that the homotopy avoids the actual zero set
$\sing_f$ and not only the homogeneous approximation to it.
In any set bounded away from the critical points, this is
automatic.  It suffices to consider what happens in a neighborhood 
of a critical point $\zz = \exp (\xx + i \ww)$.  Here, what
we want is in fact true in considerable generality: moving
the origin to $\xx + i \ww$, any projective set avoiding 
$\sing_{\ft_\ww}$ avoids $\ft$ in a neighborhood of the origin.
The range of the homotopy $\Phi_t^{\ee , \delta}$ is a 
projective set in a neighborhood of the origin, meaning 
that it is locally a closed conical set of the form
$$\{ \lambda \vv : \vv \in K , \lambda \in [0,\ee] \} \, ,$$
where $K$ is a closed subset of the unit sphere.  
The set $\sing_{\ft_\ww}$ is a also a closed conical set.  On
the unit sphere, these two closed sets do not intersect and hence
are separated sphere by a positive distance.  
When $\delta$ is small enough, the normalized points 
$\uu / |\uu|$ for $\uu \in \sing_f$ are within
$\ee / 2$ of the points of $\sing_{\ft_\ww}$ on the unit 
sphere when $|\uu| < \delta$.  Thus the homotopy $\Phi_t^{\ee , \delta}$
avoids $\sing_f$, and the theorem follows.  
$\Cox$

\subsection{Consequences and an example} \label{ss:cone and plane}

As outlined in Section~\ref{ss:overview}, the deformations 
constructed in Theorems~\ref{th:avoid f} and~\ref{th:eta cone}
allows us to localize and then compute the Cauchy integral.
These computations are carried out in the next section using
Fourier apparatus.  We record here a preliminary estimate that is 
useful in a more general context (see, e.g.,~\cite{brady-thesis}).  If
\begin{equation} \label{eq:F}
F := \prod_{j=1}^k Q_j^{s_j}
\end{equation}
is the product of $d$-variate polynomials to arbitrary real powers
and $\ZZ$ is any complex vector, the homogeneous degree
$\deg (F , \ZZ)$  of $F$ at $\ZZ$ is defined by
$$\deg (F , \ZZ) := \sum_{j=1}^k s_j \deg (Q_j , \ZZ)$$
(it is easy to check that this is independent of the
representation of $F$ as such a product).  

\begin{lem}[big-O estimate] \label{lem:big-O}
Let $F = \prod_{j=1}^k Q_j^{s_j}$ and let $H$ denote the 
product of all the $Q_j$ for which $s_j$ is not a positive 
integer, so that $\sing_H$ is the singular locus of $F$.  
Let $f := F \circ \exp$ and let $\{ a_\rr \}$ be the coefficients
of a Laurent series for $F$ corresponding to the component
$B$ of $\R^d \setminus \amoeba (H)$.  Fix a proper, 
non-obstructed direction $\rr$ in the interior of $-B^*$
and let $\xmax \in \partial B$ be the minimizing point for $\rr$.
For $\ww \in \logcrit (\rr)$, let $\wchain$ denote any of
the three chains $\proj (\ww)$, $\proj^{(\delta)} (\ww)$ or
$\dchain (\ww)$.  Then 
\begin{enumerate} \romenumi
\item If $\phi (\zz) = O(|\zz|^\beta)$ and $\beta + d > 0$ then
\begin{equation} \label{eq:big-O}
|\zz^\rr| \int_{\wchain} \exp (- \rr \cdot \zz) \phi (\zz) \, d\zz
   = O \left ( |\rr|^{-d-\beta} \right ) \, .
\end{equation}
\item Consequently, if $\deg (f , \ww) + d > 0$, then 
for any bounded function $\psi$, the following integral 
is absolutely convergent and
$$|\zz^\rr| \int_{\wchain} e^{-\rr \cdot \zz} \psi (\zz) f(\zz) \, d\zz
   = O \left ( |\rr|^{- (d + \deg (f , \ww))} \right )  \, .$$
\item It follows further that the Taylor coefficients $a_\rr$ of $F$ 
satisfy $a_\rr = O \left (|\zz|^{-\rr} |\rr|^{- \alpha} \right )$ where 
$$\alpha = d + \max_{\ww \in \logcrit (\rr)} \deg (f , \ww) \, .$$
\end{enumerate}
\end{lem}

\noindent{\sc Proof:} We prove~$(i)$ and~$(ii)$ first for $\wchain = 
\proj$.  The chain $\proj$ is a cone avoiding the singular locus of 
the homogeneous function $f$ except at zero.  Everything is homogeneous,
so it is just a matter of keeping track of degrees.  

Let $S$ denote the section $\{ \zz : \Real \{ \rhat \cdot \zz \} = 1 \}$ 
of this cone.  We have \label{chain} $\CC_\ww = [0,\infty) \times S$.  
Write $f (\zz) = |\zz|^{\deg (f , \ww)} F_0 (\zz / |\zz|)$ 
for some smooth function $F_0$ on $S$ and decompose $d\zz = 
t^{d-1} \, dt \wedge dS$ for some form $dS$ on $S$.  Let
$M := \int_S |F_0 (\uu)| \, dS(\uu)$ and $M' := \sup |\psi|$.
Integrating first over $S$ then over $[0,\infty)$ gives
\begin{eqnarray*}
\left | \zz^\rr \int_{\CC_\ww} e^{-\rr \cdot \zz} \psi (\zz) f (\zz) 
   \, d\zz \right | 
   & = & \left | \int_0^\infty t^{\deg (f , \ww)} t^{d-1} dt
   e^{- |\rr| t} \left [ \int_S \psi (\uu) F_0(\uu) \, dS \right ] \right | \\
& \leq & \int_0^\infty e^{- |\rr| t} t^{d + \deg (f , \ww) - 1} \, M \, M' 
   \, dt 
\end{eqnarray*}
which is absolutely convergent and $O(|\rr|^{- (d + \deg (f , \ww))})$
as desired.  

For sufficiently small $\delta > 0$, the chains $\proj^{(\delta)}$ 
are homotopic in the domain of holomorphy of the integrand, whence
the value of the integral is independent of the particular value
of $\delta$.  As $\delta \downarrow 0$ the integrals over over the
parts where $\proj^{(\delta)} \neq \proj$ converge to zero (by the
same estimates) while the integrals over the parts where 
$\proj^{(\delta)} = \proj$ converge to the integral over $\proj$
by the definition of the Lebesgue integral.  This proves the result
for $\wchain = \proj^{(\delta)} (\ww)$.  For $\wchain = \dchain (\ww)$,
observe that difference is the integral over a set where 
$$\Real \{ - \rr \cdot \zz \} < - \rr \cdot \xx - c$$
for some $c = c (\delta) > 0$.  This exponentially small term
is smaller than the remainder term in the conclusion, so the 
theorem holds for these chains as well.

The third conclusion now follows from localization 
(Theorem~\ref{th:localize}).
$\Cox$

We close the section with an example of the conclusion
of Lemma~\ref{lem:eta proj} and Theorem~\ref{th:eta cone} 
in the case of the product of a quadratic cone $Q$ and a 
linear function $H$.  
We give an explicit construction a vector field $\eta$ satisfying 
the conclusion of the lemma; this explicit construction will
be useful in computing an integral in Section~\ref{ss:cone plane}.

We consider the case of $F = 1/(QH)$, where $q := Q \circ \exp$
and $h := H \circ \exp$ are respectively quadratic and smooth
at the origin:
\begin{eqnarray*}
q & = & \qt + R_1 \\
h & = & \hht + R_2
\end{eqnarray*}
with $\qt$ homogeneous quadratic, $\hht$ linear, $R_1 (\yy) = O(|\yy|^2)$, 
and $R_2 (\yy) = O(|\yy|^3)$.  The signature of $\qt$ is assumed 
to be $(1,2)$ whence the zero set of $\qt$ in real space is a cone 
over a circle.  Suppose that 
the zero sets of $\qt$ and $\hht$ intersect transversally in real 
space.  In other words, the plane $\{ \hht = 0 \}$ intersects the
cone $\{ \qt = 0 \}$ in two lines; 
in projective space, the line $\{ \hht = 0 \}$ intersects the
circle $\{ \qt = 0 \}$ in two points.

The construction of $\eta$ depends on the choice of the cone of
hyperbolicity $B$ of $\qt \hht$; in the applications below, this is
the component of $\amoeba (QH)^c$ containing the negative half
of the $z$-axis. Let us assume in this example that $A := \qt \hht$ is
hyperbolic with respect to $-e_3$.  We also assume without loss of
generality that $\hht (-e_3) > 0$.  We have seen in several of the
examples that $B := B_1 \cap B_2$ where $B_1$ is a halfspace dual to
$\hht$ and $B_2$ is the projective ellipse defined by $\qt$.  
The dual cone is the cone over the teardrop pictured in 
figure~\ref{fig:teardrop}: here the conic is the dual to $B_1$, 
and the vertex is the line dual to the hyperplane $B_2$; 
as usual, the dual to the intersection of the convex
sets is the convex hull of their respective duals.  

We will construct the section $\eta$ guaranteed by 
Lemma~\ref{lem:eta proj} for which $\rr \cdot \eta (\yy) > 0$, 
along with a null section $\tilde{\eta}$ satisfying 
$\rr \cdot \tilde{\eta} = 0$ that is needed in
Section~\ref{ss:cone and plane}.

Fix $\rr \in \normal$.  There are two nearly identical cases, 
depending on whether or not $\rr \in B_2^*$.  Assume first that 
$\rr \in B_2^*$.  In figure~\ref{fig:eta} below, height is the 
linear functional defined by $-\rr$, so the plane 
$X_\rr:=\{ \rr \cdot \xx = 0 \}$ is drawn as as horizontal, 
with $\rr \cdot \xx$ increasing as one goes downward.  
This plane intersects the real cone $\{ \qt = 0 \}$ only 
at the origin, hence the cone has a positive (lower) and 
a negative (upper) half; we have assumed $-e_3$ is in 
the upper half.  The construction of $\tilde{\eta}$
is automatic if we mandate that 
$\tilde{\eta} (\yy) = -e_3 |\yy| + \lambda \yy$
for some $\lambda$.  In order to obtain $\rr \cdot 
\tilde{\eta} (\yy) = 0$, we need to take 
\begin{equation} \label{eq:lambda}
\lambda = |\yy| \frac{\rr \cdot e_3}{\rr \cdot \yy} \, .
\end{equation}

Wherever $\rr \cdot \yy \neq 0$, this is clearly smooth and
1-homogeneous.  Setting aside the points where $\rr \cdot \yy = 0$,
at every other point
$\yy$ where $\qt$ or $\hht$ vanishes but not both,
the cone $\cK^{A,B} (\yy)$ is a halfspace bounded by 
the tangent plane at $\yy$ to $\{ A = 0 \}$.  This 
halfspace contains the vector $\yy$, making it obvious that
of the two halfspaces bounded by this plane, $\tilde{\eta} (\yy)$ 
is in the one containing $-e_3$, thus is in $\cK^{A,B} (\yy)$. 
When $\yy$ is in the intersection $\qt = \hht = 0$, again
$\tilde{\eta} (\yy) \in \cK^{A,B} (\yy)$, because this cone
is the intersection of two halfspaces, each of which we have
seen to contain $\tilde{\eta} (\yy)$.  Finally, to deal with
the points where $\rr \cdot \yy = 0$, note that $\qt$ is nonvanishing
(outside of the origin) on this set.  In a neighborhood of the point $\hht = 0$, we use 
a smooth bump function $\psi_\ee : \R \to [0,1]$ that is one
on $[-\ee,\ee]$ and zero outside $[-2\ee , 2\ee]$.  Letting
$\xx$ be any vector with $\xx \cdot \rr = 0$ and $\hht(\xx) > 0$,
we take
$$\tilde{\eta} (\yy) 
   := |\yy| \psi_\ee (\rr \cdot \hat{\yy}) \xx 
   + (1 - \psi_\ee) \left ( - |\yy| e_3 + \lambda \yy \right )$$
where $\lambda$ is defined by~\eqref{eq:lambda}.  This completes
the construction of $\tilde{\eta}$.  When $A(\yy) \neq 0$, the
cone $\cK^{A,B} (\yy)$ is all of $\R^d$, so verification that
$\eta (\yy) \in \cK^{A,B} (\yy)$ is trivial; we conclude that 
$\tilde{\eta}$ is a 1-homogeneous section of $\cK^{A,B} (\cdot)$ 
with $\rr \cdot \eta \equiv 0$.  

What we have accomplished is to find a single formula
for $\eta$ that works for all strata of $\{ A = 0 \}$,
resorting to partitions of unity, only in one place, away
from $\{ \qt = 0 \}$; this will be useful in the sequel.  
For a pictorial description of the construction we have just
completed, see figure~\ref{fig:eta}.  In the upper half of
the cone, $\tilde{\eta}$ points inside, as does $-e_3$.  In the lower
half (not shown), both $\tilde{\eta}$ and $-e_3$ point outside.
To obtain the vector field of Lemma~\ref{lem:eta proj}, we
set $\eta (\yy) = \tilde{\eta} (\yy) + \ee |\yy| e_3$ for a
sufficiently small $\ee$.

\begin{figure}[ht]
\begin{center}
\includegraphics[scale=0.48]{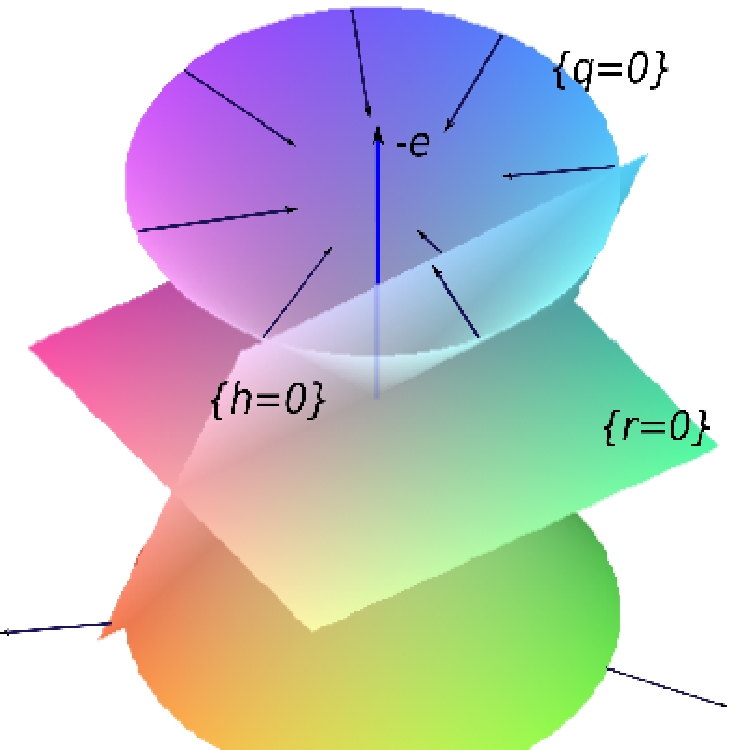}
\caption{the vector field $\tilde{\eta}$}
\label{fig:eta}
\end{center}
\end{figure}

\begin{figure}[ht]
\hspace{0.8in} \includegraphics[scale=0.48]{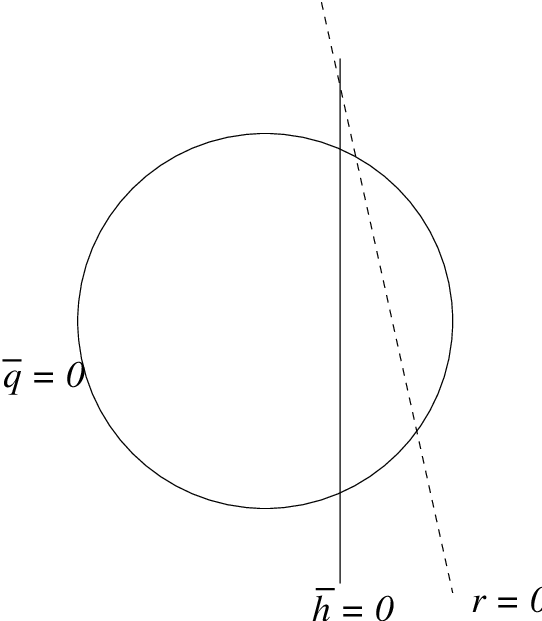}
\caption{when $\rr$ is outside $B_2^*$}
\label{fig:eta2}
\end{figure}

When $\rr$ is outside the circle, but inside the teardrop, 
the plane orthogonal to $\rr$ does intersect the real cone 
$\qt = 0$.  In projective space, the analogue of figure~\ref{fig:eta}
is figure~\ref{fig:eta2}.  This case is in some ways simpler
because one may choose a vector $\vv$ in the cone $B_2$
for which $\rr \cdot \vv = 0$.  Because $B_2 \subseteq \cK^{A,B} (\yy)$
for every $\yy$ where $\qt$ vanishes, setting $\eta (\yy) \equiv \vv$
works everywhere except where $\hht$ vanishes and $\cK^{A,B} (\yy)$
may be smaller.  In a neighborhood of this projective line,
we may instead take $\tilde{\eta} (\yy) \equiv - \yy - c e_3$
for some $c > 0$.  Piecing these together, projectively, via a 
partition of unity, finishes the construction.

Finally, we note that when $\rr$ is on the dashed boundary in
figure~\ref{fig:teardrop}, it is obstructed and the above
construction does not work.

\setcounter{equation}{0}
\section{Evaluation of integrals} \label{sec:integrals}

\subsection{Reduction to integrals of meromorphic forms over
  projective cycles}
\label{ss:replace}

In this section we prove Theorems~\ref{th:no plane} 
and~\ref{th:cone and plane}.  The first step is to show 
that the linearization in Lemma~\ref{lem:exp} may be integrated 
term by term.  With $\dchain$ and the local pieces 
$\{ \dchain (\ww) \}$ defined as in Theorem~\ref{th:eta cone}, 
we recall from~\eqref{eq:contrib-chain} the quantity
$$\contrib (\ww) := \left ( \frac{1}{2 \pi i} \right )^d
\int_{\dchain (\ww)} e^{-\rr \cdot \zz} \frac{p(\zz)}{q(\zz)^s 
   \prod_{j=1}^k h_j (\zz)^{n_j} } 
   \, d\zz  \, .$$
The next step is to replace functions in the integrand by 
an appropriate series of homogeneous functions.  Term by term
integration follows immediately from the expansion~\ref{eq:est 1}
and the big-O estimate (Lemma~\ref{lem:big-O}).
\begin{lem} \label{lem:replace}
Let $F$ satisfy the quadratic point hypotheses.  Let $\dchain (\ww)$ 
be as given in Theorem~\ref{th:eta cone}.  Let $c (\mm , l , n)$ 
be the coefficients of the expansion given in Lemma~\ref{lem:exp} for
the function $f := F \circ \exp$ at the point $\xmax + i \ww$.
Then for every $N \geq 1$ and sufficiently small $\delta > 0$,
\begin{eqnarray} \label{eq:asym proj}
\frac{\ZZ^\rr}{(2 \pi i)^d} \; \contrib (\ww) 
   & = & \sum_{ |\mm| - h \ell - kn < N } 
   \int_{\dchain (\ww)} c(\mm , \ell , n) 
   \yy^\mm \qt (\yy)^{-s-l} \prod_{j=1}^k \hht_j (\yy)^{n_j - n} 
   \, d\yy \\[3ex]
   && + O \left ( |\rr|^{2s - d - N} \right ) 
   \, .  \nonumber
\end{eqnarray}
In the case $k=0$, this reduces to
$$\frac{\ZZ^\rr}{(2 \pi i)^d} \contrib (\ww) \;
   = \sum_{ |\mm| - h \ell < N } 
   \int_{\dchain (\ww)} c(\mm , \ell) \yy^\mm \qt (\yy)^{-s-l} \, d\yy 
   + O \left ( |\rr|^{2s - d - N} \right ) \, .$$
$\Cox$
\end{lem}

\begin{unremark}
The result is also true replacing $\dchain (\ww)$ by 
the infinite chain $\proj^{(\delta)} (\ww)$ defined 
in Definition~\ref{def:conechain} with $A = \ft_\ww
= \homog (f , \ww)$, as the integrals over these
chains differ by a term exponentially smaller than
$|\ZZ^{-\rr}|$.
\end{unremark}

\subsection{Generalized functions} \label{ss:gen func}

The integrands in~\eqref{eq:asym proj} are homogeneous and
the chains of integration projective.  Many such integrals
are evaluated in~\cite{ABG} but there is a hitch:
the integral is evaluated not over $\dchain$ but over
$i \R^d$.  The latter integral is in general not convergent
over $i \R^d$ for two reasons.  First, integrability
will fail near the zeros of $\qt$ whenever $s$ is large.  
Secondly, because $|\exp (\rr \cdot \xx)| = 1$ on $i \R^d$ 
(as opposed to the exponential decay on $\wchain$), 
integrability at infinity will fail whenever $|\mm| \geq 2s - d$.  
These problems are solved respectively by moving the contour
and by inserting compactly supported functions inside the
integral.  The apparatus to do this is the theory of
generalized functions (distributions) and their Fourier transforms,
developed in~\cite{gelfand-generalized} and elsewhere.  We
summarize the facts needed from this literature.

We work with two linear spaces that are dual to each other.  
We call these $\R^d$ and $\R^{d*}$.  We fix bases dual to 
each other so that for $\rr \in \R^{d*}$ and
$\xx \in \R^d$, we have $\rr \cdot \xx = \sum_{j=1}^d r_j x_j$.
While all of the ensuing constructions could be defined on
either space, our purposes require slightly different constructions 
the two spaces and we reduce confusion by developing these 
asymmetrically. 

\label{cs}
Let $\cs$ denote the space of smooth complex valued functions 
on $\R^{d*}$ with compact support.  These are called 
\Em{test functions} in~\cite{gelfand-generalized} and the closed
support of a test function $g$ is denote by $\Sup (g)$.  
Topologize test functions by convergence
of all derivatives; this may be metrized, for example, by
$$||g|| := \sum_n 2^{-n} \sum_{|\kk| = n} \phi \left ( 
   \sup \left | \frac{\partial^\kk}{\partial \rr^\kk} \, g \right |
   \right ) $$
\label{gfspace}
where $\phi (x) = x/(x+1)$.  The space $\gfspace$ of 
\Em{generalized functions} (sometimes called \Em{distributions})
is defined to be the dual of
$\cs$, namely the space of continuous linear functions 
\label{locint}
on $\cs$.  Let $\locint$ be the space of locally integrable
functions on $\R^{d*}$, that is, functions $g$ such that
$g \in L^1 (B_N)$ for the ball $B_N$ of every radius $N$ in $\R^{d*}$.
There is a natural embedding of $\locint$ into $\gfspace$
mapping the function $f$ to the linear map $g \mapsto 
\int f(\rr) g(\rr) \, d\rr$.  We denote by $\imf$ the
image of $f$ under this identification.  Generalized functions
in the image of this identification are called \Em{standard}
functions, but there are many nonstandard functions.  One
example is the function $\delta_\rr$ defined by $\delta_\rr (g) 
= g(\rr)$.  

Sometimes a function is not standard but agrees
with a standard function on some region.  Let $\domain$
be an open set in $(\R^d)^*$ and suppose that for any $g$ 
whose closed support is contained in $\domain$, the value
of the generalized function $\cL$ is given by
$\int_{\domain} f(\rr) g(\rr) \, d\rr$ for some
function $f \in L^1 (\domain)$.  We then say that
$\cL$ is \Em{partially identified} with $f$ on $\domain$.

Differentiation may be defined on $\gfspace$ by
\begin{equation} \label{eq:derivative}
\frac{\partial}{\partial r_j} \gf := 
   g \mapsto - \gf \left ( \frac{\partial}{\partial r_j} g \right ) \, .
\end{equation}
This commutes with the identification map: integrating by parts,
\begin{eqnarray}
\left ( \frac{\partial}{\partial r_j} \imf \right ) (g) & := &
   - \imf \left ( \frac{\partial}{\partial r_j} g \right ) 
   \nonumber \\
& := & - \int f(\rr) \frac{\partial}{\partial r_j} g (\rr) \, d\rr 
   \nonumber \\
& = & \int \frac{\partial}{\partial r_j} f (\rr) g (\rr) \, d\rr 
   \label{eq:derivative 2}
\end{eqnarray}
which is evidently the embedded image of $\partial f / \partial r_j$
applied to $g$.
An example of this is the generalized function 
$(\partial / \partial r_i) \delta_\rr$ which maps
$g$ to $(\partial g / \partial r_j) (\rr)$.  A famous
result (not needed here) is that every generalized 
function is of this form: given $\gf \in \gfspace$,
there is a continuous $f \in \locint$ and a $\kk$ for 
which $\gf = \partial^\kk f / \partial \rr^\kk$.
Restricting the integral to $\domain$, we see that 
differentiation also commutes with partial identification.

\label{rd}
On $\R^d$ we define a slightly different space of test functions.
Denote by $\rd$ the space of rapidly decaying smooth functions,
meaning that they are $O(|\xx|^{-N})$ at infinity for every $N > 0$.
Again, topologize by convergence of all derivatives.
Let $\gfs$ denote the dual of $\rd$.  This space of generalized
functions is slightly smaller than the space $\gfspace$.
Let $\polybd$ denote the space of functions $f$ on $\R^d$
satisfying $|f(\xx)| \leq C (1 + |\xx|)^N$ for some $C,N > 0$.
Then the space $\polybd$ embeds in $\gfs$; again, we denote
the image of $f$ under this identification by $\imf$.

\subsection{Inverse Fourier transforms} \label{ss:fourier}

We now define Fourier transforms and their inverses.  Fourier
transforms will be defined for functions on the dual space,
while inverse Fourier transforms will be defined for functions
on ordinary space.  Fourier transforms will be defined only
for nice functions, while inverse Fourier transforms will be
defined for generalized functions. 

For $g \in \cs$, define the Fourier transform $\hat{g}$ by
\begin{equation} \label{eq:fourier 1}
\hat{g} (\xx) := \int_{\R^{d*}} g(\rr) \exp (-i \rr \cdot \xx) 
   \, d\rr \, .
\end{equation}
Observe that $\hat{g} \in \rd$ (this is the Riemann-Lebesgue Lemma).
In fact, we may extend $\hat{g}$ to a function on all of $\C^d$.
This is a holomorphic function and for every integer $N > 0$
it is shown in~\cite[(2.3)]{ABG} to satisfy an estimate
\begin{equation} \label{eq:estimate}
|\hat{g} (\xx + i \yy)| \leq C(N) (1 + |\xx + i \yy|)^{-N} 
   \exp \left ( \sup_{\rr \in \Sup (g)} \rr \cdot \yy \right ) \, .
\end{equation}

Let $\gf$ be a generalized function in $\gfs$.  We define the
inverse Fourier transform $\IFT (\gf)$ by
\begin{equation} \label{eq:IFT}
\IFT (\gf) (g) := (2\pi)^{-d} \gf (\hat{g}) \, .
\end{equation}
\label{IFT}
This is well defined because we have just seen that $\hat{g} \in \rd$
and it is easy to see that it is continuous and therefore an
element of $\gfspace$.  Suppose that $f \in L^1 (\R^d)$.  Then
\begin{eqnarray*}
\IFT (\imf) (g) & = & (2 \pi)^{-d} \imf (\hat{g}) \\
& = & (2 \pi)^{-d} \int_{\R^d} f(\xx)  
   \left ( \int_{\R^{d*}} g(\rr) \exp (i \rr \cdot \xx) \, d\rr \right )
   \, d\xx \, .
\end{eqnarray*}
Since $|f|$ and $|g|$ are integrable, we may switch the order of
integration to see that $\IFT (\imf)$ is the generalized function
identified with the actual function 
$$(2 \pi)^{-d} \int f(\xx) \exp (i \rr \cdot \xx) \, d\xx \, .$$
Although we do not need it here, we remark that $\IFT$ inverts
the Fourier transform: let $g \in \cs$ so that $\hat{g} \in L^1 (\R^d)$;
then the above computation shows that $\IFT (\hat{g})$ is equal
to the standard function $(2 \pi)^{-d} \int \hat{g} (\xx) 
\exp (i \rr \cdot \xx) \, d\xx$, which is equal to $g$ by the
usual theorem on inverting Fourier transforms.

\subsubsection*{Boundaries of holomorphic functions}

To evaluate the integrals arising in this paper, we must examine
generalized functions arising as limits of holomorphic functions.
Let $f$ be holomorphic in a domain $\R^d + i \Delta$ where 
$\zero$ is contained in the boundary of $\Delta$.  Suppose that 
$f$ satisfies an estimate 
$$\left | f (\xx + i \yy) \right | \leq C |\yy|^{-N} (1 + |\xx|)^N$$
for some $N > 0$.  The estimate~\eqref{eq:estimate} shows that 
the integral
$$\int_{\R^d + i \eta} f(\xx) \hat{g} (\xx) \, d\xx$$
exists and is independent of $\eta \in \Delta$.  The same
is true for $\int_{\R^d + i \eta} f(\xx) h (\xx) \, d\xx$
as long as the estimate~\eqref{eq:estimate} is satisfied 
with $h$ in place of $\hat{g}$.  In particular, this defines
a generalized function in $\gfs$ (see~\cite[(1.3)--(1.5)]{ABG}
and following).  We denote by $\iota_\Delta f$ the 
generalized function this defines; if $f$ has a limit
in $L^1$ as $\eta \to 0$ in $\Delta$ then $\iota_\Delta f$ is 
just this standard function.  
An example is the function $f(\xx) = A(\xx)^{-s}$
for some homogeneous polynomial, $A$.  If $\Delta$ is a cone of
hyperbolicity for $A$, then $f$ is holomorphic on $\R^d + i \Delta$
and blows up no worse than a power of the magnitude of the imaginary
part of the argument.  When $s$ is sufficiently large, this is
not a standard function.  

Two classical and useful results generalize the analogous
well known results for ordinary Fourier transforms.
\begin{pr} \label{pr:FT}
Let $f$ be a function satisfying $f(\xx + i\yy) = O(|\yy|^N)$
for some $N$, as above.  Let $\xx^\mm$ be any monomial and
let $L$ be any linear transformation.  Then the inverse
Fourier transforms of $\xx^\mm f$ and $f \circ L^{-1}$ 
are given respectively by
\begin{eqnarray}
\IFT (\xx^{\mm} f) (\rr) & = & i^{|\mm|}
   \frac{\partial^\mm \IFT (f)}{\partial \rr^\mm} 
   \, \label{eq:x m} ; \\
\IFT (f \circ L^{-1}) (\rr) & = & |L| \IFT (f) (L^* \rr) \, . 
   \label{eq:L}
\end{eqnarray}
\end{pr}

\noindent{\sc Proof:}  Pick any $h \in \cs ((\R^d)^*)$.
Integrals in the following calculation will be over $\R^d + i \joy$
in the $\xx$-domain and over $(\R^d)^*$ in the $\rr$-domain.
Using the definition of Fourier transform in the first line,
calculus in the second, integration by parts in the third line,
Fubini's Theorem in the fourth, and integration by parts once
more, we see that
\begin{eqnarray*}
\int_\xx \xx^\mm f(\xx) \fth (\xx) & = & \int_\xx \int_\rr
   \xx^\mm f(\xx) h(\rr) e^{i \rr \cdot \xx} \, d\rr \, d\xx \\
& = & \int_\xx \int_\rr f(\xx) h(\rr)
   \left ( - i \frac{\partial}{\partial \rr} \right )^\mm
   e^{i \rr \cdot \xx} \, d\rr \, d\xx \\
& = & \int_\xx \int_\rr f(\xx) e^{i \rr \cdot \xx}
   \left (i \frac{\partial}{\partial \rr} \right )^\mm h(\rr)
   \, d\rr \, d\xx \\
& = & \int_\rr \hat{f} (\rr) \left (i \frac{\partial}{\partial \rr}
   \right )^\mm h(\rr) \, d\rr \, d\xx \\
& = & \int_\rr f(\rr) \left (-i \frac{\partial}{\partial \rr}
   \right )^\mm \hat{h} (\rr) \, d\rr \, .
\end{eqnarray*}
The left-hand side of this is $(2 \pi)^d \IFT (\xx^\mm f) (h)$ while
the right-hand side is by definition (see~\eqref{eq:derivative})
equal to $(2 \pi)^d i^m [(\partial / \partial \rr)^\mm \IFT (f)] (h)$,
thus verifying~\eqref{eq:x m}.  

The second assertion of the theorem is directly verified.
Making the coordinate change $\xx = L \xx$ and using 
$\rr \cdot L \xx = (L^* \rr) \cdot \xx$ recovers~\eqref{eq:L}.   
$\Cox$

Suppose the function $E$ on $\R^{d*}$ is not locally integrable.
Then $\iota E$ is not a well defined generalized function.
Nevertheless, $\iota E$ is defined as a partial function: if
$g$ is a function supported on a set where $E$ is locally integrable 
then $\iota E (g) = \int E(\rr) g(\rr) \, d\rr$ is perfectly
well defined.  We wish to conclude that an actual Fourier integral
such as occurs in~\eqref{eq:asym proj} of Lemma~\ref{lem:replace} 
is equal to the locally integrable function of $\rr$ computed 
by~\cite{ABG} as the generalized Fourier transform of the integrand 
in~\eqref{eq:asym proj}.  We therefore require the following lemma. 

\begin{lem} \label{lem:equal}
Let $F$ satisfy the quadratic point hypotheses and let 
$\proj^{(\delta)}$ be as in Definition~\ref{def:conechain}
with $A = \ft_\ww$ for some $\ww \in logcrit$.  Let
$$\term (\xx) := \frac{\xx^\mm}{\qt^s \prod_{j=1}^k \hht_j^{n_j}}$$
be any of the terms in the series expansion $f$ at $\xx + i \ww$
as in Lemma~\ref{lem:exp}.  Suppose that the inverse Fourier transform 
$\IFT (\term )$, defined relative to the domain $- i \uu + \R^d$, is given 
by a partial function $\iota \IFT (\term)$ on the set of non-obstructed 
dual vectors in the dual cone, $\normal$ to $\tan_\xx (B)$.  Let
\begin{eqnarray*}
\psi (\rr) & := & (2\pi)^{-d} 
   \int_{- i \proj^{(\delta)}} e^{- i \rr \cdot \yy} \term (\yy) d\yy  \, \; \\
E(\rr) & := & \iota \IFT (\term) (-\rr) \, .
\end{eqnarray*}
Then
\begin{equation} \label{eq:psi=E}
\psi (\rr) \; = \; E(\rr) \hspace{1.4in}
\end{equation}
for any non-obstructed $\rr$ in the dual cone $\normal$.
\end{lem}

\noindent{\sc Proof:}
This is a matter of moving $- i \proj^{(\delta)}$ to $\fiberchain$ while
introducing appropriate smoothing functions to maintain integrability.
We fix a neighborhood $\nbd$ of $\rr$ of non-obstructed dual vectors
in $\normal$, as in the conclusion of Theorem~\ref{th:projective}, 
whose closure is in the interior of $-B^*$.  We then see that 
\begin{equation} \label{eq:unif exp}
\left | e^{-i \rr \cdot \yy} \right | \to 0
\end{equation}
exponentially fast in $|\yy|$ as $\yy \to \infty$ in $- i \proj^{(\delta)}$,
uniformly as $\rr$ varies over $\nbd$.  It follows that if $g : (\R^d)^*
\to \C$ is smooth and supported on some compact subset of $\nbd$, then
\begin{equation} \label{eq:unif g}
| \hat{g} (\yy)| \leq c ||g|| \exp (- c' |\yy|)
\end{equation}
as $\yy$ varies over $\proj^{(\delta)}$, where $||g|| := \int |g|$ 
and $c$ and $c'$ are positive constants not depending on $g$. 
Now fix $\rr \in \nbd$ and let $g_n$ be a sequence of smooth 
functions supported on $\nbd$ and converging to $\delta_{\rr}$.
Note that the estimate~\eqref{eq:unif g} holds for the function
$g = \delta_\rr$ as well as for all $g_n$, where in this
case $\hat{g} (\xx) = e^{- i \rr \cdot \xx}$.  

To establish~\eqref{eq:psi=E}, fix an $\ee > 0$.  Nonvanishing 
of $\qt$ and each $\hht_j$ on $-i \proj^{(\delta)}$, 
together with~\eqref{eq:unif g}, implies that we may 
pick a compact set $K$ such that 
\begin{equation} \label{eq:triangle 2}
\int_{-i \proj^{(\delta)} \setminus K} 
   \left | \hat{g_n} (\yy) \term(\yy) \right | \, d\yy
   \leq \frac{\ee}{4}
\end{equation}
for all $n$, and also for $\delta_{\rr}$ in place of $g_n$.
The sequence $\hat{g}_n$ converges to $\exp (i \rr \cdot \yy)$
uniformly on $K$, hence we may choose $N_0$ large enough so that
for $n \geq N_0$, 
\begin{equation} \label{eq:triangle 3}
\left | \int_K \term(\yy)
   \left | \exp ( i \rr \cdot \yy) - \hat{g}_n (\yy) \right | \, d\yy
   \right | \leq \frac{\ee}{4} \, .
\end{equation}
Increasing $N_0$ if necessary, we may also ensure that 
\begin{equation} \label{eq:E close}
\left | E(\rr) - \int E(\rr') g_n (\rr') \, d\rr' \right | 
   \leq \frac{\ee}{4}
\end{equation}
for all $n \geq N_0$.  We may now conclude that $n \geq N_0$ implies
\begin{equation} \label{eq:penult}
\left | \psi (\rr) - (2 \pi)^{-d} \int_{- i \proj^{(\delta)}} 
   \hat{g}_n (\yy) \term(\yy) \, d\yy \right | \leq \frac{3}{4} \ee \, .
\end{equation}
Indeed, the two terms we have subtracted are integrals
over $- i \proj^{(\delta)}$ of two integrands; denoting the
integrands by $\beta$ and $\beta'$, we break $- i \proj^{(\delta)}$
into $K \cup K^c$ and use the triangle inequality, viz.,
$$\left | \int \beta - \int \beta' \right | \leq \int_K |\beta - \beta'|
   + \int_{K^c} |\beta| + \int_{K^c} |\beta'|$$
and $(2 \pi)^{-d} < 1$ to obtain~\eqref{eq:penult}.

The homotopy $- i \Phi_t$ in Theorem~\ref{th:projective} 
moves $-i \proj^{(\delta)}$ to $-i \uu + \R^d$ while avoiding 
the singularities of $\term$ (by homogeneity, term is
singular at $\yy$ if and only if it is singular at $i \yy$).
We have seen in~\eqref{eq:estimate} that $\hat{g}_n$ is rapidly
decreasing on the image of this homotopy.  Truncating at the boundary
of a large ball and sending this to infinity shows that the integral
in~\eqref{eq:penult} is unaffected by applying the homotopy.  Hence,
$$(2 \pi)^{-d} \int_{-i \proj^{(\delta)}} \hat{g}_n (\yy) \term(\yy) \, d\yy
   = (2 \pi)^{-d} \int_{\fiberchain} \hat{g}_n (\yy) \term(\yy) \, d\yy 
   = \int E(\rr') g_n (\rr') \, d\rr'$$
by the definition of the inverse Fourier transform.  This identity
allows us to apply the triangle inequality to~\eqref{eq:E close}
and~\eqref{eq:penult}, yielding
$$\left | \psi (\rr) - E (\rr) \right | \leq \ee \, .$$
Since $\ee > 0$ was arbitrary, this proves the lemma.
$\Cox$

\subsection{Proof of Theorem~\protect{\ref{th:no plane}}}
\label{ss:cone only}

We begin with a result from~\cite{ABG}, evaluating the 
Fourier transform of $S^{-s}$ where $S$ is the standard
Lorentzian quadratic $x_1^2 - x_2^2 - \cdots - x_d^2$.
For this special case, we let $\Delta := \{ \xx : x_1 < 0
\mbox{ and } S(\xx) > 0 \}$ be the cone of hyperbolicity
containing the negative $x_1$ axis, we choose an 
element $\eta = -e_1$ of $\Delta$, and we let 
$\normal := \{ \rr : r_1 > 0 \mbox{ and } r_1^2 - 
\sum_{k=2}^d r_k^2 > 0 \}$ be the dual cone to $\Delta$.
In the case where $s$ is not an integer, we also need notation
to specify what is meant by $S(\xx)^{-s}$.  To specify
a branch of this is the same as specifying a branch of
the argument function $\Arg S(\xx)$.  On any simply connected
domain where $S \neq 0$, this may be accomplished by specifying
$\Arg S (\xx)$ at any point in the domain.  Therefore, we write
$\disp{S(\xx)^{-s} |_{\Arg (S(\eta)) = \theta}}$ to denote 
such a specification.  

\begin{thm}[\protect{\cite[Equation (4.20)]{ABG}}] \label{th:riesz}
Let $\symp (\xx) := x_1^2 - x_2^2 - \ldots - x_d^2$, so that
$\symp^* (\rr) = r_1^2 - r_2^2 - \ldots - r_d^2$.
Then the inverse Fourier transform of $\symp^{-s}$ exists in a 
generalized sense and, if $s \neq 0 , - 1, -2, \ldots$, it is given by 
\begin{equation} \label{eq:special riesz}
   e^{i \pi s} \frac{\symp^* (\rr)^{s - (d/2)}}
   {2^{2s-1} \pi^{(d-2)/2} \Gamma (s) \Gamma (s + 1 - (d/2))} .
\end{equation}
To be precise, let $\Delta$ be the component of the real cone 
$\{ S > 0 \}$ that contains the negative $x_1$ axis and let
$\eta \in \Delta$, for example, $\eta = e_1$.  Then if $g$ is 
supported on a compact subset of $\normal$, 
\begin{equation} \label{eq:special riesz explained}
(2 \pi)^{-d} \int_{\R^d + i \eta} \symp^{-s} (\xx) \hat{g} (\xx) 
   = C \int \symp^* (\rr) g(\rr) \, d\rr
\end{equation}
where $C = e^{i \pi s} \,  /  \, [2^{2s-1} \pi^{(d-2)/2} 
\Gamma (s) \Gamma (s+1-d/2)]$.  When the Gamma function is
infinite, the generalized Fourier transform vanishes on
the open cone $\Delta$ (it is supported on $\partial \Delta$).
$\Cox$
\end{thm}

\noindent{\sc Proof:} This result is taken from~\cite{ABG} but
with definitions spread across several sections.  So as to 
make the citation checkable (especially in light of some 
minor errors), we reference a number of passages of~\cite{ABG}.
Equation~\cite[4.13]{ABG} defines a Fourier transform of $S^{-s}$ 
(their notation for $S$ is $a$).
Then in~\cite[4.20]{ABG} they give the following
formula for this quantity, attributed to~\cite{riesz}:
\begin{equation} \label{eq:corrected}
\frac{\symp^* (\rr)^{2 - d/2}}{\pi^{d/2-1} 2^{2s-1} \Gamma (s)
   \Gamma (s+2-d/2)} \, .
\end{equation}
The argument of the $\Gamma$ function the second time is wrong: 
it should be $s+1-d/2$, agreeing with~\cite{riesz}.  They are also 
missing a factor of $e^{i \pi s}$.  To see that this factor should 
be present, note that their specification of the branch of $S^{-s}$
is given at the top of page~146: they specify this over the simply
connected set $i \eta + \R^d$ by specifying that $\Arg S(i \eta)
= \pi + \Arg (S(\eta))$.  Taking $\eta = e_1$, we see that 
$\Arg (S(i \eta))$ must be an odd multiple of $\pi$.  For such a
specification, the Fourier transform cannot be real for small real 
values of $s$.  Indeed, taking $s$ to be small and positive, and 
noting that $S$ is never positive real on $i \eta + \R^d$, we see 
that all arguments of $S(\xx)$ lie between~$2n \pi$ and $2(n+1) \pi$, 
hence all arguments of $S(\xx)^{-s}$ lie between~$2ns\pi$ and 
$(2n+2)s \pi$.  Let $\xx = i \eta + \yy$.  Then switching $\yy$
and $-\yy$ conjugates $e^{i \rr \cdot \xx}$ while reflecting 
$S(\xx)^{-s}$ about the line $\Arg = (2n+1) s \pi$.  Therefore, 
the integrand of the Fourier transform 
$$\int_{i \eta + \R^d} \uu e^{i \rr \cdot \xx} S(\xx)^{-s} \, d\xx$$
is also reflected about the line $\Arg = (2n+1) s \pi$ implying that
the integral must therefore lie on the line of reflection.  For
small values of $s$ this is not real, demonstrating the need for 
a correction.  The corrected formula~\eqref{eq:corrected} 
implies~\eqref{eq:special riesz}.
$\Cox$

\begin{cor}[Fourier transform of a cone] \label{cor:gen riesz}
For any real quadratic $A$ having signature $(1,d-1)$, any 
monomial $\xx^\mm$, and any $s \neq 0 , d/2 - 1$, the 
inverse Fourier transform of $\xx^\mm A^{-s}$ is given by
\begin{equation} \label{eq:riesz}
   e^{i \pi s} i^{|\mm|} \, 
   \frac{|M| (\partial / \partial \rr)^\mm \dual (\rr)^{s - (d/2)}}
   {2^{2s-1} \pi^{(d-2)/2} \Gamma (s) \Gamma (s + 1 - (d/2))} 
\end{equation}
where $M$ is any real linear transformation such that $\form
= \symp \circ M^{-1}$.
\end{cor}

\noindent{\sc Proof:} Pick a linear transformation $M$ such
that $\form = \symp \circ M^{-1}$.  Recall that $\dual (\rr)
= \symp^*(L^* \rr)$.  Use the second part of Proposition~\ref{pr:FT}
and then the first to obtain~\eqref{eq:riesz}.  
$\Cox$

\noindent{\sc Proof of Theorem}~\ref{th:no plane}:
We are required to prove the given asymptotic expansion 
of $\contrib (\ww)$.  We have assumed no linear factors,
so the second, simpler formula from Lemma~\ref{lem:replace}
applies, giving 
$$\contrib (\ww) = \left ( \frac{1}{2 \pi i} \right )^d
   \ZZ^{-\rr} \int_{\dchain (\ww)} e^{-\rr \cdot \yy} 
   \sum_{|\mm| - 2 n < N} c(\mm , n) \yy^\mm \qt (\yy)^{-s-n}
   \; d\yy \;  + O \left ( |\rr|^{-(N+d-2s)} \right )$$
as long as $N > 2s-d$.  

Now integrate term by term.  We see that
\begin{equation} \label{eq:summation}
\contrib (\ww) = \ZZ^{-\rr} \left [ O (|\rr|^{-d+2s-N}) +
   \sum_{|\mm| - 2 n < N} \left [
   \left ( \frac{1}{2 \pi i} \right )^d
   \int_{\dchain (\ww)} e^{-\rr \cdot \yy}
   c(\mm , n) \yy^\mm \qt (\yy)^{-s-n}
   \right ] \; d\yy \right ] \, .
\end{equation}
The specification of $-s$ power for this generating function is
that the argument of $\qt (\uu)$ is zero.  To turn this into a 
Fourier transform, the change of variables $\yy = i \yy'$ is 
needed.  Under this change of variables, $d\yy = i^d \, d\yy'$,
and the summand in~\eqref{eq:summation} becomes
$$ c(\mm , n) (2 \pi )^{-d}
   \int_{- i \dchain (\ww)} e^{- i \rr \cdot \yy'}
   (i \yy')^\mm \qt (i \yy')^{-s-n} \, d\yy' \, ,$$
Now the argument of $\qt (i \yy')$ is still continued from 
initial data $\Arg (\qt (\uu)) = 0$, and this argument may also be
written as $i \pi$ plus the argument of $\qt (\yy')$ continued from 
initial data $\Arg (\qt ( - i \uu)) = - \pi$ as $\yy'$ varies over
$-i \dchain (\ww)$.  The summand in~\eqref{eq:summation} 
now becomes
$$c(\mm , n) (2\pi)^{-d} i^{|\mm|} e^{-i \pi s}
   \int_{-i \dchain (\ww)} e^{-i \rr \cdot \yy'} |\yy'|^\mm 
   \qt (\yy')^{-s-n} \, d\yy' \, .$$
Everything is now lined up.  Lemma~\ref{lem:equal} shows that
this integral, $\psi (\rr)$, is equal to the partial function
$E$ defined by the Fourier transform of the integrand relative
to the domain $-i \uu + \R^d$.  Corollary~\ref{cor:gen riesz}
computes this partial function (recalling that the inverse 
Fourier transform builds in the factor $(2 \pi)^{-d}$) and 
yields the summand
$$c(\mm , n) (-1)^{|\mm|}
   \frac{|M| (\partial / \partial \rr)^\mm \qt^* (-\rr)^{s - (d/2)}}
   {2^{2s-1} \pi^{(d-2)/2} \Gamma (s) \Gamma (s + 1 - (d/2))} \, .$$
Multiplying by the $\ZZ^{-\rr}$ in front of the 
right-hand side of~\eqref{eq:summation} establishes
desired conclusion~\eqref{eq:s power}, still under the 
assumption $N > 2s-d$, which was used to bound the 
remainder term.

Finally, if $N \leq 2s-d$ then the foregoing argument
may be applied with $N$ replaced by the least integer $N'$
greater than $2s - d$.  Each term in the sum with $|\mm| - 2n > N$
is $O(|\rr|^{-d+2s-N})$ by~\eqref{eq:size}; the remainder 
term satisfies this bound as well because it is $O(|\rr|^{-d+2s-N'})$
with $N' > N$.  The theorem is therefore proved for every $N$. 
$\Cox$

\subsection{Extra linear factors give rise to integral operators} 
\label{ss:many planes}

Let $F$ satisfy the quadratic point hypotheses and let $\rr$ 
be a non-obstructed vector in the dual cone to $\normal 
:= \tan_{\xmax} (B)$.  
Equation~\eqref{eq:x m} of Proposition~\ref{pr:FT} has a moral
inverse: just as multiplication by $\xx$ turns into differentiation 
in the $\rr$-domain, division by a linear function in $\xx$ 
should turn into integration in the $\rr$-domain.  This subsection 
proves a theorem along these lines.  However, because 
anti-differentiation is not well defined, the resulting 
formula~\eqref{eq:ders} fails to specify which iterated
anti-derivative will result.  We show that the correct choice 
can be determined under the additional assumption $2s > d+1$.
This is not, however, the case with the bulk of our examples,
whence our alternative analysis in Section~\ref{ss:cone plane}.

Let $L$ be any linear function, with coefficients $L(\xx) 
= a_1 x_1 + \cdots + a_d x_d$.  We may view $L$ as a
vector in $(\R^d)^*$.  The notation $\partial / \partial L$
will be used to denote the differential operator 
$\disp{\sum_{j=1}^d a_j \frac{\partial}{\partial r_j}}$ on $\rr$-space.
Also, $\disp{\frac{\partial^{\nn}}{\partial L^\nn}
:= \prod_{j=1}^k \left ( \frac{\partial}{\partial L_n} 
\right )^{n_j}}$ denotes the corresponding sum of
monomial operators $(\partial / \partial \rr)^{\nn}$.  

\begin{pr} \label{pr:many planes}
Let $F = P / (Q^s \prod_{j=1}^k H_j^{n_j})$ with $n_j$
positive integers, $\qt$ quadratic and each $\hht_j$ 
linear.  Let $\nn = (n_1 , \ldots , n_k)$ be the multiexponent 
of the functions $H_1 , \ldots , H_k$ in the denominator of $F$ and 
denote $\pp := \nn + n \one = (n_1 + n , \ldots , n_k + n)$.  Then
\begin{equation} \label{eq:ders}
\contrib (\ww) = \sum_{|\mm| - 2 \ell - k n < N} 
   \zz^{-\rr} (-1)^{|\pp| + |\mm|} c(\mm , \ell , n)  
   \frac{\partial^\mm}{\partial \rr^\mm} \term_{\ell , n} (\rr) \, ,
\end{equation}
where $\term_\ell$ is an iterated anti-derivative of $\qt^{-s-\ell}$ satisfying
\begin{eqnarray} \label{eq:h}
i^{|\pp|} \left ( \frac{\partial^\pp}{\partial \hht^{\pp}} \right )
   \, \term_{\ell , n} & = & \IFT (\qt^{-s-\ell}) \\[1ex]
& = & e^{i \pi (s + \ell)} i^{|\mm|} \,
   \frac{|M| (\partial / \partial \rr)^\mm \dual (\rr)^{s + \ell - (d/2)}}
   {2^{2(s + \ell)-1} \pi^{(d-2)/2} \Gamma (s + \ell) 
   \Gamma (s + \ell + 1 - (d/2))} \, . \nonumber
\end{eqnarray}
\end{pr}

\noindent{\sc Proof:} 
The proof proceeds analogously to the proof of Theorem~\ref{th:no plane}.
Using the full expansion, Lemma~\ref{lem:exp}, instead
of Lemma~\ref{lem:linearization} leads to the following
generalization of~\eqref{eq:summation}:
\begin{eqnarray} \label{eq:gen sum}
\contrib (\ww) & = & \ZZ^{-\rr} \left ( \frac{1}{2 \pi i} \right )^d
   \sum \int_{\dchain (\ww)} e^{-\rr \cdot \yy}
   c(\mm , \ell , n) \yy^\mm \qt (\yy)^{-s-\ell} 
   \prod_{j=1}^k \hht_j (\yy)^{-n_j - n} \; d\yy  \\[1ex]
&& + O \left (|\zz^{-\rr}| |\rr|^{-d+2s-N} \right ) 
\nonumber
\end{eqnarray}
where the sum is over the finitely many terms with
$|\mm| - h \ell - k n < N$ and terms with 
$|\mm| - h \ell - k n \geq N'$ are seen by the big-O lemma 
to contribute $O(|\zz|^{-\rr} |\rr|^{-d+2s-N'})$.  Again, 
changing variables to $\yy = i \yy'$ shows the summand to be equal to
$$\ZZ^{-\rr} i^{|\mm| - |\pp|} e^{-i \pi s} 
   c(\mm , \ell , n) \IFT \left ( \frac{\yy^\mm}{\qt (y)^{-s-\ell} 
   \prod_{j=1}^k \hht_j(\yy)^{-n_j - n} } \right ) \, .$$
The first conclusion of Proposition~\ref{pr:FT}
identifies this inverse Fourier transform as
an iterated derivative (introducing a factor of $i^{|\mm|}$), 
hence the summand becomes
\begin{equation} \label{eq:ders 2}
\ZZ^{-\rr} e^{-i \pi s} i^{- |\pp|} (-1)^{\ell + |\mm|}
   c(\mm , \ell , n) \frac{\partial^\mm}{\partial \rr^\mm}
   \IFT \left ( \frac{1}{\qt (y)^{-s-\ell} 
   \prod_{j=1}^k \hht_j(\yy)^{-n_j - n} } \right ) \, .
\end{equation}
 From Corollary~\ref{cor:gen riesz} we have
\begin{equation} \label{eq:g_ell}
\IFT \left ( \qt^{-s-\ell} \right )  = 
   e^{i \pi s} (-1)^\ell \, \frac{|M| \dual (\rr)^{s + \ell - (d/2)}}
   {2^{2(s + \ell)-1} \pi^{(d-2)/2} \Gamma (s + \ell) 
   \Gamma (s + \ell + 1 - (d/2))} \, .
\end{equation}
Multiplying the numerator and denominator of $1 / \qt^{s + \ell}$
by $\prod_{j=1}^k \hht_j^{n_j+n}$ and applying once more
the first part of Proposition~\ref{pr:FT}, we see that
$$\IFT (\qt^{-s-\ell}) = i^{|\pp|} \, \frac{\partial^\pp}{\partial \hht^\pp}
   \IFT \left ( \frac{1}{\qt(y)^{-s-\ell} 
   \prod_{j=1}^k \hht_j(\yy)^{-n_j - n} } \right ) \, ,$$
proving the proposition.    $\Cox$

The rest of the work is in determining $h$ from the
derivative~\eqref{eq:h}.  Begin with two classical
regularity lemmas.
\begin{lem} \label{lem:lipschitz}
Let $f \in L^1 (i \Delta + \R^d)$.  Then $\IFT (\iota_\Delta f)$ 
is standard and locally Lipschitz.
\end{lem}

\noindent{\sc Proof:} Let $g$ have compact support in $(\R^d)^*$.  
\begin{eqnarray*}
\IFT (\iota_\Delta f) (g) & := & (2 \pi)^{-d} (\iota_\Delta f) (\hat{g}) \\
& := & (2 \pi)^{-d} \int_{i \eta + \R^d} f(\xx) \hat{g} (\xx) \, d\xx \\
& := & (2 \pi)^{-d} \int_{i \eta + \R^d} f(\xx) \left [ 
   \int_{(\R^d)^*} g(\rr) e^{- i \rr \cdot \xx} \, d\rr \right ] \, d\xx \\
& = & \int_{(\R^d)^*} g(\rr) \left [ (2 \pi)^{-d} \int_{i \eta + \R^d} f(\xx) 
   e^{- i \rr \cdot \xx} \, d\xx \right ] \, d\rr 
\end{eqnarray*}
by Fubini's theorem, since $e^{-i \rr \cdot \xx}$ is bounded
as $\rr$ varies over the support of $g$ and the imaginary part
of $\xx$ varies over any bounded subset of $\Delta$.  This shows that
$\IFT (\iota_\Delta f)$ is the standard function 
$\iota h$ where $h(\rr) = (2 \pi)^{-d} \int_{i \eta + \R^d}
f(\xx) e^{-i \rr \cdot \xx} \, d\xx$.  To check the local Lipschitz
condition on $h$, note that
$$\left | h(\rr) - h(\rr') \right | = (2 \pi)^{-d} \int_{i \eta + \R^d}
   f(\xx) \left ( e^{- i \rr \cdot \xx} 
   - e^{- i \rr' \cdot \xx} \right ) \, d\xx \, .$$
If $\rr , \rr'$ vary over a compact set $K$, and 
$\xx = i \eta + \xi$ then there is a bound independent of $\xi$:
$$\left | e^{- i \rr \cdot \xx} - e^{-i \rr' \cdot \xx} \right |
   \leq C_K |\rr - \rr'| \, ,$$
which implies $|h (\rr) - h(\rr')| \leq C_K \cdot 
||f||_1 \cdot |\rr - \rr'|$.
$\Cox$

We also require the Paley-Wiener Theorem, stated 
as~\cite[Theorem~2.5]{ABG}.  A generalized function is 
said to have support in a closed set $K$ if it annihilates 
test functions vanishing off of $K$.  
\begin{lem}[Paley-Wiener Theorem] \label{lem:PW}
Suppose that $\Delta$ contains the convex cone $K$.
Then the support of $\IFT (\iota_\Delta f)$ is contained
in the negative dual cone $-K^*$.
$\Cox$
\end{lem}

With these in hand, let $\cKd$ be a connected component of
the non-obstructed subset of $\normal$.  Fix $\rr \in \cKd$.  
Suppose that for each $1 \leq j \leq k$, there is a line
segment $\{ \rr + \lambda L_j : \lambda \in [0 , \lambda_*] \}$ 
(where $\lambda_*$ could be negative) such that $\rr + \lambda_* L_j$
is on the boundary of $\normal$ and $\rr + \beta \lambda^* L_j$
is in $\cKd$ for all $0 \leq \beta < 1$.  In other words, 
for each $j$, traveling from $\rr$ in the directions $\pm L_j$, 
we come to the boundary of $\cKd$ at the same time as 
we come to the boundary of $\normal$.  Define the integral
operator $I_j$ on functions on $\rr$-space by 
$$I_j (g) (\rr) = \int_{\lambda_*}^0 g(\rr + \lambda L_j) 
   \, d\lambda \, .$$
Let ${\bf I}^\pp$ denote the composition over $j$ of powers $I_j^{p_j}$.
We then have the following result.

\begin{thm} \label{th:general}
Under the above geometric conditions on the component
$\cKd$ of the non-obstructed set of $\normal$, if
$2s > d+1$, then $h(\rr)$ in~\eqref{eq:ders} is given by
$$i^{-|\pp|} \; {\bf I}^\pp \; \left [ \IFT ( \qt^{-s-\ell}) \right ]$$
where $\pp = \nn + n \one$ as in Proposition~\ref{pr:many planes}.
\end{thm}

\noindent{\sc Proof:} Under the condition $2s > d+1$, 
one has $|\qt(\xx)^{-s}| = O(|\xx|^{-d-\ee})$ for some $\ee > 0$.
Hence the function $\qt^{-s}$ is integrable away from its poles,
as is therefore $\qt^{-s-\ell} \prod_{j=1}^k \hht_j^{-n_j - n}$.  
Moving to $i \eta + \R^d$, we avoid all poles and hence
$\qt^{-s-\ell} \, \hht^{\pp} \in L^1$.  By
Lemma~\ref{lem:lipschitz}, the Fourier transform is a
standard, locally Lipschitz function, hence continuous.  
The domain of analyticity of $f$ contains every cone
whose closure is in the interior of $\tan_{\xmax} (B)$.
By the Paley-Wiener Theorem, therefore, the inverse Fourier transform
vanishes outside of the closed negative dual cone, $\normal$.
Each differential operator $\partial / \partial L_j$
in~\eqref{eq:h} may now be inverted uniquely, due to the 
boundary condition of vanishing at $\rr + \lambda_* L_j$.  
The unique inverse is $I_j$.  Together with~\eqref{eq:h},
this proves the theorem.   $\Cox$. 

\begin{unremark}
Without the hypothesis $2s > d+1$ we still have 
$$h(\rr) = i^{-|\pp|} \; {\bf \tilde{I}}^\pp \; 
   \left [ \IFT ( A^{-s-\ell}) \right ]$$
where $\tilde{I}_j$ are anti-derivative operators whose boundary
conditions are not determined by continuity from the 
Paley-Wiener Theorem.
\end{unremark}

\subsection{Proof of Theorem~\protect{\ref{th:cone and plane}}}
\label{ss:cone plane}

We have seen that a linear factor $L(\xx)$ in the denominator 
corresponds to a convolution with a Heaviside function, or
equivalently, to an integral operator $I_{L}$.  The integrability
hypothesis in this result is unfortunately somewhat restrictive,
ruling out, for example, the case $s=1, d=3$.
Moreover, from a computational viewpoint, it is not desirable
to have the answer represented as an (iterated) integral.  It
is therefore worth exploring a general method for reducing the
dimension of the integral in question.  In~\cite{ABG}, homogeneity
of the integrand is exploited: integrating out the radial part,
the Fourier transform is reduced to an integral over a cycle in 
$(d-1)$-dimensional projective space, which will be either the
Leray cycle or the Petrovsky cycle.  To this device, we add a
residue computation that further reduces the dimension by one.
Evaluation of the resulting one-dimensional integral 
leads to Theorem~\ref{th:cone and plane}.  Computations
in projective space rely on some standard constructions and
notational conventions which we now introduce.

Let $\pi: \C^d \setminus \{ \zero \} \to \CP^{d-1}$ be the 
projection map.  Any meromorphic form $\omega$ on 
$\CP^{d-1}$ pulls back to a form $\pi^* \omega$ 
on $\C^d \setminus \{ \zero \}$.  The pullback $\pi^*$ is 
one to one onto its range.  It is well known 
e.g.,~\cite[page~409]{griffiths-harris-principles},
that the range is the set of all meromorphic $(d-1)$-forms
on $\C^d$ whose contraction with the Euler vector field
$\sum x_i \partial / \partial x_i$ is zero and that are
homogenous of degree zero.  Here, the degree of $f dx_{j_1}
\wedge \cdots \wedge dx_{j_k}$ is $\deg f + k$ and for
$(d-1)$-forms, those forms killing the Euler field are in a 
one-dimensional subspace of the $(d-1)$-dimensional cotangent
space at each point.  Forms on $\CP^{d-1}$ have no natural
names of their own, so we name them by identifying with their 
pullbacks to $\C^d$, as is done in~\cite{ABG} and elsewhere.  
For computational purposes, when integrating over a chain
$\CC$ in $\CP^{d-1}$, we usually use an elementary chart map
from a slice of $\C^d$, such as $\pi$ restricted to 
$(z_1 , \ldots , z_{d-1} , 1)$.  If $f_j$ are homogeneous of 
degree $1-d$, for example, pulling back by this chart map yields
$$\int_\CC \sum_{j=1}^d f_j (\zz) \, dz_1 \wedge \cdots \wedge
   \widehat{dz_j} \wedge \cdots \wedge dz_d
   = \int_{\CC'} f_d \, dz_1 \wedge \cdots \wedge dz_{d-1}$$
where $\CC'$ is the (unique) lifting of $\CC$ to the 
(simply connected) slice.

The proof of Theorem~\ref{th:cone and plane} begins analogously to
the proof of Theorem~\ref{th:no plane}.  Use the general expansion 
in Lemma~\ref{lem:exp}, just for the leading term, to write 
$$\frac{p (\xmax + i \ww + \yy)}{q (\xmax + i \ww + \yy) 
   h (\xmax + i \ww + \yy)} 
   = \frac{p(\zz)}{\qt (\yy) \hht (\yy)} + R$$
where $R = O(|\yy|^{-2})$ on $\wchain$.  Use Lemma~\ref{lem:replace}
and the big-O lemma to see that 
$$\contrib (\ww) = \frac{\ZZ^{-\rr} P(\ZZ)}{(2 \pi i)^3}
   \int_{\proj^{(\delta)}} \exp (- \rr \cdot \yy) 
   \frac{1}{\qt (\yy) \hht (\yy)} \, d\yy \, + O(|\rr|^{-1}) .$$
Changing variables by $\yy = i \yy'$ and noting that
$d\yy / [\qt (\yy) \hht (\yy)] = d\yy' / [\qt (\yy') \hht (\yy')]$ gives
\begin{eqnarray*}
\contrib (\ww) & = & \frac{\ZZ^{-\rr} P(\ZZ)}{(2 \pi i)^3} 
   \int_{-i \uu + \proj^{(\delta)}} 
   \exp (- i \rr \cdot \yy) \frac{1}{\qt (\yy) \hht (\yy)} \, d\yy 
   \, + O(|\rr|^{-1}) \\[1ex]
& = & \ZZ^{-\rr} P(\ZZ) \, i^{-3} \, 
   \IFT \left ( \frac{1}{\qt \hht} \right ) \, + O(|\rr|^{-1}) \, .
\end{eqnarray*}
Comparing to~\eqref{eq:atn}, it suffices, therefore, to show that
\begin{equation} \label{eq:suffices 1}
   \IFT \left ( \frac{1}{\qt \hht} \right )
   = i^3 \, \frac{\dblres}{\pi} \arctan \left ( 
   \frac{\sqrt{\qt^* (\rr) \qt(\hht^*)}}{\hht^* (\rr)} \right ) \, .
\end{equation}

\subsubsection*{The Leray and Petrovsky cycles}

We are left to compute the inverse Fourier transform of 
$1/(\qt \hht)$.  The first part of this computation, reducing 
to the Leray cycle, is valid for any hyperbolic polynomial
in any non-obstructed direction, so we do it in this generality.
Let $P/H$ be the ratio of two homogeneous polynomials and assume 
$H$ is hyperbolic.  Later we will specialize to the case where 
$H = \qt \hht$.  Denote by $d_* := \deg H - \deg P - d$ the 
inverse degree of homogeneity of the form $(P/H) \, d\zz$.  

Let $B$ be a cone of hyperbolicity for $H$ and fix $\uu \in B$.  
Fix a non-obstructed vector $\rr \in \normal := - B^*$.  Recall 
from our homotopy constructions that there is a vector field $\eta$ 
on $\R^d$ with the following properties.
\begin{enumerate}
\item $\eta$ is homogeneous of degree $+1$.
\item There is a 1-homogeneous homotopy $\{ \eta_t : 0 \leq t \leq 1 \}$
   between $\eta_0 \equiv \uu$ and $\eta_1 = \eta$ such that
   for all $t$ and all nonzero $\yy$, $H(i \yy + \eta_t (\yy) \neq 0$.
\item For all $y \neq 0$, $\rr \cdot \eta (\yy) = 0$.
\end{enumerate}
Indeed, a similar homotopy with the third condition replaced by 
$\rr \cdot \eta (\yy) < 0$ is constructed in Section~\ref{sec:homotopies}.
Stopping the homotopy at the instant, depending on $\yy$, that it 
crosses the hyperplane orthogonal to $\rr$, yields the desired $\eta$
along with a homotopy as prescribed.

Let $S_+$ denote the hemisphere $\{ \yy \in \R^d : |\yy| = 1, 
\rr \cdot \yy \geq 0 \}$.  Let $S_-$ denote the other
hemisphere, where $\rr \cdot \yy \leq 0$.  Let cycles
$\sigma_{\pm}$ be (singular triangulations of) $S_\pm$ oriented
in such a way that $\partial (\sigma_+ + \sigma_-) = 0$.  
Then $\sigma := \sigma_+ - \sigma_-$ is a $(d-1)$-chain
supported on $S^{d-1}$ whose boundary is supported on the 
equator $\{ \yy \in S^{d-1} : \rr \cdot \yy = 0 \}$.

The map $\phi$ defined by $\phi (\yy) := i \yy + \eta (\yy)$ 
induces a covariant map $\phi_*$ on cycles and
homology.  The chain $\phi_* (\sigma)$ maps to a cycle in $\C^d$
with boundary in the complex hyperplane $X^\rr := \{ \zz : 
\rr \cdot \zz = 0 \}$.  Hence $\phi_* (\sigma)$ represents
a homology class in $H_{d-1} (\C^d , X_\rr)$.  For any
homogeneous set $W \subseteq \C^d$, let $\overline{W}$
denote the projection $\pi W$ of $W$ to $\CP^{d-1}$.  The sets
$\sing_\qt$ and $X_\rr$ are homogenous, therefore the pair
$(\C^d - \sing_\qt , (\C^d - \qt) \cap X_\rr)$ projects radially
under $\pi$ to $(\CP^{d-1} - \overline{\sing_\qt} , 
(\CP^{d-1} - \overline{\sing_\qt}) \cap \overline{X_\rr})$.\footnote{The bars denoting projective varieties are about to proliferate; we apologize for the mess, but we tried dropping them but we became confused about which varieties were projective and which were affine.}
\begin{defn}[Leray and Petrovsky cycles] \label{def:leray}
\label{leray}
The chain $\alpha = \alpha (\rr) := \pi \phi (\sigma)$ 
is called the Leray cycle and its class $\pi_* \phi_* [\sigma]
\in H^{d-1} (\CP^{d-1} - \overline{\sing_\qt} , 
(\CP^{d-1} - \overline{\sing_\qt}) \cap \overline{X_\rr})$
is called the Leray class.  

\label{petrovsky}
The boundary of $\alpha$ is a cycle $\beta$ representing a class 
in $H_{d-2} (\CP^{d-1} - \overline{\sing_\qt}) \cap \overline{X_\rr})$. 
Define the Petrovsky cycle $\gamma$ to be a tubular neighborhood
around $\beta$ orthogonal to $X_\rr$.  This is the image under $\pi$
of a cycle supported on $\{ \yy : |\rr \cdot \yy| = \ee \}$ and
avoiding $\sing_\qt$.
\end{defn}
The following result is proved in~\cite[Theorem~7.16]{ABG};
here we correct a typo: the second appearance of $\chi_q^0$,
namely the one in ($7.17'$), should be just $\chi_q$; see 
equation~(1.6) on page 122 of~\cite{ABG}.  
Define a $(d-1)$-form $\omega$, killing the Euler vector field
and having homogeneous degree $d$, by
$$\omega := \frac{1}{d} \sum_{j=1}^d (-1)^{j+1} z_j \, dz_1 
   \wedge \cdots \wedge \widehat{dz_j} \wedge \cdots \wedge dz_d \, .$$
To explain what is about to appear 
in~\eqref{eq:leray}~--~\eqref{eq:petrovsky} below, we must see why
the integral of a meromorphic 2-form on $\CP^2$ over 
a relative homology class in $H_2 (\CP^2 , \M)$ is well defined. 
Indeed, integration over a relative homology class with respect
of a complex submanifold of positive co-dimension is always well 
defined, for the following reason.  Let $\CC$ be a representing cycle 
for the class, that is a 2-chain with boundary in $\M$.  The integral 
over any other relative cycle differs from the integral over $\CC$
by the integral over a relative boundary, a relative boundary 
being an absolute boundary plus something in $\M$.  
Since $d\omega = 0$ for any meromorphic 2-form on $\CP^2$,
the integral over the boundary vanishes, and since $\M$ has
positive complex co-dimension, the second part of the integral
vanishes as well.

\begin{thm}[Reducing Fourier integrals to Leray/Petrovsky cycles]
\label{th:leray petrovsky}
Let $P/H$ be hyperbolic and fix $\uu$ in a cone $B$ of hyperbolicity 
for $H$ as above.  Let $d_*$ be the inverse degree of homogeneity 
of the form $(P/H) \omega$, that is, $d_* = \deg H - \deg P - d$.  
Let $\alpha$ be the Leray cycle and $\gamma$ be the Petrovsky cycle.  
If $d_* \geq 0$ then 
\begin{equation} \label{eq:leray}
\IFT \left ( \frac{P}{H} \right )
    = \frac{i^{d_*+1}}{(2 \pi)^{d-1} d_*!} \int_\alpha
   (\rr \cdot \zz)^{d_*} \frac{P}{H} \, \omega \, ,
\end{equation}
while if $d_* < 0$ then 
\begin{equation} \label{eq:petrovsky}
\IFT \left ( \frac{P}{H} \right )
   = \frac{i^{d_*}}{(2\pi)^d (|d_*|-1)!} \int_\gamma
   (\rr \cdot \zz)^{d_*} \frac{P(\zz)}{H(\zz)} \, \omega \, .
\end{equation}
$\Cox$
\end{thm}
\begin{unremarks}
$(i)$ Note that the introduction of the factor $(\rr \cdot \zz)^{d_*}$
makes the integrand 0-homogeneous, which is exactly what we need
to interpret it as a form on $\CP^{d-1}$. 
{$(ii)$} In the case $d_* < 0$, the integrand contains a negative
power of $\rr \cdot \zz$.  Let $\Res (\omega)$ be 
$|d_*|^{th}$ residue of the integrand $(\rr \cdot \zz)^{d_*} (P/H) 
\, \omega$ along the projective hyperplane $\overline{X_\rr}$.  
The product structure in the Petrovsky cycle immediately reduces 
the integral one dimension further to $\int_\beta \Res (\omega)$.  
In the case $d_* \geq 0$, the integral does not localize to 
the boundary cycle $\beta$ and one must work harder to kill 
one more dimension.  
\end{unremarks}

\subsubsection*{Residue reduction}

The second step is to reduce by one further dimension
via a residue computation.  The first half of this step
still works in any dimension.  Begin by observing:
\begin{lem} \label{lem:contract equator}
The homology group $H_{d-1} (\CP^{d-1} , \overline{X_\rr})$ vanishes.
\end{lem}

\noindent{\sc Proof:}  $\CP^{d-1} \cap X_\rr$ is homeomorphic to 
$\CP^{d-2}$.  If $p \leq 2(d-2)$ then this inclusion induces
an isomorphism $H_p (\CP^{d-1} \cap X_\rr) \to H_p (\CP^{d-1})$,
where both groups have rank~1 if $p$ is even and vanish otherwise.
It follows that the first and last arrows are isomorphisms 
in the exact sequence
$$H_{d-1} (\CP^{d-1} \cap X_\rr) \to H_{d-1} (\CP^{d-1}) \to 
   H_{d-1} (\CP^{d-1} , X_\rr) \to 
   H_{d-2} (\CP^{d-1} \cap X_\rr) \to H_{d-2} (\CP^{d-1})$$
hence the middle group vanishes.
$\Cox$

We now make use of the Thom isomorphism to ``pass $\alpha$
through $\overline{\sing_H}$'' and obtain a $(d-2)$-chain,
whose tubular neighborhood is homologous to $\alpha$.  The
details are as follows.  By Lemma~\ref{lem:contract equator} 
any chain representing the Leray class is a boundary of a $d$-chain
$\beta\in(\CP^{d-1} , \overline{X_\rr})$.
We may therefore choose a $d$-chain $C$ in $\CP^{d-1}$ whose 
boundary is $\alpha$ plus something in $\overline{X_\rr}$
(because $\partial C$ is a cycle and $\alpha$ is not, the
part in $\overline{X_\rr}$ will be nonzero).  Perturbing $C$ if
necessary, we can assume that $C$ intersects $\overline{\sing_H}$ 
transversely.  The dimension of the intersection of the $d$-chain
$C'$ with the surface $\overline{\sing_H}$ having co-dimension~2
is a $(d-2)$-chain $\thom$, whose orientation is prescribed by
the orientations of $C, \CP^{d-1}$ and $\overline{\sing_H}$.
The chain $\thom$ has boundary in $\overline{\sing_H} \cap 
\overline{X_\rr}$, and is therefore a relative cycle in
$(\overline{\sing_H} , \overline{\sing_H} \cap \overline{X_\rr})$.

Now we are at a point where we require the dimension to be~3.
The chain $C$ has (real) dimension~3 and the surface 
$\overline{\sing_H}$ has dimension~4.  The projective
variety $\overline{\sing_H}$ may not be smooth, but its
singular set has complex co-dimension at least~1, hence
has real dimension at most~2.  Generically perturbed
$C$ therefore does not intersect the singular set of
$\overline{\sing_H}$, and hence $\thom$ is supported on
the set of smooth points of $\overline{\sing_H}$.  We may 
define the tubular neighborhood $T(\thom)$ supported 
on the set $\{ |H| = \ee \}$, which is locally a product of 
$\thom$ with a small circle about the origin in $\C^1$ with 
the standard orientation.  

Integration around this circle is computed by taking a residue.  
We recall a definition of the residue form on any complex space.
Let $\theta$ be a meromorphic form with a pole on the set 
$\overline{\sing_H}$, the pole being simple except on a 
proper subvariety $\overline{\sing_H \cap \sing_K}$.  
Then the residue is defined as follows.  Write $\theta 
= H \cdot \omega$ where $\omega$ is holomorphic away from 
on $\overline{\sing_K}$.  
\label{Res}
Define $\Res [\theta , \overline{\sing_H}]$ to be the unique 
form on $\overline{\sing_H}$ that satisfies
$$\Res [\theta , \overline{\sing_H}] \wedge dH = \omega$$
away from $\overline{\sing_K}$.  In coordinates, the residue of 
$(G/H) dz_1 \wedge \cdots \wedge dz_d$ is given by 
$(G / (\partial H / \partial z_1)) dz_2 \wedge \cdots \wedge dz_d$; 
(with $z_j$ in place of $z_1$, we have the alternative expression
$(-1)^{j-1} (G / \partial H / \partial z_j) dz_1 \wedge \cdots
\wedge dz_{j-1} \wedge dz_{j+1} \wedge \cdots \wedge dz_d$).
The following well known result may be demonstrated by 
expressing the integral over the tube as an iterated integral, 
first around a circle.
\begin{lem} \label{lem:thom}
Suppose $d=3$.
The Leray class $\alpha$ is homologous to the tube $T(\thom)$
around $\thom$.  Consequently, for any meromorphic form $\theta$ 
on $H$ with a simple pole at $H$, 
$$\int_\alpha \theta = \int_{T(\thom)} \theta 
   = (2 \pi i) \int_{\thom} \Res [\theta ; \overline{\sing_H}] \, .$$
In particular, when $d_* = 0$ in Theorem~\ref{th:leray petrovsky},
putting this together with~\eqref{eq:leray} and specializing to 
$H = \qt \hht$ and $\theta = \omega / (\qt \hht)$ yields
$$\IFT \left ( \frac{1}{\qt \hht} \right )
   = \frac{i}{(2 \pi)^2} \int_{\alpha} \frac{\omega}{\qt \hht}
   = \frac{-1}{2 \pi} \int_{\delta} \omega_L$$
where
$$\omega_L := \Res \left [ 
   \frac{\omega} {\qt \hht} ; \overline{\sing_{\qt \hht}} 
   \right ]  \, .$$
$\Cox$
\end{lem}

\begin{unremark}
In higher dimensions, it may happen that $\thom$ intersects the
singular set of $\overline{\sing_H}$.  In that case, one might
expect a version of Lemma~\ref{lem:thom} showing the Leray cycle
to be homologous to the sum of a tubular neighborhood of 
$\thom$ away from the singular set and a cycle supported on 
a neighborhood of the singular set.  
\end{unremark}

\subsubsection*{The special case of a cone and a plane}

For the remainder of this section, we specialize to the case in 
Theorem~\ref{th:cone and plane}.  In addition to $d=3$, 
we suppose $k=s=1$, whence the inverse degree, 
$d_*$, of homogeneity is zero.  We further specialize to the case
where $\sing_\qt \cap \sing_\hht \cap \R^d \neq \emptyset$ and
where the variety $\overline{\sing_{\qt \cdot \hht}}$ has precisely 
two points, each of multiplicity one.  Thus we are in the case of 
Section~\ref{ss:cone and plane}, which arises in several of our 
applications and is illustrated in figure~\ref{fig:cubeplot}.  
The normal cone is shaped as in figure~\ref{fig:teardrop} and the 
vector field $\eta_1$ constructed in Section~\ref{ss:cone and plane} 
is one we can use to construct the Leray cycle.  
In this case ($d=3$), since $\delta$ avoids $\overline{\sing_\hht}$, 
we may restrict to $\overline{\sing_\qt}$ and write 
$$\omega_L := \Res \left [ \frac{\omega}{\qt \hht} ; \overline{\sing_\qt} 
   \right ] \, .$$

Let $p$ be a common zero of $\qt$ and $\hht$.  
The space of 2-forms on $\CP^2$ is one-dimensional, 
hence at $p$, the form $\omega$ is a 
multiple of the form $d\qt \wedge d\hht$.   The double residue 
\label{dblres 2}
$$\dblres := \Res \left [ \frac{\omega}{\qt \hht} , 
   \overline{\sing_\qt} \cap \overline{\sing_\hht} \right ]$$ 
is the value of this ratio, that is, 
$\omega = \dblres (p) \, d\qt \wedge d\hht$ at $p$.  We have the following 
explicit description of the residue form $\omega_L$.

\begin{lem} \label{lem:t}
Let $t : \CP^1 \to \overline{\sing_\qt}$ be any local 
parametrization.  Then 
\begin{eqnarray*}
\omega_L & = & \dblres (t_3) \, \cdot \,  \left ( \frac{dt}{t - t_3} 
   - \frac{dt}{t - t_4} \right ) \\
& = & \dblres (t_4) \, \cdot \,  \left ( \frac{dt}{t - t_4} 
   - \frac{dt}{t - t_3} \right ) 
\end{eqnarray*}
where $t_3$ and $t_4$ are the values of the parameter $t$ for which
$\hht$ vanishes.  
\end{lem}

\noindent{\sc Proof:} The form $\omega_L$ is meromorphic on 
$\overline{\sing_\qt}$ with precisely two simple poles.  Therefore, it
may be written as $C [dt/(t - t_3) - dt / (t - t_4)]$.  Taking 
the residue at $t=t_3$ yields $\Res (\omega_L ; t_3) = C$.  
But iterated resides are the same as multiple residues, hence
\begin{eqnarray*}
C & = & \Res (\omega_L ; t_3) \\
& = & \Res \left [ \Res \left ( 
   \frac{\omega}{\qt \hht} ; 
   \overline{\sing_\qt} \right ) ; \overline{\sing_\hht} \right ] \\
& = & \Res^{(2)} \left [ \frac{\omega}{\qt \hht} ; 
   \overline{\sing_\qt \cap \sing_\hht} \right ]  \\
& = & \dblres \, .
\end{eqnarray*}

To carry out the rest of the computation, let us choose
coordinates for $\C^3$ in which $\qt = x^2 - y^2 - z^2$.
Although these are unrelated to the coordinates in 
which $\hht$ and $\rr$ are described, we will use them 
to compute a coordinate-free description of the
integral.  

The cones of hyperbolicity of $\qt$ are the two components of 
$x^2 > y^2 + z^2$, one containing the negative $x$-axis and one 
containing the positive $x$-axis.  Recall that the cones of 
hyperbolicity for $\qt \hht$ are these cones, bisected by 
the plane $\hht = 0$.  Recall we have fixed $\uu \in B$, where $B$ 
is one of these sliced cones and its dual has a teardrop shape.  

\begin{figure}[ht!]
\hspace{1.8in} \includegraphics[scale=0.40]{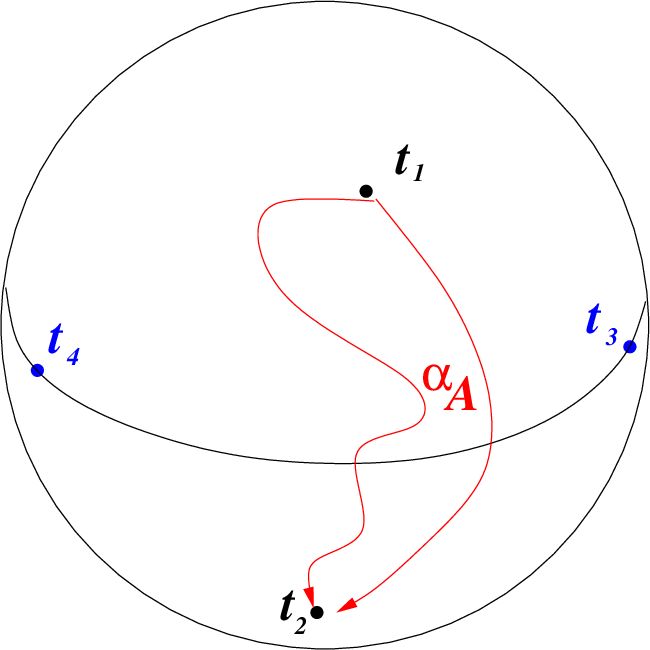}
\caption{the topological sphere $\overline{\sing_\qt}$, its real
part (the equator) and its intersections with the planes $\{ \hht = 0 \}$
and $\overline{X_\rr}$}
\label{fig:rel_cycle}
\end{figure}

The space $\overline{\sing_\qt}$ is a quadratic curve in $\CP^2$.  
We may choose the explicit parametrization 
$t : \CP^1 \to \overline{\sing_\qt}$ by 
$$\overline{(x , tx)} \mapsto \overline{(t^2 + 1 , 2t , t^2 - 1)}$$
in the $(x,y,z)$-coordinate system.
Topologically, $\CP^1$ is a sphere with $\RP^1$ as its equator,
and our parametrization has  the nice feature that the copy of 
$\RP^1$ inside $\CP^1$ maps to the real part of the quadric 
$\overline{\sing_\qt}$.  Thus figure~\ref{fig:rel_cycle} depicts
$\overline{\sing_\qt}$ as a sphere and the real part as the
equator.  The points of $\overline{\sing_\qt}$ where $\rr \cdot \xx$
vanishes are denote $t_1$ and $t_2$ and the points where 
$\hht$ vanishes are denoted by $t_3$ and $t_4$.  
The assumption that $\sing_\hht$ and $\sing_\qt$ intersect in real
space is equivalent to $\hht$ lying outside the dual cone, which
is equivalent to $t$ having real roots.  Thus $t_3$ and $t_4$
are shown on the equator in figure~\ref{fig:rel_cycle}.
Note that the two arcs into which $t_3, t_4$ separate the equator
differ in their position with regard to the cone of hyperbolicity: 
one of them bounds it (we will call this arc {the active} one), 
while the other not.  

If $\rr$ is in the dual cone then $t_1$ and $t_2$
will be complex conjugates, while if $\rr$ is in the pointy
region of the teardrop then $t_1$ and $t_2$ will be real.

\subsubsection*{Case 1: $t_1$ and $t_2$ are complex}

We are nearly ready to evaluate the integral, but we need first
to understand the cycle $\thom$.
The intersection class $\thom$ is a relative cycle 
in $(\overline{\sing_\qt} , \overline{\sing_\qt} \cap \overline{X_\rr})$.
Thus we may draw a representative of this class as a path,
beginning and ending in the set $\{ t_1 , t_2 \}$.  The
meromorphic residue form $\omega_L$ is holomorphic away from
$t_3$ and $t_4$ where it has simple poles.  The integral of 
this form over $\delta$ is therefore determined by combinatorial 
invariants of $\delta$: the positions of the endpoints and the 
number of signed intersections with the two equatorial arcs 
bounded by $t_3$ and $t_4$.  

\begin{lem} \label{lem:cycle}
The homology class of $\thom$ in $(\sing_\qt \setminus \{ t_3 , t_4 \} ,
\{ t_1 , t_2 \})$ is that of an oriented path from $t_2$ to $t_1$,
intersecting one equatorial arc exactly once.
\end{lem}

\noindent{\sc Proof:}  It is shown in~\cite{ABG} (see in particular
figure~6b there 
and the paragraph preceding it) that one can find a representative of
the Leray class such that its boundary  (in
$\overline{X_\rr}-\overline{\sing_\qt} \cap \overline{X_\rr}$)  
is localized near the {\it
  complex points of ${\sing_\qt} \cap \overline{X_\rr}$}, i.e. in our situation 
is the sum of small circles around
the complex zeros of $\qt$ in $\overline{X_\rr}$, oriented according to their
imaginary parts (note that $\hht$ has no non-real zeros there).  
It follows that the boundary of the (relative) cycle $\thom$ is given
by
$$
\partial \delta = [t_1] - [t_2]
$$ 
where $t_1$ has the positive imaginary part, and $\delta$ is the
claimed path, plus one or more absolute cycles (i.e. oriented closed loops)
in ${\sing_\qt}$ and  ${\sing_\hht}$.

To find the homology classes represented by these arcs and loops, we
recall the definition of the Leray class: the vector field $\eta$ constructed
at the end of section \ref{ss:cone and plane} is restricted to the
unit sphere, defining the relative class, the sum of oppositely
oriented hemispheres separated by the hyperplane $X_\rr$.
To evaluate $\thom$ we find a homotopy shrinking these hemispheres 
to a point, keeping their boundaries in $X_\rr$ and tracking where 
the resulting 3-chain hits $\sing_\qt$ and $\sing_\hht$.
This homotopy proceeds in two stages: first we take the linear
homotopy of (the restriction to the unit sphere of) $\eta$, the
vector field constructed in Section~\ref{ss:cone and plane} to the
(restriction to the unit sphere of the) constant vector
field $\xx$ defined in the same place.  In the second stage we 
collapse the sphere to a point, keeping the constant vector field. 

\begin{figure}
  \begin{center}
    \includegraphics[height=1.6in]{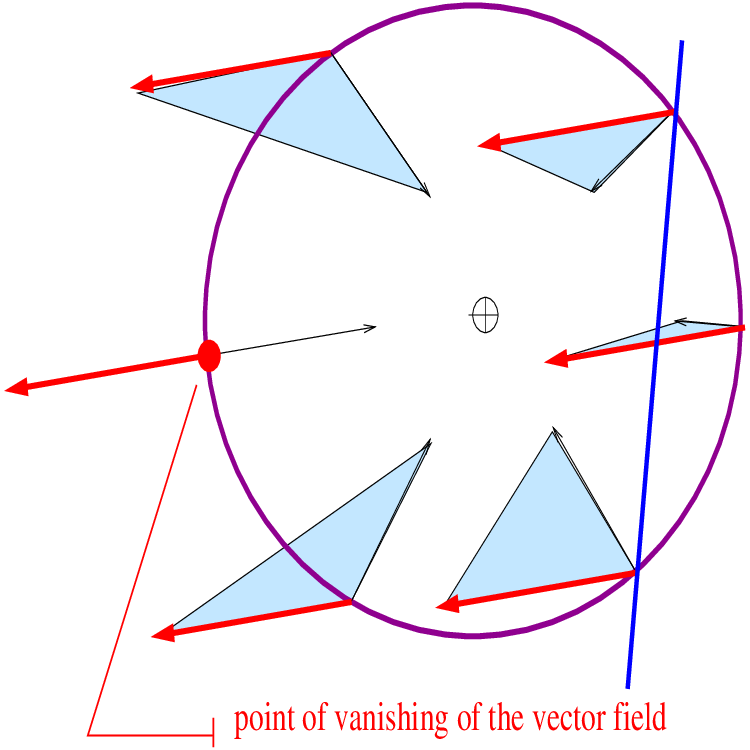}
  \end{center}
\caption{showing where the homotopy intersects the quadric}
\label{fig:intersection count}
\end{figure}

We note first that in ${\sing_\hht}$ this deformation yields the 
empty set: at no instant are the deformed vectors tangent to $\hht$, 
which would be necessary if the deformation were to intersect $\hht$.
This is not the case for $\qt$, and indeed, we know already that the
boundary of $\thom$ there is nontrivial.  The class of $\thom$ is
completely determined by the index of intersection with the active arc
between $t_3$ and $t_4$.  The intersection number of $\thom$ with the
active arc is just the number of points in the real part of
$\sing_\qt$ where our deformation results, at some instant of the
homotopy, in a vector field tangent to the (real part of) the
quadric.  It is immediate that there is a single point and 
a single time in the homotopy where this occurs (in fact, 
the vector field vanishes at this place and time),
and it is easy to check that this uniqueness survives small
perturbations (see figure \ref{fig:intersection count}). 
$\Cox$

Having this geometric understanding of $\overline{\sing_\qt},
t_1, t_2, t_3, t_4$ and $\delta$, we may now compute the
integral.  We find that
\begin{eqnarray*}
\int_{\thom} \dblres \, \cdot \, 
   \left ( \frac{dt}{t - t_3} - \frac{dt}{t - t_4} 
   \right ) & = & \dblres \, . \, 
   \log \frac{(t_2 - t_3)(t_1 - t_4)}{(t_1 - t_3)(t_2 - t_4)} \\[2ex]
& = & i \, (\alpha + \beta) \, \dblres \, .
\end{eqnarray*}
\begin{figure}[ht!]
\hspace{1.8in} \includegraphics[scale=0.30]{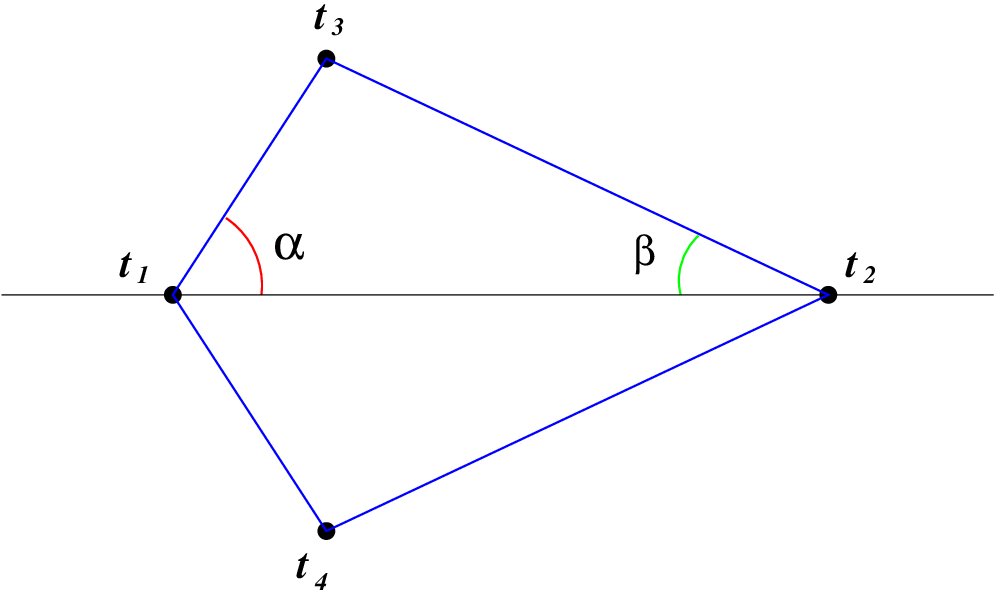}
\caption{the logarithm as an arctangent}
\label{fig:quad}
\end{figure}
Here, the fact that $t_1, t_2$ are complex conjugates while
$t_3, t_4$ are real implies that the numerator and denominator 
are complex conjugates and the logarithm is purely imaginary,
being in fact $2i$ times the arctangent of the sum $\alpha
+ \beta$ of the angles shown in figure~\ref{fig:quad}.
The logarithm is therefore given by twice the arctangent 
of the ratio of imaginary to real parts in the numerator.  

Whenever $t_1, t_2$ satisfy a quadratic equation $t^2 + a t + b$ 
with real coefficients while $t_3 , t_4$ satisfy a quadratic equation 
$t^2 + a' t + b' = 0$, also with real coefficients,  
then simple algebra shows the cross ratio to be given by 
\begin{equation}\label{eq:ratio}
\frac {(t_1-t_4)(t_2-t_3)} {(t_1-t_3)(t_2-t_4)}
   = \frac{b+b'-aa'/2+i\sqrt{a^2-4b}\sqrt{4b'-a'^2}}
   {b+b'-aa'/2-i\sqrt{a^2-4b}\sqrt{4b'-a'^2}} \, .
\end{equation}
The ratio of the imaginary to real parts simplifies considerably, 
so we obtain the equivalent expressions 
\begin{equation} \label{eq:arctan}
2i \arctan \frac{\sqrt{a^2-4b}\sqrt{4b'-a'^2}} {b+b'-aa'/2} \, .
\end{equation}
Here, we recall the definition of the range of the arctangent function
in Theorem~\ref{th:cone and plane}, namely $0 \leq \arctan x < \pi$.

Let $L = \ell_1 x + \ell_2 y + \ell_3 z$ describe $\hht$ in our
coordinate system.  Then $t_3 , t_4$ solve $\qt = L = 0$.  The
minimal polynomial for $t_3$ and $t_4$ (produced, for
example, in Maple as an elimination polynomial for the ideal
$\langle x-(t^2+1) , y-2t , z-(t^2-1) , x^2-y^2-z^2, 
\ell_1 x + \ell_2 y , \ell_3 z$) is given by
$$t^2 + \frac{2 \ell_2}{\ell_1 + \ell_3} t 
   + \frac{\ell_1 - \ell_3}{\ell_1 + \ell_3} \, .$$
Similarly, let $\rr = r_1 x + r_2 y + r_3 z$, giving
the minimal polynomial for $t_1$ and $t_2$ as
$$t^2 + \frac{2 r_2}{r_1 + r_3} t 
   + \frac{r_1 - r_3}{r_1 + r_3} \, .$$
Plugging in $a =  2 r_2 / (r_3 + r_1)$, $b = (r_1 - r_3)/(r_1 + r_3)$,
$a' =  2 \ell_2 / (\ell_3 + \ell_1)$ and
$b' = (\ell_1 - \ell_3)/(\ell_1 + \ell_3)$ to~\eqref{eq:arctan} now gives
$$ \int_{\thom} \omega_L = 2 i \, \dblres \, 
   \arctan \left ( \frac{\sqrt{r_1^2 - r_2^2 - r_3^2} 
   \sqrt{-\ell_1^2 + \ell_2^2 + \ell_3^2}}
   {r_1 \ell_1 - r_2 \ell_2 - r_3 \ell_3} \right )$$
where the two quantities under the radical signs are both positive.
Writing the right-hand side as a combination of coordinate-free
quantities, this becomes 
\begin{equation} \label{eq:invariant}
\int_{\thom} \omega_L = 2 i \, \dblres \, 
   \arctan \frac{\sqrt{\qt^* (\rr , \rr)} \sqrt{-\qt^* (\hht , \hht)}}
   {\qt^*(\rr , \hht)} \, .
\end{equation}
Combining this with the result of Lemma~\ref{lem:thom} shows that
$$\IFT \left ( \frac{1}{\qt \cdot \hht} \right ) = \frac{-i}{\pi}
   \dblres 
   \arctan \frac{\sqrt{\qt^* (\rr , \rr)} \sqrt{-\qt^* (\hht , \hht)}}
   {\qt^*(\rr , \hht)} \, ,$$
and checking this against~\ref{eq:suffices 1} proves 
Theorem~\ref{th:cone and plane} in the case where $\rr$
is inside the dual cone.

\subsubsection*{Case 2: $t_1$ and $t_2$ are real}

In this case, $\rr$ is in the pointy region of the teardrop.
The result~(6.23) from~\cite{ABG} tells us that the Leray
cycle, which is by definition a relative cycle, is an
absolute cycle.  It follows that the intersection class
$\delta$ is an absolute cycle in the twice punctured sphere.
Thus $t_1 = t_2$ when the roots are real, and the cycle
$\delta$ is represented by a circle.  Geometrically, in
the previous case ($t_1, t_2$ complex), as $\rr$ approaches 
the boundary of the dual to the conic, the points $t_1$
and $t_2$ converge to a single point on the equator and 
the arc $\delta$ closes up into a circle.  Provided this
point of convergence is not one of the poles, $t_3$ or $t_4$,
the integral will approach a well defined limit, which is
the integral over the absolute intersection cycle.  By
continuity, the homology class of this absolute cycle cannot 
vary as the limit point varies over the common boundary
of the two regions of the teardrop, nor can this class vary
as $\rr$ varies over the pointy region of the teardrop.

To summarize, there is a constant $c$ such that for all $\rr$ 
in the pointy region, the integral is equal to $c$.  This
is also the limiting value if $\rr$ approaches any point
of the common boundary from inside the other region, and
thus coincides with the limit of the quantity in the previous
case, as $\rr$ approaches any ray in the common boundary;
in the limit the arctangent is $\pi$ and we obtain simply
$P(\ZZ) \dblres \ZZ^{-\rr}$.
$\Cox$

\begin{unremark}
If $\rr$ crosses out of the dual conic at a point $\alpha$
not on the boundary of the pointy region, it exits the 
normal cone.  We know the integral must become zero in this case.
Geometrically, this corresponds to $t_1$ and $t_2$ coming together
in a cycle homologous to zero.  There is a discontinuity if 
$\alpha$ is one of the two projective points of tangency in 
figure~\ref{fig:disk-teardrop}.  Near these two points, $t_1$ and 
$t_2$ approach $t_3$ and $t_4$ respectively.  Crossing out
out of the dual cone on one side or the other will cause
$\delta$ to close up to a null or non-null cycle, in the former
case the integral is zero; the difference between the two integrals
is the residue at the pole $t_3$ or $t_4$.  
\end{unremark}

\setcounter{equation}{0}
\section{Further questions} \label{sec:further}

\begin{enumerate}
\item
It would be nice to remove the integrability hypothesis $2s > d+1$
from Theorem~\ref{th:general}.  Doing so would necessitate a 
specification of which antiderivative $\int_{\lambda_*}^0 
g(\rr  +\lambda L_j) \, d\lambda$ is meant when the integral 
in question is not convergent.  There are cases when the is
one ``obvious'' interpretation of this as a closed form function,
but proving this to be correct requires better understanding of
the generalized function partially identified as $g$ on the
Paley-Wiener cone.  

\item When $d=3$ but $Q$ has an isolated singularity of degree 
greater than~2, the techniques of Section~\ref{ss:cone plane} are still
applicable through Lemma~\ref{lem:thom}.  The representation in
Lemma~\ref{lem:t} must be replaced by one with four poles, and
the intersection class in Lemma~\ref{lem:cycle} correspondingly
specified.  Work is in progress on details of this computation
and its application to the Fortress tiling ensemble.

\item In principle, generalized Fourier transform theory should
give us some information on asymptotics in scaling windows near
the boundary or in obscured directions.  An additional complication
is that, because projective homotopies do not exist giving exponential
decay of $\exp (- \rr \cdot \xx)$ for these $\rr$, one must verify 
the existence of chains on which $\exp (- \rr \cdot \xx)$ decays
sufficiently rapidly to justify the exchanges of limits in 
Lemma~\ref{lem:equal} and elsewhere.  This, together with the
increased complexity of Fourier integrals with varying parameters, 
has kept us thus far from obtaining limit theorems near these
boundaries.  This is perhaps the most broad and challenging 
open problem pertaining to the results in this paper.
\end{enumerate}

\section*{Acknowledgement}

Thanks to J. Borcea for providing a self-contained proof of 
Proposition~\ref{pr:relinearization}.

\bibliographystyle{alpha}
\bibliography{RP}

\end{document}